\documentclass[11pt,final]{iuthesis}
\usepackage{latexsym,amsfonts,amssymb,amsmath}
\newtheorem{theorem}{Theorem}[chapter]
\newtheorem{definition}[theorem]{Definition}
\newtheorem{proposition}[theorem]{Proposition}
\newtheorem{lemma}[theorem]{Lemma}
\newtheorem{remark}[theorem]{Remark}
\newtheorem{corollary}[theorem]{Corollary}

\def\Utimes{\cup\kern-.6383em\lower-.7ex\hbox{$_\times$}}

\def\utimes{\cup\kern-.86em\lower-.7ex\hbox{$_\times$}}

\begin{document}

\title{complex analysis methods in noncommutative probability}
\author{Serban Teodor Belinschi}
\advisor{Hari Bercovici}
\secondreader{Scott W. Brown}
\thirdreader{Ayelet Lindenstrauss}
\fourthreader{Zhenghan Wang}

\department{Mathematics}
\defensedate{17 May 2005}
\submitdate{July 2005}
\copyrightyear{2005}



\begin{dedication}
LV
\end{dedication}

\begin{acknowledgements}
First of all, my deepest thanks go to Professor Hari Bercovici. Without his 
teaching 
and constant, careful, and patient guidance, nothing of my PhD, including this 
thesis, could have been accomplished. I am deeply grateful also to Professor 
Zhenghan Wang. I have learned a lot from him, and his help and support 
(including financial, through an RAship) have been determinant for the 
completion of my work towards a PhD. I regret very much I wasn't able to
do more!

I thank Professors Ayelet Lindenstrauss and Scott W. Brown for their active
participation in my PhD committee. 

I have learned very much from many faculty members of the IU Department of 
Mathematics, either during courses, or in private discussions. Among them, 
I would like to mention Professors Steen Andersson, Eric Bedford, Hari 
Bercovici, Michael Larsen, Ayelet Lindenstrauss, Charles Livingston, Russel
Lyons, Kevin Pilgrim, and Zhenghan Wang. I also thank Misty Cummings for
all the help she gave me with the so-many bureaucratic problems.

Special thanks are due to Professor \c{S}erban Str\v{a}til\v{a}, who first 
introduced me to the realm of advanced mathematics, and whose help, trust, and
guidance made possible for me to become a student at IU. I also thank 
Professor George Dinc\v{a} for his help and teaching.

Finally, I thank my parents and my grandparents, to whom this thesis is 
dedicated, and all my friends.

\end{acknowledgements}

\begin{abstract}
We study convolutions that arise from noncommutative probability theory. 
In the case of free convolutions, we prove that the absolutely continuous part,
with respect to the Lebesgue measure, of the free convolution of two 
probability measures is always nonzero, and has a locally analytic density. 
Under slightly less general hypotheses, we show that the singular continuous 
part of the free additive convolution of two probability measures is zero. We 
also show that any probability measure belongs to a partially defined 
one-parameter free convolution semigroup. In this context, we find a connection
between free and boolean infinite divisibility.
For monotonic convolutions, we prove that any infinitely divisible probability
measure with respect to monotonic additive or multiplicative convolution
belongs to a one-parameter semigroup with respect to the corresponding
convolution. Our main tools are several subordination and inversion theorems 
for analytic functions defined in the upper half-plane. We prove these theorems
using the theory of Denjoy-Wolff fixed points.

\end{abstract}

\frontmatter
\maketitle
\signaturepage
\copyrightpage
\makededication
\makeack
\makeabstract
\tableofcontents

\mainmatter

\section[introduction]{Introduction}

The idea of noncommutative analogues of classical mathematical notions like 
topology, geometry, measure theory, or probability, derives from operator 
algebras. In the classical (commutative) context, one can usually associate to 
a given space an appropriate set of complex-valued functions defined on that 
space. This set becomes a commutative algebra under the operations of
pointwise addition and multiplication of functions. 
Properties of the original space can be then expressed in terms of the algebra 
associated to that space. In the noncommutative context, one replaces 
the commutative algebra of functions with a noncommutative algebra (usually, 
but not necessarily, of operators). We exemplify with the case of interest for 
us.

To define a probability space, one specifies a triple $(\Omega,\Sigma,P)
$, where $\Omega$ is a set, $\Sigma$ is a $\sigma-$algebra,
and $P$ is a probability measure.
Consider the 
space $L^\infty(\Omega)=L^\infty(\Omega,\Sigma,P)$ of essentially bounded
measurable functions $f\colon\Omega\longrightarrow\mathbb C,$ endowed with the
linear functional 
$$\varphi\colon L^\infty(\Omega)\longrightarrow\mathbb C,\quad
\varphi(f)=\int_\Omega f\,dP.$$
It is known that $L^\infty(\Omega)$ is a $W^*$-algebra.
Let us observe that for any set $A\in\Sigma$, the characteristic function of $
A$, 
$$\chi_A(\omega)=\left\{\begin{array}{lcc}
1 & {\rm if} & \omega\in A\\
0 & {\rm if} & \omega\not\in A
\end{array}\right.
$$
belongs to $L^\infty(\Omega),$ and $\varphi(\chi_A)=P(A).$
Conversely, if $f\in L^\infty(\Omega)$ satisfies $f(\omega)^2=f(\omega)=
\overline{f(\omega)}$ for $P-$almost all $\omega\in\Omega$, then $f^{-1}(\{1\})
\in\Sigma$ and $f=\chi_{f^{-1}(\{1\})}$ almost surely. Moreover, $\varphi(f)=P
(f^{-1}(\{1\})).$

Thus, it appears natural to give the 
definition of a noncommutative probability space as it follows:
\begin{definition}\label{Noncommutativeprobabiltyspace}
A noncommutative probability space is a pair $(\mathcal A,\varphi)$, where
$\mathcal A$ is a unital algebra over $\mathbb C$, and $\varphi$ is a linear
functional on $\mathcal A$ 
such that $\varphi(1)=1.$ An element $x\in\mathcal A$ is called a random 
variable.
The distribution of $x$ is the linear functional $\mu_x$ on 
$\mathbb C[X]$ (the algebra of complex polynomials in one variable), defined
by $\mu_x(P)=\varphi(P(x))$, $P\in\mathbb C[X].$
\end{definition}
For an introduction to the subject, we refer to \cite{VDN}.
The pair $(\mathcal A,\varphi)$ will be called a $W^*$-probability space if,
in addition, $\mathcal A$ is a $W^*$-algebra, and $\varphi$ is positive and
weakly continuous. We shall generally consider the case when $\varphi$ is also
a trace (i.e. $\varphi(xy)=\varphi(yx)$ for all $x,y\in\mathcal A$) in which
case we call $(\mathcal A,\varphi)$ a tracial $W^*$-probability space.
Let us observe that $(L^\infty(\Omega),\varphi)$ is 
a tracial $W^*$-probability space, and the notion of random
variable in the noncommutative context translates in the classical case as 
``bounded random variable''. 
Moreover, if $x=x^*$ is a selfadjoint random variable in a $W^*-$probability 
space, then $\mu_x$ extends to a compactly supported probability measure on
the real line, i.e. there exists a unique probability measure $d\mu_x$ such
that 
$$\int_\mathbb R P(t)d\mu_x(t)=\varphi(P(x))\quad{\rm for\ all}\quad P\in
\mathbb C[X].$$

While the classical notion of independence still has a meaning in the 
noncommutative context (see \cite{VDN} for details), new types of independence
arise naturally in noncommutative probability theory. We record here 
the defintions of the three types of independence that we shall encounter
in this paper.
In the following three definitions, $(\mathcal A,\varphi)$ is a noncommutative 
probability space, and $\mathcal A_1,\mathcal A_2$ are two subalgebras of 
$\mathcal A.$
\begin{definition}\label{free}
The unital algebras $\mathcal A_1,\mathcal A_2$ are called free (or freely 
independent) if $\varphi(a_1a_2\cdots a_n)=0$ whenever $a_j\in\mathcal A_{i_j},
i_j\in\{1,2\},$ $i_1\neq i_2,i_2\neq i_3,\dots,i_{n-1}\neq i_n,$ $n\in\mathbb N
,$ and $\varphi(a_{i_j})=0$ for $1\leq j\leq n.$
\end{definition}
\begin{definition}\label{boolean}
The (usually not unital) algebras $\mathcal A_1,\mathcal A_2$ are said to be 
boolean independent if 
   $$\varphi(a_1a_2\cdots a_n)=\varphi(a_1)\varphi(a_2)\cdots\varphi(a_n)$$
for any $a_j\in\mathcal A_{i_j}$, $i_j\in\{1,2\},$ $i_1\neq i_2,i_2\neq i_3,
\dots,i_{n-1}\neq i_n,$ $n\in\mathbb N,$ $1\leq j\leq n.$
\end{definition}
\begin{definition}\label{monotonic}
The (usually not unital) algebras $\mathcal A_1,\mathcal A_2$ are called
monotonically independent if the following two conditions are satisfies:
\begin{enumerate}
\item[{\rm(1)}] for all $a_1,b_1\in\mathcal A_1$ and $a_2\in\mathcal A_2,$ 
we have $a_1a_2b_1=\varphi(a_2)a_1b_1$;
\item[{\rm(2)}] for all $a_1\in\mathcal A_1$ and $a_2,b_2\in\mathcal A_2,$
we have $\varphi(a_2a_1b_2)=\varphi(a_2)\varphi(a_1)\varphi(b_2),$
$\varphi(a_2a_1)=\varphi(a_2)\varphi(a_1),$ and $\varphi(a_1b_2)=\varphi(a_1)
\varphi(b_2).$ 
\end{enumerate}
\end{definition}

A pair of random variables $a_1,a_2\in\mathcal A$ is said to be independent in 
any of the above senses if the subalgebras $\mathcal A_1,\mathcal A_2$ (unital
in the case of free independence) generated by $a_1,a_2$, are independent.
Definition \ref{free} is due to Voiculescu \cite{VDN}, definition \ref{monotonic}
to Muraki \cite{Muraki}, while the notion of boolean independence was first 
explicitly formulated by Speicher and Woroudi \cite{SpeicherWouroudi}, but its
origins can be traced to Bo\.zejko \cite{Bozejko86}.

These notions of independence make it possible to define several types of
convolutions of probability measures.
It is known from the classical probability theory that any compactly 
supported probability measure $\mu$ on $\mathbb R$ can be realized as the 
distribution of a selfadjoint random variable $X_\mu$ 
belonging to some $W^*$-probability space.  Roughly speaking, to define 
a convolution of two compactly supported probability measures, one considers a 
pair of subalgebras of a noncommutative probability space which are
independent in the required sense, and finds in each of the algebras a 
selfadjoint random variable distributed according to each of the two 
probability measures; the distribution of the sum will be equal to the 
convolution of the two probability measures. The following definitions 
make the above description more precise.

\begin{definition}\label{Free}
\ 
\begin{enumerate}
\item[{\rm (a)}]
Let $\mu,\nu$ be two compactly supported Borel probability measures on the real
line. The free additive convolution of $\mu$ and $\nu$,
denoted by $\mu\boxplus\nu$, is defined as the distribution of $x_\mu+x_\nu$,
where $x_\mu$ and $x_\nu$ are selfadjoint random variables belonging to some
tracial $W^*$-probability space $(\mathcal A,\varphi)$, free from each other,
and distributed according to $\mu$ and $\nu$, respectively.
\item[{\rm(b)}] Let $\mu,\nu$ be two compactly supported probability measures 
on the positive half-line $[0,+\infty).$ The free multiplicative convolution
$\mu$ and $\nu$, denoted by $\mu\boxtimes\nu$, is the distribution of 
$x_\mu^{1/2}x_\nu x_\mu^{1/2}$, where  $x_\mu$ and $x_\nu$ are positive random 
variables belonging to some tracial $W^*$-probability space 
$(\mathcal A,\varphi)$, free from each other, and distributed according to 
$\mu$ and $\nu$, respectively. 
\item[{\rm(c)}] Let $\mu,\nu$ be two probability measures
on the unit circle $\mathbb T$ in the complex plane. 
 The free multiplicative convolution of
$\mu$ and $\nu$, denoted by $\mu\boxtimes\nu$, is the distribution of
$x_\mu x_\nu$, where  $x_\mu$ and $x_\nu$ are unitary random
variables belonging to some tracial $W^*$-probability space
$(\mathcal A,\varphi)$, free from each other, and distributed according to
$\mu$ and $\nu$, respectively.
\end{enumerate}
\end{definition}
Definitions (b) and (c) have been extended to measures with noncompact
support in \cite{BercoviciVoiculescuIUMJ}. For more details, we refer to 
\cite{VDN} and \cite{BercoviciVoiculescuIUMJ}.

\begin{definition}\label{Boolean}
\
\begin{enumerate}
\item[{\rm(a)}] Let $\mu,\nu$ be two compactly supported Borel probability measures on the real
line. The boolean additive convolution of $\mu$ and $\nu$,
denoted by $\mu\uplus\nu$, is defined as the distribution of $x_\mu+x_\nu$,
where $x_\mu$ and $x_\nu$ are selfadjoint random variables belonging to some
tracial $W^*$-probability space $(\mathcal A,\varphi)$, boolean independent,
and distributed according to $\mu$ and $\nu$, respectively.
\item[{\rm(b)}] Let $\mu,\nu$ be two probability measures
on the unit circle $\mathbb T$ in the complex plane.
 The boolean multiplicative convolution of
$\mu$ and $\nu$, denoted by $\mu\utimes\nu$, is the distribution of
$x_\mu x_\nu$, where  $x_\mu$ and $x_\nu$ are unitary random
variables belonging to some tracial $W^*$-probability space
$(\mathcal A,\varphi)$,  distributed according to
$\mu$ and $\nu$, respectively, and such that $x_\mu-1$ and $x_\nu-1$
are boolean independent.
\end{enumerate}
\end{definition}
These definitions appear in \cite{SpeicherWouroudi} and \cite{Franz}, 
respectively. The extension of part (a) of the above definition to
measures with noncompact support is done using analytic methods (see
Theorem \ref{booleananalitic} below).

\begin{definition}\label{Monot}
\
\begin{enumerate}
\item[{\rm(a)}] Let $\mu,\nu$ be two compactly supported Borel probability measures on the real
line. The monotonic additive convolution of $\mu$ and $\nu$,
denoted by $\mu\rhd\nu$, is defined as the distribution of $x_\mu+x_\nu$,
where $x_\mu$ and $x_\nu$ are selfadjoint random variables belonging to some
tracial $W^*$-probability space $(\mathcal A,\varphi)$, monotonic independent,
and distributed according to $\mu$ and $\nu$, respectively.
\item[{\rm(b)}] Let $\mu,\nu$ be two compactly supported Borel probability measures on
the positive half-line. The monotonic multiplicative convolution of $\mu$ and $\nu$,
denoted by $\mu\circlearrowright\nu$, is defined as the distribution of 
$x_\mu x_\nu$,
where $x_\mu$ and $x_\nu$ are selfadjoint random variables belonging to some
tracial $W^*$-probability space $(\mathcal A,\varphi)$, distributed according 
to $\mu$ and $\nu$, respectively, such that $x_\mu-1$ and $x_\nu-1$ are 
monotonic independent.
\end{enumerate}
\end{definition}
Observe that, unlike the previous convolutions, monotonic 
convolutions are {\it not} commutative. Monotonic additive convolution
has been introduced by Muraki \cite{Muraki}, while its multiplicative 
version has been defined by Bercovici in \cite{BercoviciMonot}. As in the
case of boolean convolution, one can define monotonic convolutions of
probability measures with noncompact support with analytic methods, as
it will be seen in Theorems \ref{MurakiMonot} and \ref{BercoviciMonot}.

It is not hard to see that the independence condition makes each
of the above seven operations well-defined.
The way to construct independent random variables as in the above definitions
can be found in \cite{VDN} for
the case of free independence, \cite{SpeicherWouroudi} for the case
of boolean independence, and \cite{Muraki} and \cite{BercoviciMonot} for 
monotonic independence.

The present thesis focuses on an analytic perspective on the 
convolutions defined above, which will be described in the main body 
of the thesis. Our main interest will be in the free convolutions. We 
will also point out connections between the three types of independence 
defined above. These connections have not yet been explored to their
fullest extent. 

The thesis has three chapters. The first two chapters have three sections
each, while the much shorter third chapter is not divided into sections.
In the first chapter we focus on regularity properties of free convolutions.
The first section of Chapter 1 introduces a number of definitions and 
theorems from complex analysis that will be used all along the thesis, and
describes the absolutely continuous, singular continuous and atomic parts 
of the free additive convolution of two probability measures, with respect 
to the Lebesgue measure on $\mathbb R.$ Specifically, we show that the 
absolutely continuous part is always nonzero, the density function is locally 
analytic, and, under slightly less general conditions, the singular continuous 
part is zero (Theorem \ref{Reg+}). Our main tool is a subordination result 
for analytic functions (Theorem \ref{Firstsubordinationresult}).
In Sections 2 and 3 we prove similar results for free multiplicative 
convolutions of probability measures supported on the positive half-line and on
the unit circle in the complex plane, respectively (Theorems \ref{RegX1} and 
\ref{RegX2}). Note however that the existence of a singular continuous part 
remains open in these cases.

The second chapter is dedicated to proving the existence of partially
defined semigroups with respect to free convolutions (Theorems \ref{mut+},
\ref{mutX1}, and \ref{mutX2}). Our main tools are three inversion theorems for analytic functions.
In Theorems \ref{RegPtmut+}, \ref{regmutX1},
and \ref{regmuX2}, we describe
the absolutely continuous, singular continuous, and atomic parts of such 
measures. A brief mention of a connection with boolean convolution is also 
made.

In the third chapter we generalize results of Muraki and Bercovici
on monotonic infinite divisibility. Namely, we show that any probability 
measure, not necessarily with compact support, which is infinitely
divisible with respect to monotonic additive, and, respectively, 
multiplicative, convolution belongs to a one-parameter semigroup with respect 
to monotonic convolution.


\chapter[Regularity for free convolutions...]{Regularity properties for
free convolutions of two probability measures}

This chapter is dedicated to describing the absolutely continuous and
atomic part of free convolutions of measures.

In the first section we discuss the behavior of the absolutely
continuous part with respect to the Lebesgue measure of the free additive
convolution of two Borel probability measures supported on the real line, none
of them concentrated in one point.
Namely, we show that the the absolutely continuous part is never zero,
its density function is locally analytic, and continuous everywhere where 
it is finite. 
The results are derived by studying the boundary behavior 
of the subordination functions. As a byproduct, we also obtain a new proof
of Biane's subordination result.

The second and third sections are dedicated to proving similar results for 
the free multiplicative convolution of probability measures on the unit circle
and the positive half-line, respectively. In addition, we also describe the 
atomic part of free multiplicative convolution of two probability measures.
(The analogous result for free additive convolution is already known.)

\section[regularity for free additive convolution...]{Regularity for the free 
additive convolution of probability measures}

For any finite positive measure $\sigma$ on $\mathbb R$, define its Cauchy 
transform
$$G_\sigma(z)=\int_\mathbb R \frac{d\sigma(t)}{z-t},\quad z\in\mathbb C
\setminus\mathbb R.$$
Since $G_\sigma(\overline{z})=\overline{G_\sigma(z)},$ we shall consider from 
now on only the restriction of $G_\sigma$ to the upper half-plane 
$\mathbb C^+=\{z\in\mathbb C\colon \Im z>0\}.$ 

For given $\alpha\geq0,\beta>0,$ let us denote 
$\Gamma_{\alpha,\beta}=\{z\in\mathbb C^+\colon\Im z>\alpha,|\Re z|<\beta\Im z
\}.$
The analytic function $G_\sigma$ satisfies the following two properties:
\begin{trivlist}
\item[{\ \ (i)}] $G_\sigma(\mathbb C^+)\subseteq-\mathbb C^+;$
\item[{\ \ (ii)}] For any $\alpha,\beta>0,$
$$\lim_{z\to\infty,z\in\Gamma_{\alpha,\beta}}zG_\sigma(z)=\sigma(\mathbb R).$$
\end{trivlist}
Remarkably, as the following theorem shows, any function satisfying these two 
properties is the Cauchy transform of some finite positive measure on 
$\mathbb R$ (for proof and details, we refer to \cite{Achieser}).

\begin{theorem}\label{Cauchy}
Let $G\colon\mathbb C^+\longrightarrow-\mathbb C^+$ be an analytic function.
The following statements are equivalent:
\begin{enumerate}
\item There exists a unique positive measure $\sigma$ on $\mathbb R$ such that 
$G=G_\sigma$;
\item  For any $\alpha,\beta>0,$ we have that
$$\lim_{z\to\infty,z\in\Gamma_{\alpha,\beta}}zG(z)$$ exists and is finite.
\item The limit $\lim_{y\to+\infty}iyG(iy)$ exists and is finite.
\end{enumerate}
Moreover, the two limits from 2 and 3 both equal $\sigma(\mathbb R).$
\end{theorem}
Observe also that 
$$-\frac1\pi\Im G_\sigma(x+iy)=\frac1\pi\int_\mathbb R\frac{y}{(x-t)^2+y^2}d\sigma(t),
\quad x\in\mathbb R,y>0,$$
is the Poisson integral of $\sigma.$

It turns out that in many situations it is much easier to deal with the 
reciprocal $F_\sigma=1/G_\sigma$ of the Cauchy transform of the measure $\sigma
.$ The following proposition is an obvious consequence of Theorem \ref{Cauchy}:
\begin{proposition}\label{Cauchyalt}
Let $F\colon\mathbb C^+\longrightarrow\mathbb C^+$ be an analytic self-map of 
the upper half-plane. The following statements are equivalent:
\begin{enumerate}
\item There exists a positive measure $\sigma$ on $\mathbb R$ such that 
$F=1/G_\sigma$;
\item For any $\alpha,\beta>0,$ the limit 
$\lim_{z\to\infty,z\in\Gamma_{\alpha,\beta}}\frac{F(z)}{z}$ exists and belongs
to $(0,+\infty)$;
\item The limit $\lim_{y\to+\infty}\frac{F(iy)}{iy}$ exists and belongs to 
$(0,+\infty)$.
\end{enumerate}
Moreover, both limits form 2. and 3. equal $\sigma(\mathbb R)^{-1}.$
\end{proposition}
In general, analytic self-maps of the upper half-plane can be represented 
uniquely by a triple $(a,b,\rho)$, where $a$ is a real number, $b\in[0,+\infty
),$ and $\rho$ is a positive finite measure on $\mathbb R.$
This representation is called the Nevanlinna representation (see \cite{Achieser}).
\begin{theorem}\label{Nevanlinna}
Let $F\colon\mathbb C^+\longrightarrow\mathbb C^+$ be an analytic function.
Then there exists a triple $(a,b,\rho)$, where $a\in\mathbb R$, $b\geq0,$ and
$\rho$ is a positive finite measure on $\mathbb R$ such that
$$F(z)=a+bz+\int_{\mathbb R}\frac{1+tz}{t-z}d\rho(t),\quad z\in\mathbb C^+.$$
The triple $(a,b,\rho)$ satisfies $a=\Re F(i),$ $b=\lim_{y\to+\infty}\frac{F(iy
)}{iy},$ and $b+\rho(\mathbb R)=\Im F(i)$.
\end{theorem}
The converse of Theorem \ref{Nevanlinna} is obviously true.
\begin{remark}\label{ImF>Imz}
{\rm
An immediate consequence of Proposition \ref{Cauchyalt}  and Theorem \ref{Nevanlinna} 
is that for 
any finite measure $\sigma$ on $\mathbb R$, we have $\Im F_\sigma(z)\geq
\sigma(\mathbb R)^{-1}\Im z$
for all $z\in\mathbb C^+,$ with equality for any value of $z$ if and only if 
$\sigma$ is a point mass. In this case, the measure $\rho$ in the statement of
 Theorem \ref{Nevanlinna} is zero.}
\end{remark}

As observed above, any finite measure $\sigma$ on the real line is uniquely 
determined by its
Cauchy transform. Moreover, regularity properties of $\sigma$ can be deduced
from the behavior of $G_\sigma$, and hence of $F_\sigma$, near the boundary 
of its domain. In the following we shall state several classical theorems 
concerning analytic self-maps of the unit disk $\mathbb D=\{z\in\mathbb C\colon
|z|<1\}$ and their boundary behaviour, i.e. the behavior near points in
to the boundary $\mathbb T=\{z\in\mathbb C\colon|z|=1\}$ of $\mathbb D$. 
Because the upper half-plane is 
conformally equivalent to the unit disc via the rational transformation
$z\mapsto\frac{z-i}{z+i},$ most of these theorems will have obvious
formulations for self-maps of the upper half-plane. 

For a function $f\colon\mathbb C^+\longrightarrow\mathbb C\cup\{\infty\}$, 
and a point $x\in\mathbb R$, we say that the nontangential limit of $f$ at 
$x$ exists if the limit 
$\lim_{z\to x,z\in\Gamma_\alpha(x)}f(z)$ exists for all $\alpha>0$, where 
$\Gamma_\alpha(x)=\{z\in\mathbb C^+\colon|\Re z-x|<\alpha\Im z\}.$
A similar definition holds for functions defined in the unit disc. We shall 
denote nontangential limits by $\sphericalangle\lim_{z\to\alpha}f(z),$ or
$$\lim_{\stackrel{ z\longrightarrow\alpha}{{
\sphericalangle}}}f(z).$$

The following three theorems describe properties of meromorphic functions in 
the unit disc related to their nontangential boundary behavior.

\begin{theorem}\label{Fatou}
Let $f\colon\mathbb D\longrightarrow\mathbb C$ be a bounded analytic function.
Then the set of points $x\in\mathbb T$ at which the nontangential limit of 
$f$ fails to exist is of linear measure zero.
\end{theorem}

\begin{theorem}\label{Lindelof}
 Let $f\colon\mathbb D\longrightarrow\mathbb C\cup\{\infty\}$ be a 
meromorphic
function, and let $e^{i\theta}\in\mathbb T$. Assume that the set 
$(\mathbb C\cup\{\infty\})\setminus f(\mathbb D)$ contains more than three 
points. If there exists a path $\gamma\colon[0,1)\longrightarrow\mathbb D$
such that $\lim_{t\to1}\gamma(t)=e^{i\theta}$ and $\ell=\lim_{t\to1}f(\gamma
(t))$ exists in $\mathbb C\cup\{\infty\}$, then the nontangential limit 
of $f$ at $e^{i\theta}$ exists, and equals $\ell$
\end{theorem}

\begin{theorem}\label{Privalov}
Let $f\colon\mathbb D\longrightarrow\mathbb C$ be an analytic function.
Assume that there exists a set $A$ of nonzero linear measure in $\mathbb T$ 
such that the nontangential limit of $f$ exists at each point of $A$, and
equals zero. Then the function $f(z)=0$ for all $z\in\mathbb D$.
\end{theorem}
Theorem \ref{Fatou} is due to Fatou, and Theorem \ref{Privalov} to Privalov.
Theorem \ref{Lindelof} is an extension of a result by Lindel\"{o}f. For 
proofs, we refer to \cite{CollingwoodL}.

The following lemma will be useful in the proof of the main result
of this section.

\begin{lemma}\label{nontangentialconvergenceatinfinity}
Let $F\colon\mathbb C^+\longrightarrow\mathbb C^+$ be an analytic self-map of
the upper half-plane.
Suppose
$\gamma\colon[0,1)\longrightarrow\mathbb C^+$ is a continuous path such that
$\lim_{s\to1}\gamma(s)=\infty.$ Then the limit
$\lim_{s\to1}\frac{F(\gamma(s))}{\gamma(s)}$, if it exists, cannot be infinite.
\end{lemma}
\begin{proof}
Assume that there exists a path $\gamma$ as in the statement of the lemma
such that $\lim_{s\to1}\frac{F(\gamma(s))}{\gamma(s)}=\infty.$
Ofserve that $F(z)/z$ never equals a negative number. Therefore, the range of
the function
$$\Phi\colon\mathbb D\longrightarrow\mathbb C,\quad\Phi(w)=\frac{F\left(i\frac{
1-w}{1+w}\right)}{i\frac{1-w}{1+w}}$$
omits more than three points from $\mathbb C$.
Since $\sphericalangle\lim_{z\to\infty}F(z)/z$ exists and is finite, Theorem 
\ref{Lindelof} provides a contradiction.

\end{proof}

Consider a domain (i.e. an open connected set) $D\subseteq\mathbb C\cup\{\infty
\}$, and a function $f\colon D\longrightarrow\mathbb C\cup\{\infty\}$.
The cluster set $C(f,x_0)$ of the function $f$
at the point $x_0\in\overline{D}$ is
$$\{z\in\mathbb C\cup\{\infty\}\ |\  \exists
\{z_n\}_{n\in\mathbb N}\subset D\setminus x_0\ {\rm such\ that }
\lim_{n\to\infty}z_n=x_0, \ \lim_{n\to\infty}f(z_n)=z\}.{\rm }$$
The following result is immediate.
\begin{lemma}\label{ClusterConex}
Let $D\subset\mathbb C\cup\{\infty\}$ be a domain and let
$f\colon D\longrightarrow\mathbb C\cup\{\infty\}$ be continuous. If
$D$ is locally connected at $x\in\overline{D}$, then $C(f, x)$ is connected.
\end{lemma}
This result appears in \cite{CollingwoodL}, as Theorem 1.1.

The following theorem of Seidel
describes the behavior of certain analytic functions near the boundary of
their domain of definition.
For proof, we refer to
\cite{CollingwoodL}, Theorem 5.4.

\begin{theorem}\label{Seidel}
Let $f\colon\mathbb D\longrightarrow\mathbb D$ be an analytic function such 
that the radial limit $f(e^{i\theta})=\lim_{r\to1}f(re^{i\theta})$ exists and
has modulus $1$ for almost every
$\theta$ in the interval $(\theta_1,\theta_2)$. If $\theta\in(\theta_1,\theta_2
)$ is such that $f$ does not extend analytically through $e^{i\theta},$ then
$C(f,e^{i\theta})=\overline{\mathbb D}.$
\end{theorem}

This theorem can be applied to self-maps of the upper half-plane $\mathbb C^+$,
via a conformal mapping, but in that case one must consider meromorphic,
instead of analytic, extensions.

A second result refering to the behavior of $C(f,x)$ for bounded analytic 
functions $f$ is the following theorem of Carath\'{e}odory. (This result appears in 
\cite{CollingwoodL}, Theorem 5.5.)

\begin{theorem}\label{Caratheodory}
Let $f\colon\mathbb D\longrightarrow\mathbb C$ be a bounded analytic function. 
Assume that for almost every $\theta\in(\theta_1,\theta_2)$
the radial limit $f(e^{i\theta})$ belongs to a set $W$ in the plane. Then,
for every $\theta\in(\theta_1,\theta_2)$
the cluster set $C(f, e^{i\theta})$ is contained in the
closed convex hull of $W$.
\end{theorem}

The previous two theorems allow us to prove the following
\begin{proposition}\label{SeidelCaratheodory}
\ 
\begin{enumerate}
\item[{\rm (a)}]
Let $f$ be an analytic self-map of $\mathbb D$ such that $|\lim_{r\to1}f(re^{i
\theta})|=1 $ for almost every $\theta\in(\theta_1,\theta_2)$.
Suppose that there is a point $\theta_0\in(\theta_1,\theta_2)$ such that the
function $f$ cannot be continued analytically through $\theta_0.$
Then for any $t_1<t_2$ there is a set $E\subset(\theta_1,\theta_2)$
of nonzero Lebesgue measure such that
$\lim_{r\to1}f(re^{i\theta})$ exists for all $\theta\in E$, and
the set ${\{\lim_{r\to1}f(re^{i\theta})\colon \theta\in E\}}$ is dense in
the arc $A=\{e^{it}\colon t_1<t<t_2\}.$
\item[{\rm (b)}] Let $f$ be an analytic self-map of $\mathbb C^+$ such that
$\lim_{y\to 0}f(x+iy)$ exists and belongs to $\mathbb R$ for almost every
$x\in (a,b).$ Suppose that
$x_0\in(a,b)$ is such that $f$ cannot be continued meromorphically through $x_0
$.
Then for any $c<d$ there is a set
$E\subseteq(a,b)$ of nonzero Lebesgue measure such that
$\lim_{y\to0}f(x+iy)$ exists for all points $x\in E,$ and the set
${\{\lim_{y\to0}f(x+iy)\colon x\in E\}}$ is dense in the interval 
$(c,d).$
\end{enumerate}
\end{proposition}
\begin{proof}
Let $f$ and $\theta_0$ be as in the hypothesis of (a).
Using Theorems \ref{Seidel} and \ref{Caratheodory}, we conclude that,
on the one hand, $C(f,e^{i\theta_0})=\overline{\mathbb D},$ and on the other, 
that $C(f,e^{i\theta_0})$
equals the closure of the convex hull of the set
$${\{\lim_{r\to1}f(re^{i\theta})\colon \theta\in G\}}$$ for any set $G\subset
(\theta_1,\theta_2)$ with the property that $\lim_{r\to1}f(re^{i\theta
})$ exists for all $\theta\in G$, and $(\theta_1,\theta_2)\setminus G$
has zero linear measure. This implies that the set
$$\{\lim_{r\to1}f(re^{i\theta})\colon\theta_1<\theta<\theta_2,\  \lim_{r\to 1}
f(re^{i\theta})\ {\rm exists\ and\ belongs\ to}\  A\}$$
is dense in $A$. It remains to
prove that the set $E$ of those $\theta\in(\theta_1,\theta_2)$ such that
$\lim_{r\to1}f(re^{i\theta})$ exists and belongs to $A$ has nonzero linear
measure. If this were not true, then,
according to Theorem \ref{Caratheodory}, we could replace $G$ by $G\setminus E$ in
the previous argument, and obtain a contradiction. This proves (a).

Part (b) follows directly from (a), by using the conformal automorphism of the
extended complex plane
$z\mapsto\frac{z-i}{z+i}$ and its inverse.
\end{proof}

We will next focus on boundary behaviour of derivatives of analytic
self-maps of the unit disk and of the upper half-plane. These results are 
described in detail
by Shapiro \cite{Shapiro}; see also Exercises 6 and 7 in Chapter I of Garnett's
book \cite{Garnett}. 

\begin{theorem}\label{JuliaCaratheodory}
Let $f\colon\mathbb D\longrightarrow\mathbb D$ be an analytic function, and 
let $w\in\mathbb T$. The following statements are equivalent:
\begin{enumerate}
\item We have
 $$\liminf_{z\to w}\frac{|f(z)|-1}{|z|-1}<\infty;$$
\item  There exists a number $\zeta\in\overline{\mathbb D}$ such that
$$\lim_{\stackrel{ z\longrightarrow w}{{\sphericalangle}}}f(z)=\zeta,$$
and the limit
\begin{equation}
\ell=
\lim_{\stackrel{ z\longrightarrow w}{{\sphericalangle}}}\frac{f(z)-\zeta}{z-w}
\label{eq1.1}
\end{equation}
exists and belong to $(0,+\infty)$.
\end{enumerate}
Moreover, if the equivalent conditions above are satisfied, the limit
$\sphericalangle\lim_{z\to w}f'(z)$ exists, and the 
following equality holds:
 $$\ell=\frac{w}{\zeta}\lim_{\stackrel{ z\longrightarrow w}{{\sphericalangle}}}
f'(z)=\liminf_{z\to w}\frac{|f(z)|-1}{|z|-1}.$$

If $$\liminf_{z\to w}\frac{|f(z)|-1}{|z|-1}=\infty$$ and 
$$\lim_{\stackrel{ z\longrightarrow w}{{\sphericalangle}}}f(z)=\zeta,$$ then
the limit in equation \ref{eq1.1} exists and equal infinity.
\end{theorem}
The number $\ell$ from the above theorem is called the Julia-Carath\'{e}odory
derivative of $f$ at $w$.
The formulation of the preceding theorem for self-maps of the upper half-plane 
is similar; however, for the point at infinity, the usual formula
of the Julia-Carath\'{e}odory derivative must be replaced by the limit
$$\lim_{\stackrel{z\longrightarrow\infty}{\sphericalangle}}\frac{z}{f(z)}.$$

Consider now an analytic function $f\colon\mathbb D\longrightarrow\overline{
\mathbb D}.$ A point $w\in\overline{\mathbb D}$ is called a Denjoy-Wolff
point for $f$ if one of the following two conditions is satisfied:
\begin{enumerate}
\item[(1)] $|w|<1$ and $f(w)=w$;
\item[(2)] $|w|=1,$ $\sphericalangle\lim_{z\to w}f(z)=w, $ and 
$$\lim_{\stackrel{z\longrightarrow w}{\sphericalangle}}\frac{f(z)-w}{z-w}\leq1.
$$
\end{enumerate}
The following result is due to Denjoy and Wolff.
\begin{theorem}\label{DenjoyWolff}
Any analytic function $f\colon\mathbb D\longrightarrow\overline{\mathbb D}$ has
a Denjoy-Wolff point. If $f$ has more than one such point, then $f(z)=z$ for all
$z$ in the unit disk. If $z\in\mathbb D$ is a Denjoy-Wolff point for $f$, then
$|f'(z)|\leq1;$ equality occurs only when $f$ is a conformal automorphism of
the unit disc.
\end{theorem}

The Denjoy-Wolff point of a function $f$ is characterized also by the fact that
it is the uniform limit on compact subsets of the iterates $f^{\circ n}=
\underbrace{f\circ f\circ\cdots\circ f}_{n\ {\rm times}}$ 
 of $f$. We state the 
following theorem for the sake of completeness (for the original statements,
see \cite{Denjoy} and \cite{Wolff}):

\begin{theorem}
Let $f\colon\mathbb D\longrightarrow\overline{\mathbb D}$ be an analytic
function. If $f$ is not a conformal automorphism of $\mathbb D$, then the 
functions $f^{\circ n}$ converge uniformly on compact subsets of $\mathbb D$ to
the Denjoy-Wolff point of $f$.
\end{theorem}

We now prove a general subordination result for analytic functions.

\begin{theorem}\label{Firstsubordinationresult}
Let $F_j\colon\mathbb C^+\longrightarrow\mathbb C^+$, $j\in\{1,2\}$ be two
analytic functions, which are not conformal self-maps of
$\mathbb C^+,$
satisfying $\lim_{y\to+\infty}\frac{F_j(iy)}{iy}=1,$ $j\in\{1,2\}$,
and let $f\colon(\mathbb C^+\cup\mathbb R)\times\mathbb C^+\longrightarrow
\mathbb C^+$ be defined by $f(z,w)=F_2(F_1(w)-w+z)-(F_1(w)-w)$. Then
there
exists a unique analytic self-map $\omega$ of $\mathbb C^+$ such that
$f(z,\omega(z))=\omega(z)$ for all $z\in\mathbb C^+.$ Moreover,
$\lim_{y\to+\infty}\frac{\omega(iy)}{iy}=1.$
\end{theorem}

\begin{proof}
Let us first observe that $f$ is well defined. Indeed, since $F_j$ is not
a conformal automorphism of $\mathbb C^+$, by Remark \ref{ImF>Imz} we have
 $\Im F_j(w)>\Im w$ for all $w\in
\mathbb C^+,$ $j=1,2$. Thus, for any $z$ with $\Im z\geq0,$ $\Im(F_1(w)-w+z)>0
,$ and $\Im F_2(F_1(w)-w+z)>\Im(F_1(w)-w+z)\geq\Im(F_1(w)-w).$
We conclude that the formula defining $f$ makes sense, and $f$ takes values in 
the upper half-plane.

We shall now prove that for any fixed $z\in\mathbb C^+$,
$$\lim_{v\to+\infty}\frac{F_2(F_1(iv)-iv+z)}{iv}\not\in[1,+\infty].$$
Assume to the contrary that there exists an $\ell\in[1,+\infty]$ such that
$$\lim_{v\to+\infty}\frac{F_2(F_1(iv)-iv+z)}{iv}=\ell.$$
Define $\theta\colon[0,1)\longrightarrow\mathbb C^+$, $\theta(t)=\frac{i}{
1-t}.$ Then
\begin{eqnarray*}
\ell & = & \lim_{v\to+\infty}\frac{F_2(F_1(iv)-iv+z)}{iv}\\
& = & \lim_{t\to1}\frac{F_2(F_1(\theta(t))-\theta(t)+z)}{F_1(\theta(t))-
\theta(t)+z}\cdot\frac{F_1(\theta(t))-
\theta(t)+z}{\theta(t)}
\end{eqnarray*}
Since $\lim_{t\to1}\frac{F_1(\theta(t))-\theta(t)+z}{\theta(t)}=0$, we conclude
that
$$\lim_{t\to1}\frac{F_2(F_1(\theta(t))-\theta(t)+z)}{F_1(\theta(t))-
\theta(t)+z}=\infty.$$
This is possible only if $\lim_{t\to1}[{F_1(\theta(t))-\theta(t)+z}]=\infty.$
The path $\gamma(t)=F_1(\theta(t))-\theta(t)+z,$ $t\in[0,1)$, 
satisfies the conditions of Lemma \ref{nontangentialconvergenceatinfinity}, 
and we have 
$\lim_{t\to1}\frac{F_2(\gamma(t))}{\gamma(t)}=\infty,$ and
$\lim_{v\to+\infty}\frac{F_2(iv)}{iv}=1.$ This contradiction 
implies that
$$\lim_{v\to+\infty}\frac{F_2(F_1(iv)-iv+z)}{iv}\not\in[1,+\infty].$$

Using this fact, we have
\begin{eqnarray*}
\lim_{v\to+\infty}\frac{f(z,iv)}{iv} & = &
\lim_{v\to+\infty}\frac{F_2(F_1(iv)-iv+z)-(F_1(iv)-iv)}{iv}\\
& = & \lim_{v\to+\infty}\frac{F_2(F_1(iv)-iv+z)}{iv}\\
& \neq & \ell
\end{eqnarray*}
for all $\ell\in[1,+\infty].$
Also, as seen above, $\Im f(z,w)> \Im z>0.$ Thus, by Theorem \ref{DenjoyWolff}, 
we conclude that the Denjoy-Wloff point of the function $f_z(\cdot)=f(z,\cdot)$
belongs to $\mathbb C^++z$. Denote it by $\omega(z)$.

Since, by Theorem \ref{DenjoyWolff}, $|f_z'(\omega(z))|<1$, we can apply the 
implicit
function theorem to conclude that the function $z\mapsto\omega(z)$ is
an analytic self-map of $\mathbb C^+.$ The uniqueness of the Denjoy-Wolff
point of $f_z$ for each $z$ guarantees the uniqueness of $\omega.$ This
proves the existence and uniqueness of $\omega.$

To complete the proof, we need to show that $\lim_{y\to+\infty}\frac{\omega(iy)
}{iy}=1.$ As
noted above, $\omega(z)\in\mathbb C^++z$, so $\Im\omega(z)>\Im z,$ for all $z
\in\mathbb C^+.$ 
By Theorem \ref{Nevanlinna},
$\lim_{y\to+\infty}\frac{\omega(iy)}{iy}\geq1.$
The same argument applied to $(F_2, F_1)$ in place of $(F_1,F_2$) provides
a unique analytic self-map of the upper half-plane $\tilde{\omega}$ such 
that
$$F_1(F_2(\tilde{\omega}(z))-\tilde{\omega}(z)+z)-(F_2(\tilde{\omega}(z))
-\tilde{\omega}(z))=\tilde{\omega}(z),
\quad z\in\mathbb C^+.$$
We claim that the analytic function $g\colon\mathbb C^+\longrightarrow\mathbb
C^+,$ defined by $g(z)=F_1({\omega}(z))-{\omega}(z)+z,$
$z\in\mathbb C^+,$ also satisfies
$$F_1(F_2(g(z))-g(z)+z)-(F_2(g(z))-g(z))=g(z).$$
Indeed,
$$F_2(g(z))=f(z,\omega(z))+g(z)-z=\omega(z)+g(z)-z=F_1(\omega(z)),$$
so that
\begin{eqnarray}
\lefteqn{F_1(F_2(g(z))-g(z)+z)-(F_2(g(z))-g(z))}\nonumber \\
& = & F_1(F_1(\omega(z))-g(z)+z)-(F_1(\omega(z))-g(z))\nonumber \\
& = & F_1(F_1(\omega(z))-(F_1(\omega(z))-\omega(z)+z)+z)\nonumber \\
& & \mbox{} -(F_1(\omega(z))-(F_1(\omega(z))-\omega(z)+z))\nonumber \\
& = & F_1(\omega(z))-\omega(z)+z\nonumber \\
& = & g(z)\nonumber.
\end{eqnarray}
The uniqueness of $\tilde{\omega}$ implies that $g=\tilde{\omega},$ and thus
$\omega(z)+\tilde{\omega}(z)=F_1(\omega(z))+z.$ 
Since $\Im \tilde{\omega}(z)>\Im z$, $z\in\mathbb C^+,$ the limit 
$\ell=\lim_{y\to+\infty}\tilde{\omega}(iy)/iy$ is at least one. In particular, 
$\tilde{\omega}(iy)$ tends to infinity nontangentially as $y\to+\infty.$ 
Therefore
$$\ell+\lim_{y\to+\infty}\frac{\omega(iy)}{iy}=
\lim_{y\to+\infty}\frac{F_2(\tilde{\omega}(iy)}{\tilde{\omega}(iy)}\cdot\frac{
\tilde{\omega}(iy)}{iy}+1=\ell+1.$$
We conclude that
$\lim_{y\to+\infty}\omega(iy)/iy=1,$
as desired.

\end{proof}

Let us observe that, given functions $F_1$ and $F_2$ as in the above theorem,
the pair of analytic functions $(\omega,\tilde{\omega})$ is uniquely
determined by the system of equations
\begin{equation}
F_1(\omega(z))=F_2(\tilde{\omega}(z))
=\omega(z)+\tilde{\omega}(z)-z,\quad z\in\mathbb C^+.
\label{eq1.2}
\end{equation}
The following theorem describes the boundary behavior of 
the function $\omega.$
\begin{theorem}\label{Boundary+}
Let $F_1,F_2,f,$ and $\omega $ be as in Theorem \ref{Firstsubordinationresult}.
Fix $a\in\mathbb R$, 
and define the following two self-maps
of the upper half-plane: $u_j(w)=F_j(w)-w+a,$ $j\in\{1,2\}$. Assume that $u_2
\circ u_1$
is not a conformal automorphism of $\mathbb C^+.$Then:
\begin{enumerate}
\item[{\rm (1)}] If $C(\omega,a)\cap\mathbb C^+\neq\varnothing,$
then the function $\omega$ extends analytically in a neighbourhood of $a$.
\item[{\rm (2)}] Assume that there exist open intervals $I_1$ and $I_2$ in $
\mathbb R$ such that $F_j$ extends analytically through $I_j,$ $j\in\{1,2\}.$
Then the limit $\lim_{z\to a}\omega(z)$
exists in $\mathbb C\cup\{\infty\}.$
\end{enumerate}
\end{theorem}
\begin{proof}
Consider a sequence $\{z_n\}_{n\in\mathbb N}\subset\mathbb C^+$ and a number
$\ell\in\mathbb C^+$ with the proprety
that $\lim_{n\to\infty}z_n=a$ and $\ell=\lim_{n\to\infty}\omega(z_n).$ We 
have
 $$\ell=\lim_{n\to\infty}\omega(z_n)=\lim_{n\to\infty}f(z_n,\omega(z_n))=
f(a,\ell).$$
If $f(a,\cdot)$ is not an automorphism of the upper half-plane, then
we can use the implicit function argument
from the proof of Theorem \ref{Firstsubordinationresult} to conclude that there
exists a neighborhood 
of $a$ in $\mathbb C$ on which the function $\omega$ extends analytically.
On the other hand, if $f_a=f(a,\cdot)$ is a conformal automorphism of 
$\mathbb C^+$, denote by $k(\cdot)$ its inverse. 
The equality $$w=f_a(k(w))=f(a,k(w))=F_2(F_1(k(w))-k(w)+a)-(F_1(k(w))-k(w))$$
can be rewritten as
                  $$u_2(u_1(k(w)))=w,\quad w\in\mathbb C^+.$$
We conclude that $u_1\circ k$ (and hence 
$u_1$) must be injective. 
We apply $u_1\circ k$ to this equality to obtain

                             $$u_1(k(u_2(z)))=z,$$
for all $z$ in the open set
$(u_1\circ k)(\mathbb C^+),$ and by analytic continuation, for all
$z\in\mathbb C^+.$ We conclude that 
$u_1$ is a
conformal automorphism of the upper half-plane. But then $u_1\circ u_2$ is
also a conformal automorphism of $\mathbb C^+,$ contradicting the hypothesis.
This proves (1).

Assume now that the hypothesis in (2) is satisfied and yet $C(\omega,a)$ 
contains more than one point.
By part (1), $C(\omega,a)\subseteq
\mathbb R\cup\{\infty\},$ and by Lemma \ref{ClusterConex}, either $C(\omega,a)
\setminus\{
\infty\}$ is a closed interval in $\mathbb R$ (possibly all of $\mathbb R$), or
$\mathbb R\setminus C(\omega,a)$ is an open interval in $\mathbb R.$

We claim that for any $c$ in $C(\omega,a)\setminus\{\infty\},$ with the 
possible exception of three points,
there exists a sequence $\{z_n^{(c)}\}_{n\in\mathbb N}$ converging to $a$
such that
$\lim_{n\to\infty}\omega(z_n^{(c)})=c,$ and $\Re\omega(z_n^{(c)})=c$ for all 
$n$. 
Let $\{c_n\}_{n\in\mathbb N}$ be a dense sequence in $C(\omega,a)$, and
consider $z_n\in\mathbb C^+,$ such that $|z_n-a|<1/n,$ and $|\omega(z_n)-c_n|<
1/n,$ $n\in\mathbb N.$ We 
define a path $\gamma\colon[0,1]\longrightarrow\mathbb C^+\cup\{a\}$ 
such that $\gamma(1-1/n)=z_n,$ $\gamma(1)=a,$ and $\gamma$ is linear on the 
intervals $\left[1-\frac1n,1-\frac{1}{(n+1)}\right],$ $n\in\mathbb N.$ It will 
suffice to 
show that there exists at most one 
point $c$ in the interior of $C(\omega,a)$ such
that $\omega(\gamma([0,1)))
\cap\{c+it\colon t\in[0,\varepsilon)\}=\varnothing$  for some $\varepsilon>0$.
Indeed, assume to the contrary that $c<c'$ are two such points.
The set
$$K=\{c+it\colon t\in(0,1)\}\cup\{c'+it\colon t\in(0,1)\}\cup\{s+i,c\leq s
\leq c'\}$$
separates $\mathbb C^+$ into two components, and the path $\omega(\gamma(t))$
contains infinitely many points in either of the components, hence it crosses
$K$ infinitely many times. By our assumption, crossings cannot be close to
$c$ or $c'$, and this implies the existence of a point in $C(\omega,a)\cap K
\subset\mathbb C^+,$ contradicting the hypothesis that $C(\omega,a)\cap
\mathbb C^+=\varnothing.$ This proves our claim.

By Fatou's theorem (Theorem \ref{Fatou}), the limit $\lim_{n\to\infty}
F_1(\omega(
z_n^{(c)}))$ exists for almost all $c\in C(\omega,a).$ Denote it by $F_1(c).$
We shall prove that for every $c\in C(\omega,a),$ with at most three exceptions, $F_1(c)\not\in\mathbb C^+.$ Indeed, 
suppose that $F_1(c)\in\mathbb C^+$ for some $c\in C(\omega,a),$ and assume 
that $\omega(z_n^{(c)})$ converges nontangentially to $c$. Then,
using the relation $f(z,\omega(z))=\omega(z),$ $z\in\mathbb C^+,$  Remark 
\ref{ImF>Imz}, and the fact that $F_2$ is not a conformal self-map of 
$\mathbb C^+,$ we obtain
\begin{eqnarray*}
\Im F_1(c) & = & \lim_{n\to\infty}\Im F_1(\omega(z_n^{(c)}))\\
& = & \lim_{n\to\infty}\Im F_2(F_1(\omega(z_n^{(c)}))-\omega(z_{n}^{(c)})+
z_{n}^{(c)})\\
& = & \Im F_2(F_1(c)-c+a)\\
& > & \Im (F_1(c)-c+a)\\
& = & \Im F_1(c),
\end{eqnarray*} 
which is a contradiction.

Assume that $c_0\in{\rm Int}C(\omega,a)$ is a point where $F_1$ does not
continue meromorphically. Proposition \ref{SeidelCaratheodory} (b) shows that
the set 
$$E=\{c\in C(\omega,a)\colon a+F_1(c)-c\in I_2\}$$ 
has nonzero Lebesgue measure.
In particular, for all $c\in E$,
\begin{eqnarray*}
F_1(c) & = & \lim_{n\to\infty}F_1(\omega(z_n^{(c)}))\\
& = & \lim_{n\to\infty}F_2(F_1(\omega(z_n^{(c)}))-\omega(z_n^{(c)})+z_n^{(c)}
)\\
& = & F_2(F_1(c)-c+a),
\end{eqnarray*}
where we used the analiticity of $F_2$ on $I_2.$
Privalov's theorem implies that $F_1(z)=
F_2(F_1(z)-z+a)$ for all $z\in\mathbb C^+.$
Rewriting this equality gives
                     $$F_2(F_1(z)-z+a)-(F_1(z)-z+a)+a=z,$$
or, equivalently,
$$u_2(u_1(z))=z,\quad z\in\mathbb C^+,$$
contradicting the hypothesis.

We conclude that $F_1$ extends meromorphically through the whole interval 
${\rm Int}C(\omega,a).$ Let, as in the proof of 
Theorem \ref{Firstsubordinationresult}, 
$\tilde{\omega}(z)=F_1(\omega(z))-\omega(z)+z,$ $z\in\mathbb C^+.$ We shall 
argue that 
the set $C(\tilde{\omega},a)\subseteq\mathbb R\cup\{\infty\}$ must be also 
infinite.
Suppose this were not the case. Then for any $c\in{\rm Int}C(\omega,a),$ 
$$\lim_{z\to a}\tilde{\omega}(z)=\lim_{n\to\infty}F_1(\omega(z_n^{(c)}))-
\omega(z_n^{(c)})+z_n^{(c)}=F_1(c)-c+a,$$
so that, by analytic continuation, $F_1(z)=z-a+\lim_{z\to a}\tilde{\omega}(z)$
for all $z\in\mathbb C^+.$
This contradicts the fact that $F_1$ is not a conformal automorphism of $
\mathbb C^+.$ So $C(\tilde{\omega},a)$ must be an infinite set. 
As before,
$F_2$ must continue meromorphically through all of 
${\rm Int}C(\tilde{\omega},a)$.

By passing to a subsequence, 
we may assume
that $\lim_{n\to\infty}\tilde{\omega}(z_n^{(c)})$ exists for every $c\in{\rm In
t}C(\omega,a)$, with at most one exception.
Suppose there were a point $d\in C(\tilde{\omega},a)$ and a set $V_{d}
\subset{\rm Int}C(\omega,a)$ of nonzero Lebesgue measure
such that $\lim_{n\to\infty}\tilde{\omega}(z_n^{(c)})=d$ for all $c\in V_{d}.$
Taking limit as $n\to\infty$ in the equality
$$F_1(\omega(z_n^{(c)}))+z_n^{(c)}=
\omega(z_n^{(c)})+\tilde{\omega}(z_n^{(c)})$$
gives
$$F_1(c)+a=c+d,\quad{\rm for\ all}\  c\in V_{d}.$$
Applying again Privalov's theorem, we obtain that $F_1(z)=z-(a-d)$ for all
$z\in\mathbb C^+.$ This contradicts the fact that $F_1$ is not a conformal 
automorphism of $\mathbb C^+.$
Thus, there exists a set $E\subseteq{\rm Int}C(\omega,a)$ of positive Lebesgue 
measure such that $\{\tilde{c}=\lim_{n\to\infty}\tilde{\omega}(z_n^{(c)})
\colon c\in E\}\subseteq{\rm Int}C(\tilde{\omega},a).$
Then, since $F_2$ extends analytically through ${\rm Int}C(\tilde{\omega},a)$, 
by (\ref{eq1.2}) we conclude that
\begin{eqnarray*}
F_1(c) & = &
\lim_{n\to\infty}F_1(\omega(z_n^{(c)}))\\
& = & \lim_{n\to\infty}F_2(\tilde{\omega}(z_n^{(c)}))\\
& = & \lim_{n\to\infty}F_2\left({F_1(\omega
(z_n^{(c)}))}-\omega(z_n^{(c)})+z_n^{(c)}\right)\\
& = & F_2(a+F_1(c)-c)
\end{eqnarray*}
for all $c\in E$. Privalov's theorem implies that $F_1(z)=F_2(a+F_1(z)-
z)$ for all $z\in\mathbb C^+.$ As we have proved already, this implies that 
$u_1$ and $u_2$ are conformal automorphisms of the upper half-plane,
contradicting the hypothesis.
This concludes the proof.
\end{proof}

In the following, we prove subordination and regularity results for free 
additive convolutions of probability measures.
We shall start by stating (without proof) some well-known fundamental facts.

Voiculescu \cite{Voiculescu1} has showed that 
there exists an analogue of the logarithm of the Fourier transform for free
additive convolution, called the $R-$transform. The $R-$transform is defined in
terms of the inverse ``near infinity'' of the reciprocal of the Cauchy
transform. For the proofs of the following results we refer to
\cite{BercoviciVoiculescuIUMJ} and \cite{VDN}.
\begin{proposition}\label{inversion+}
Let $\mu$ be a probability on $\mathbb R$. There exists a nonempty domain
$\Omega$ in $\mathbb C^+$ of the form $\Omega=\cup_{\alpha>0}\Gamma_{\alpha,
\beta_\alpha}$ such that $F_\mu$ has a right inverse with respect to
composition $F_\mu^{-1}$ defined on $\Omega$. In addition, we have
$\Im F_\mu^{-1}(z)\leq\Im z$ and $$\lim_{z\to\infty,z\in
\Gamma_{\alpha,\beta}}\frac{F_\mu^{-1}(z)}{z}=1$$ for every $\alpha,\beta>0.$
\end{proposition}
(We recall that $\Gamma_{\alpha,\beta}=\{z\in\mathbb C^+\colon\Im z
>\beta,|\Re z|<\alpha\Im z\}.$)
Let $\varphi_\mu(z)=F_\mu^{-1}(z)-z$, $z\in\Omega.$ The basic property of the
function $\varphi_\mu$ is described in the following theorem of Voiculescu:
\begin{theorem}\label{R-transform}
Let $\mu,\nu$ be two probability measures supported on the real line.
Then $\varphi_{\mu\boxplus\nu}(z)=\varphi_\mu(z)+\varphi_\nu(z)$ for $z$ in the
common domain of the three functions.
\end{theorem}
The function $R_\mu(z)=\varphi_\mu(1/z)$ is called the $R-$transform
of the probability measure $\mu$.
\bigskip

Let $\mu,\nu$ be two Borel 
probability measures, neither of them a point mass. Theorem \ref{Nevanlinna} 
assures us that the 
functions $F_\mu=F_1$ and $F_\nu=F_2$ satisfy the conditions of Theorem 
\ref{Firstsubordinationresult}.
With the notations from the proof of Theorem \ref{Firstsubordinationresult}, 
denote $F_3(z)=F_1(\omega(z))=F_2(\tilde{\omega}(z)),z\in\mathbb C^+.$ 
By Theorem \ref{Nevanlinna} and Proposition \ref{inversion+},
there exists a cone $\Gamma_{\alpha,\beta}$ such that $F_j$, $j\in\{1,
2,3\}$ have right inverses on $\Gamma_{\alpha,\beta}$ and $F_3^{-1}(\Gamma_{\alpha,\beta})\subseteq\mathbb C^+,$ 
Let us replace $z$ in the relation
$$F_3(z)+z=\omega(z)+\tilde{\omega}(z),\quad z\in\mathbb C^+$$
(equation (1.2)) by $F^{-1}_3(z).$
We obtain
$$z+F_3^{-1}(z)=\omega(F_3^{-1}(z))+\tilde{\omega}(F_3^{-1}(z))=
F_{1}^{-1}(z)+F_2^{-1}(z),$$
or, equivalently,
$$F_3^{-1}(z)-z=\varphi_\mu(z)+\varphi_\nu(z),\quad z\in\Gamma_{\alpha,\beta}.$$
By Theorem \ref{R-transform}, we conclude that $F_3=F_{\mu\boxplus\nu}.$

This allows us to give a new proof for the subordination result of Biane \cite{Biane1} (see also \cite{V3}, \cite{V4}, and \cite{V5}).
\begin{corollary}\label{Biane+}
Given two probability measures $\mu,\nu$ on $\mathbb R$, there exists a 
unique pair of analytic functions 
$\omega_1,\omega_2\colon\mathbb C^+\longrightarrow\mathbb C^+$ such that 
$F_{\mu\boxplus\nu}(z)=F_\mu(\omega_1(z))=F_\nu(\omega_2(z)),$ $z\in\mathbb C^+
$. Moreover, 
$\lim_{y\to+\infty}\omega_j(iy)/iy=1,$ $j\in\{1,2\},$ and
\begin{equation}
F_{\mu\boxplus\nu}(z)+z=\omega_1(z)+\omega_2(z),\quad z\in\mathbb C^+.
\label{eq1.3}
\end{equation}
\end{corollary}
\begin{proof} If one of the two measures, say $\mu$ is a point mass, $\mu=
\delta_a$, then the statements of the corollary are obviously true, with 
$\omega_1(z)=z-a$ and $\omega_2(z)=F_\nu(z-a)+a,$ $z\in\mathbb C^+.$
Otherwise, as we have already seen, if we let $F_1=F_\mu$ and $F_2=F_\nu$ in 
Theorem \ref{Firstsubordinationresult}, the functions 
$\omega_1=\omega$ and $\omega_2=\tilde{\omega}$ will
satisfy the requirement of the Corollary. Uniqueness is guaranteed by Theorem
\ref{Firstsubordinationresult}, while the equation (\ref{eq1.3}) has been 
established before.
\end{proof}

Next, we use the boundary properties of the subordination functions shown in 
Theorem \ref{Boundary+}  to describe the atomic, singular continuous, and 
absolutely continuous parts, with respect to the Lebesgue measure on $\mathbb R
$, of the free convolution $\mu\boxplus\nu$ of two probability measures $\mu,
\nu$ on the real line. Given a finite measure $\sigma$ on $\mathbb R$, denote
by $\sigma^{ac}$, $\sigma^{sc},$ and $\sigma^p$ the absolutely continuous, 
singular continuous, respectively atomic, parts of $\sigma$ with respect to the
Lebesgue measure.
The following lemma describes the behaviour of the Cauchy transform $G_\mu$
near points belonging to the support of the singular continuous part of the
probability measure $\mu.$ The result is not new, but since we don't know a
reference, we shall give its proof.
\begin{lemma}\label{DeLaValee}
Let $\mu$ be a Borel probability measure on $\mathbb R$. For $\mu^{sc}-$almost 
all $x\in\mathbb R$, the nontangential limit of the Cauchy transform $G_\mu$ of
$\mu$ at $x$ is infinite.
\end{lemma}
\begin{proof} 
It is enough to show that for $\mu^{sc}$-almost all $x\in\mathbb R$, the 
imaginary part of $G_\mu(x+iy)$ tends to infinity as $y$ approaches zero.
An elementary calculation shows that, for any $y>0$ and $x\in\mathbb R$,
$$\Im G_{\mu}(x+iy)  =  -\int_{\mathbb R}\frac{y}{(x-t)^2+y^2}d\mu(t).$$
Let us observe that $\frac{y}{(x-t)^2+y^2}\geq\frac{1}{2y}$ for all
$t\in\mathbb R$ such that $|x-t|\leq y.$
 Thus,
for any given $y>0$,
the following holds:
\begin{eqnarray*}
-\Im G_{\mu}(x+iy) & = & \int_{\mathbb R}\frac{y}{(x-t)^2+y^2}d\mu(t)\\
& = & \frac{1}{2y}\int_{\mathbb R}
\frac{2y^2}{(x-t)^2+y^2}d\mu(t)\\
& \geq & \frac{1}{2y}\int_{x-y}^{x+y}\frac{2y^2}{(x-t)^2+y^2}d\mu(t) \\
& \geq & \frac{\mu^{sc}((x-y,x+y])}{2y}.
\end{eqnarray*}
Now we can apply de la Vall\'{e}e Poussin's theorem (\cite{Saks}, Theorem 9.6) 
to
conclude that $\lim_{y\to0}\frac{\mu^{sc}((x-y,x+y])}{2y}=\infty$ for
$\mu^{sc}-$almost all $x\in\mathbb R.$ (see also \cite{Billingsley}, Theorem 
31.6)
\end{proof}

We shall first take care of the particular case not covered by Theorem
\ref{Boundary+}.

\begin{proposition}\label{Specialcase+} 
Let $\mu,\nu$ be two Borel probability measures on the real line, neither of 
them a point mass. Let $a\in\mathbb R$, and define 
$u_\mu(z)=F_\mu(z)-z+a$, $u_\nu(z)=F_\mu(z)-z+a$, $z\in\mathbb C^+.$
The function $u_\mu\circ u_\nu$ is a conformal automorphism of the upper 
half-plane if and only if 
$\mu$ and $\nu$ are convex combinations of two point masses. Moreover, in this 
case, $(\mu\boxplus\nu)^{sc}=0$, and the 
density of $(\mu\boxplus\nu)^{ac}$ with respect to the Lebesgue measure is 
analytic everywhere, with the exception of at most four points.
\end{proposition}
\begin{proof}
We claim that the analytic functions
$u_\mu,u_\nu\colon\mathbb C^+\longrightarrow\mathbb C^+,$ are both
conformal automorphisms of the upper half-plane.
 Indeed, observe first that both $u_\nu$ and $u_\mu$ are nonconstant, by
Remark \ref{ImF>Imz}. Let $k$ be the inverse with respect to
composition of $u_\mu\circ u_\nu,$
so that 
          $$u_\mu(u_\nu(k(z)))=z,\quad z\in\mathbb C^+.$$
This implies that $u_\nu$ is injective.
Applying $u_\nu$ to both sides of the above equality, we conclude 
that $u_\nu(u_\mu(w))=w$ for all $w$ in the open subset $(u_\nu\circ k)(
\mathbb C^+)$ of the upper half-plane, so, by analytic continuation, for all 
$w\in\mathbb C^+.$ This proves surjectivity of $u_\nu$ and thus our claim is 
proved.

Since 
          $$\lim_{y\to+\infty}\frac{u_\mu(iy)}{iy}=0,$$
we conclude that there exist $p,q,r\in\mathbb R$ such that 
$${\rm det}\left[\begin{array}{cc}
p & q \\
1 & r 
\end{array}\right]>0,\quad{\rm and}\quad u_\mu(z)=\frac{pz+q}{z+r},\quad z
\in\mathbb C^+.$$ Modulo a translation of $\mu$ and $\nu$, we may assume that 
$a=0$. This implies that 
   $$F_\mu(z)=\frac{z^2+z(p+r)+q}{z+r},$$
and thus $\mu$ is the convex combination of two point masses, at
$$\frac{-p-r+\sqrt{(p+r)^2-4q}}{2}\quad{\rm and}\quad
\frac{-p-r-\sqrt{(p+r)^2-4q}}{2}$$ with weights $$
\frac{r-p+\sqrt{(p+r)^2-4q}}{2\sqrt{(p+r)^2-4q}}\quad{\rm and}\quad
\frac{p-r+\sqrt{(p+r)^2-4q}}{2\sqrt{(p+r)^2-4q}}$$ respectively.
A similar statement holds for $\nu.$ 

Conversely, if $\nu=t\delta_u+(1-t)\delta_v,$ then a direct computation
shows that $F_\nu(z)=(z-v)(z-u)(z-tv-(1-t)u)^{-1},$ so that
$$u_\nu(z)=\frac{-(tu+(1-t)v)z+uv}{z-(tv+(1-t)u)},\quad z\in\mathbb C^+,$$
and
$${\rm det}\left[\begin{array}{cc}
-tu-(1-t)v & uv \\
1 & -tv-(1-t)u
\end{array}\right]=t(1-t)(u-v)^2>0,$$
for all $0<t<1,u\neq v$.
This proves the first statement of the 
proposition.

Assume now that $\mu=s\delta_\alpha+(1-s)\delta_u$ and $\nu=t\delta_\beta+
(1-t)\delta_v,$ with $0<s,t<1,$ and $\alpha,u,\beta,v\in\mathbb R$, $\alpha
\neq u,$ $\beta\neq v.$ Modulo a translation, we may assume that $\alpha=\beta=
0,$ and $u,v>0.$ 
Thus, we just need to compute 
$(s\delta_0+(1-s)\delta_u)\boxplus(t\delta_0+(1-t)\delta_v).$
As we have seen before,
$$F_\mu(z)=\frac{z^2-zu}{z-su},\quad F_\nu(z)=\frac{z^2-zv}{z-tv},\quad z\in
\mathbb C^+.$$
From equation (\ref{eq1.2}), we conclude that $\omega_1(z)$ must satisfy the
equation
$$F_\mu(\omega_1(z))=F_\nu(F_\mu(\omega_1(z))-\omega_1(z)+z),$$ which means
\begin{eqnarray}
\lefteqn{\frac{\omega_1(z)^2-\omega_1(z)u}{\omega_1(z)-su}=}\nonumber\\
&  & \frac{\left(\frac{\omega_1(z)^2-\omega_1(z)u}{\omega_1(z)-su}-\omega_1(z)
+z\right)^2-
\left(\frac{\omega_1(z)^2-\omega_1(z)u}{\omega_1(z)-su}-\omega_1(z)+z\right)v
}{\frac{\omega_1(z)^2-\omega_1(z)u}{\omega_1(z)-su}-\omega_1(z)+z-tv}\nonumber
\end{eqnarray}

Solving for $\omega_1$ gives as possible solutions
\begin{eqnarray}
\lefteqn{}\nonumber\\
& & -\frac{-z^2+z(u(1-2s)+v)+uv(s+t-1)}{2(z-tv-(1-s)u)}\nonumber\\
& & \mbox{}-\frac{\pm\sqrt{[z(z-u-v)+uv(1-s-t)]^2+4stuvz(z-u-v)}}{2(z-tv-(1-s)
u)}.\nonumber
\end{eqnarray}
Since $\lim_{y\to+\infty}\omega_1(iy)/iy=1$, we have
\begin{eqnarray}
\lefteqn{\omega_1(z)=}\nonumber\\
& & \frac{z^2-z(u(1-2s)+v)-uv(s+t-1)}{2(z-tv-(1-s)u)}\nonumber\\
& & \mbox{}+\frac{\sqrt{[z(z-u-v)+uv(1-s-t)]^2+4stuvz(z-u-v)}}{2(z-tv-(1-s)
u)},\nonumber
\end{eqnarray}
for all $z\in\mathbb C^+.$
By relation (\ref{eq1.2}),
$$G_{\mu\boxplus\nu}(z)=G_\mu(\omega_1(z))=
\frac{s}{\omega_1(z)}+\frac{1-s}{\omega_1(z)-u}.
$$
The only points where $G_{\mu\boxplus\nu}$ may fail to extend analytically
to $\mathbb R$ are the zeros of the equation
$$[z(z-u-v)+uv(1-s-t)]^2+4stuvz(z-u-v)=0$$
(i.e. when the expression under the square root sign is zero), and the zeros of
the equations 
\begin{eqnarray}
\lefteqn{[z^2-z(u(1-2s)+v)-uv(s+t-1)]^2=}\nonumber\\
& & [z(z-u-v)+uv(1-s-t)]^2+4stuvz(z-u-v),\nonumber
\end{eqnarray}
and
\begin{eqnarray}
\lefteqn{[z^2-z(u(1-2s)+v)-uv(s+t-1)-2u(z-tv-(1-s)u)]^2=}\nonumber\\
& & [z(z-u-v)+uv(1-s-t)]^2+4stuvz(z-u-v),\nonumber
\end{eqnarray}
which occur at the denominator of $G_{\mu\boxplus\nu}.$
We have used Maple to determine that the four solutions for the first equation
(denote them by $r_j,j\in\{1,2,3,4\}$) are 
$$\frac{1}{2}\left(u+v\pm\sqrt{(u-v)^2+4uv(t+s-2st)\pm
8uv\sqrt{ts(1-t)(1-s)}}\right).$$
Whether all four of these roots are distinct or not depends on the
actual values of $s,t,u,v.$ However, it is trivial to observe that they 
are all real; indeed, $u,v>0$ by hypothesis, and an elementary computation
shows that $t+s-2st>0$ and $t+s-2st\geq2\sqrt{ts(1-t)(1-s)},$ with equality
if and only if $s=t$. Notice also that it is not possible that exactly three 
of the roots are distinct. 
Also by using Maple, we have determined that the second and third equations 
have as solutions the numbers $0,v,(1-s)u+tv$, and $u,u+v,(1-s)u+tv,$ 
respectively. We can thus write (with the help of Maple):
\begin{eqnarray}
\lefteqn{G_{\mu\boxplus\nu}(z)=}\nonumber\\
& & \frac{z^2-z(u(1-2s)+v)-uv(s+t-1)-\sqrt{
\prod_{j=1}^4(z-r_j)}}{
2uz(z-v)}\nonumber\\
& & \mbox{}-\frac{z^2-z((3-2s)u+v)+uv(t-s+1)+2u^2(1-s)}{2u(z-u)(z-u-v)}\nonumber\\
& & \mbox{}+\frac{\sqrt{\prod_{j=1}^4(z-r_j)}}{2u(z-u)(z-u-v)}.\nonumber
\end{eqnarray}
As it can be seen from the above formula, the support of $(\mu\boxplus\nu)^{ac}
$ coincides with the set where the product $\prod_{j=1}^4(z-r_j),$ $z\in
\mathbb R$, is negative, and, since $G_{\mu\boxplus\nu}$ can be infinite only
at the points  $0,u,v,u+v$,  the measure $\mu\boxplus\nu$ has at most four
atoms, which may occur exactly in $0,u,v$, or $u+v.$ This concludes the proof.
\end{proof}

We can now prove the main regularity result of this section.

\begin{theorem}\label{Reg+}
Let $\mu,\nu$ be two Borel probability measures on the real line, neither of 
them a point mass. Then the following statements hold:
\begin{enumerate}
\item[{\rm (1)}] The point $a\in\mathbb R$ is an atom of the measure $\mu
\boxplus\nu$ if and only if there exist $b,c\in\mathbb R$ such that
$a=b+c$ and $\mu(\{b\})+\nu(\{c\})>1.$ Moreover, $(\mu\boxplus\nu)(\{a\})=
\mu(\{b\})+\nu(\{c\})-1.$ 
\item[{\rm (2)}] The absolutely continuous part of $\mu\boxplus\nu$ with 
respect to the Lebesgue measure is always nonzero and its density is locally
analytic on the complement of a closed set of zero Lebesgue measure.
\item[{\rm(3)}] Assume in addition that 
${\rm supp}(\mu)$ is compact and $\mathbb R\setminus{\rm supp}(\nu)$ is
nonempty. Then the singular continuous part of $\mu\boxplus\nu$ is zero.
\end{enumerate}
\end{theorem}

\begin{proof}
Part (1) of the theorem is due to Bercovici and Voiculescu (see  
\cite{BercoviciVoiculescuRegQ}, Theorem 7.4).
We shall proceed with the proof of part (2).
Suppose that $\mu\boxplus\nu$ is purely singular, and thus for almost
all $x\in\mathbb R$ with respect to the Lebesgue measure, we have
$$\lim_{y\to0}\Im F_{\mu\boxplus\nu}(x+iy)=\lim_{y\to0}\Im 
G_{\mu\boxplus\nu}(x+iy)=0.$$ 
By part (1), we are assured that $\mu\boxplus\nu$ cannot be purely atomic, so 
we must have $(\mu\boxplus\nu)^{sc}\neq0,$
and hence, by Lemma \ref{DeLaValee}, $\lim_{y\to0} F_{\mu\boxplus\nu}(x+iy)=0$ 
for uncountably many $x\in\mathbb R.$
Theorem \ref{Seidel} applied to the function $F_{\mu\boxplus\nu}$ yields
a point $x_0\in\mathbb R$ such that $C(F_{\mu\boxplus\nu},x_0)=
\overline{\mathbb C^+}.$ Using relation (\ref{eq1.3}) we conclude that at least
one of 
$C(\omega_1,x_0),C(\omega_2,x_0)$, will intersect the upper
half-plane. (Recall that $\omega_1,\omega_2$ are the
subordination functions from Corollary \ref{Biane+}.)
But now we can apply Theorem \ref{Boundary+} (1) and Proposition \ref{Specialcase+} to obtain a contradiction. 
Thus, $\mu\boxplus\nu$ cannot be purely singular.

Next we prove that there exists a closed subset of $\mathbb R$ of zero Lebesgue
measure on whose complement
the density $f(x)=\frac{d(\mu\boxplus\nu)^{ac}(x)}{dx}$ is analytic. 
By Theorem \ref{Fatou}, there exists a subset $E$ of $\mathbb R$ of zero 
Lebesgue measure such that for all $x\in\mathbb R\setminus E$ the limits
$\lim_{y\to0}F_{\mu\boxplus\nu}(x+iy)$, $\lim_{y\to0}\omega_j(x+iy),$ 
$j\in\{1,2\}$ exist and are finite. Also, for almost all $x\in{\rm supp}((\mu
\boxplus\nu)^{ac})\setminus E$, with respect to $(\mu\boxplus\nu)^{ac},$ we 
have $\lim_{y\to0}F_{\mu\boxplus\nu}(x+iy)\in\mathbb C^+.$ By relation (\ref{eq1.3}), at 
least one of $\lim_{y\to0}\omega_j(x+iy),$ $j\in\{1,2\},$ 
must also be in $\mathbb C^+$. For definiteness,
assume $\omega_1(x)=\lim_{y\to0}\omega_1(x+iy)\in\mathbb C^+.$ Part (1) of
Theorem \ref{Boundary+} and Proposition \ref{Specialcase+} assure us that 
$\omega_1$ extends analytically through
$x$. Since $\omega_1(x)\in\mathbb C^+,$ the function 
$F_{\mu\boxplus\nu}=F_{\mu}\circ\omega_1$ also extends analytically through
$x$. We conclude that the density of $(\mu\boxplus\nu)^{ac}$ must be analytic 
in $x$. On the other hand, if there exists an interval $J\subseteq\mathbb R
\setminus{\rm supp}((\mu\boxplus\nu)^{ac}),$ then $\lim_{y\to0}\Im 
F_{\mu\boxplus\nu}(x+iy)=0$ for almost all $x\in J$ with respect to the 
Lebesgue emasure. The same argument used in the proof of the existence
of an absolutely continuous part of $\mu\boxplus\nu$ shows that $F_{\mu\boxplus
\nu}$ extends meromorphically through $J$. Of course,
the set $A$ of points $x$
where $f(x)=\frac{d(\mu\boxplus\nu)^{ac}(x)}{dx}$ is analytic is open in 
$\mathbb R$, and its complement
is of zero Lebesgue measure.

To prove (3), observe that our hypotheses on the supports of $\mu$ and $\nu$
allow us to apply either part (2) of Theorem \ref{Boundary+}, or 
Proposition \ref{Specialcase+} to conclude that the subordination functions 
$\omega_1$ and $\omega_2$ extend continuously to $\mathbb R$ as functions with 
values in the extended complex plane $\mathbb C\cup\{\infty\}$. By relation 
(\ref{eq1.3}), the same must hold for $F_{\mu\boxplus\nu}.$ We shall denote by 
the same symbols the continuous extensions of the three functions to 
$\mathbb C^+\cup\mathbb R.$

Let us assume that $(\mu\boxplus\nu)^{sc}\neq0.$ Lemma \ref{DeLaValee} allows 
us to find an uncountable set $H\subseteq{\rm supp}((\mu\boxplus\nu)^{sc})$ 
such that no $x\in H$ is an atom of $\mu\boxplus\nu$, and
$$\lim_{\stackrel{z\longrightarrow x}{\sphericalangle}}F_{\mu\boxplus\nu}(z)
=\lim_{z\to x}F_{\mu\boxplus\nu}(z)=0.$$
We claim that for any $x\in H$
\begin{enumerate}
\item[({\it i})] $\omega_1(x),\omega_2(x)\neq\infty,$
\item[({\it ii})] $\sphericalangle\lim_{z\to\omega_1(x)}F_\mu(z)=
\sphericalangle\lim_{z\to\omega_2(x)}F_\nu(z)=0,$
\item[({\it iii})] $\mu(\{\omega_1(x)\})+\nu(\{\omega_2(x)\})=1.$
\end{enumerate}
Part ({\it i}) is an easy consequence of Theorems \ref{Nevanlinna} and 
\ref{Lindelof}. Indeed, by Theorem \ref{Nevanlinna}, we have that 
$\sphericalangle\lim_{z\to\infty}F_\mu(z)=\infty$.
Assume that $\omega_1(x)=\infty.$ Then, by our assumption on $x$, Corollary 
\ref{Biane+}, and Theorem \ref{Lindelof}, we have
\begin{eqnarray*}
0 & = & \lim_{z\to x}F_{\mu\boxplus\nu}(z)\\
& = & \lim_{z\to x}F_\mu(\omega_1(z))\\
& = & \lim_{\stackrel{w\to\infty}{w\in\omega_1(x+i(0,1))}}F_\mu(w)\\
& = & \lim_{\stackrel{w\longrightarrow\infty}{\sphericalangle}}F_\mu(w)\\
& = & \infty,
\end{eqnarray*}
which is absurd. A similar argument shows that $\omega_2(x)\neq\infty.$

By Theorem \ref{Lindelof}, we have
$$\lim_{\stackrel{w\longrightarrow\omega_1(x)}{\sphericalangle\ \ \ }}F_\mu(w)=
\lim_{\stackrel{w\to\omega_1(x)}{w\in\omega_1(x+i(0,1))}}F_\mu(w)=
\lim_{y\to 0}F_\mu(\omega_1(x+iy))=\lim_{z\to x}
F_{\mu\boxplus\nu}(z)=0.$$ The same argument applied to $F_\nu$ and $\omega_2$
in place of $F_\mu$ and $\omega_1$ yields $\sphericalangle\lim_{w\to\omega_2(x)
}F_\nu(w)=0.$ This proves part ({\it ii}) of our claim.

We prove now part ({\it iii}). Observe first that, by dominated convergence and
by Theorem \ref{JuliaCaratheodory} applied to self-maps of $\mathbb C^+,$ for 
any $v\in\mathbb R$ such that $\sphericalangle\lim_{z\to v}F_\mu(z)=0,$ we have
\begin{eqnarray*}
\frac{1}{\mu(\{v\})} & = & \lim_{\stackrel{z\longrightarrow v}{\sphericalangle}}\left[(z-v)\int_\mathbb R\frac{d\mu(t)}{z-t}\right]^{-1}\\
& = & \lim_{\stackrel{z\longrightarrow v}{\sphericalangle}}\frac{F_\mu(z)}{z-v}\\
& = & \liminf_{z\to v}\frac{\Im F_\mu(z)}{\Im z}.
\end{eqnarray*}
Applying this to $v=\omega_1(x)$ yields
$$\frac{1}{\mu(\{\omega_1(x)\})}=\liminf_{z\to\omega_1(x)}\frac{\Im F_\mu(z)}{
\Im z}\leq\liminf_{z\to x}\frac{\Im F_\mu(\omega_1(z))}{\Im\omega_1(z)},$$
and similarily
$$\frac{1}{\nu(\{\omega_2(x)\})}\leq\liminf_{z\to x}\frac{\Im 
F_\nu(\omega_2(z))}{\Im\omega_2(z)},$$
for all $x\in H$ (we use here the convention $1/0=+\infty$).

From relation (\ref{eq1.3}) we obtain
$$\frac{\Im F_{\mu\boxplus\nu}(z)}{\Im \omega_1(z)}+\frac{\Im z}{\Im\omega_1(z
)}-1=\frac{\Im\omega_2(z)}{\Im\omega_1(z)},\quad z\in\mathbb C^+.$$

The following chain of inequalities holds:
\begin{eqnarray*}
\frac{1}{\nu(\{\omega_2(x)\})}-1 & \leq & \liminf_{z\to x}\frac{\Im
F_\nu(\omega_2(z))}{\Im\omega_2(z)}-1\\
& \leq & \limsup_{z\to x}\left(\frac{\Im F_\nu(\omega_2(z))}{\Im\omega_2(z)}+
\frac{\Im z}{\Im \omega_2(z)}-1\right)\\
& = & \limsup_{z\to x}\frac{\Im \omega_1(z)}{\Im \omega_2(z)}\\
& = & \left(\liminf_{z\to x}\frac{\Im\omega_2(z)}{\Im \omega_1(z)}\right)^{-1}\\
& = & \left(\liminf_{z\to x}\left(\frac{\Im F_\mu(\omega_1(z))}{\Im
\omega_1(z)}+\frac{\Im z}{\Im \omega_1(z)}\right)-1\right)^{-1}\\
& \leq & \left(\liminf_{z\to x}\frac{\Im F_\mu(\omega_1(z))}{\Im \omega_1(z)}+
\liminf_{z\to x}\frac{\Im z}{\Im \omega_1(z)}-1\right)^{-1}\\
& \leq & \left(\liminf_{z\to x}\frac{\Im F_\mu(\omega_1(z))}{\Im \omega_1(z)}-1
\right)^{-1}\\
& \leq & \left(\frac{1}{\mu(\{\omega_1(z)\})}-1\right)^{-1}.
\end{eqnarray*}
We have assumed that $\mu$ and $\nu$ are not point masses, so the above 
implies that $1<\frac{1}{\mu(\{\omega_1(z)\})},\frac{1}{\nu(\{\omega_2(z)\})}<
+\infty.$ Thus, multiplying the above inequality by 
$\frac{1}{\mu(\{\omega_1(z)\})}-1$ will give 
$$(1-\mu(\{\omega_1(x)\}))(1-\nu(\{\omega_2(x)\}))\leq\mu(\{\omega_1(x)\})\nu(
\{\omega_2(x)\}),$$ or, equivalently,
$$\mu(\{\omega_1(x)\})+\nu(\{\omega_2(x)\})\geq1.$$
Since $x\in H$, part (1) of the theorem tells us that the inequality above must
be an equality. This proves the last point of our claim.

Using relation (\ref{eq1.3}) and the fact that $F_{\mu\boxplus\nu}(x)=0$,
we obtain $$\omega_1(x)+\omega_2(x)=x\quad {\rm for\ all}\quad x\in H.$$
Since any probability can have at most countably many atoms, this, together 
with part ({\it iii}) of our claim contradicts the fact that $H$ is uncountable
and concludes the proof.
\end{proof}

\section[]{
 Regularity
for the free multiplicative convolution of probability measures on
the positive half-line}

Let $\sigma$ be a finite Borel measure on $[0,+\infty).$ We define the 
analytic function $$\psi_\sigma(z)=\int_{[0,+\infty)}\frac{zt}{1-zt}d\sigma(t),
\quad z\in\mathbb C\setminus[0,+\infty).$$
The function $\psi_\sigma$ is called the moment generating function of $\sigma.
$ A simple computation shows that $G_\sigma(1/z)=z(\psi_\sigma(z)+\sigma([0,+
\infty))),$ $z\in\mathbb C\setminus[0,+\infty).$
Using this formula,
one can prove the following analogue of Theorem \ref{Cauchy} for measures 
supported on the positive half-line (see 
\cite{BercoviciVoiculescuIUMJ}, Proposition 6.1).

\begin{proposition}\label{Moment}
Let $\psi:\mathbb C\setminus[0,+\infty)\longrightarrow\mathbb C$ be an analytic
function such that $\psi(\overline{z})=
\overline{\psi(z)}$ for all $z\in\mathbb C\setminus[0,+\infty)$. The following 
two conditions are equivalent.
\begin{enumerate}
\item[{\rm (1)}] There exists a finite measure $\sigma$ on $[0,+\infty)$ such
that
$\psi=\psi_\sigma$.
\item[{\rm(2)}] $\lim_{y\uparrow0}\psi(y)=0$ and there exists a number $C>0$
such that $z(\psi(z)+C)\in\mathbb C^+$ for all $z\in\mathbb C^+$.
\end{enumerate}
Moreover, $C=\sigma([0,+\infty)).$
\end{proposition}
Observe in particular that, unless $\sigma$ is concentrated at zero, we have
$0\not\in\psi_\sigma(\mathbb C\setminus[0,+\infty)),$ $\psi_\sigma((-\infty,0))
\subseteq(-\infty,0),$ and $\psi_\sigma(\mathbb C^+)\subseteq\mathbb C^+.$
For analytic functions $\psi$ satisfying these three conditions it is possible 
to define
a continuous function $\arg \psi:\mathbb C\setminus[0,+\infty)\longrightarrow
\mathbb R$ such that
$\exp(i\arg \psi(z))=\psi(z)/|\psi(z)|$, and $\arg \psi=\pi$ on the negative 
half-line.  We will also consider
the limit $\psi(0-)=\lim_{z\to0,z<0}\psi(z)$.  Observe that
$\psi(0-)$ exists
if $\psi(\mathbb C^+)\subseteq\mathbb C^+$.
In this case it is easy to see that $\psi'(z)>0$ for all $z\in(-\infty,0)$;
indeed, if $z<0$, $\psi'(z)$ is a real number so that
$$\psi'(z)=\lim_{y\downarrow0}\Re\frac{\psi(z+iy)-\psi(z)}{iy}=
\lim_{y\downarrow0}\frac{\Im \psi(z+iy)}{y}\ge0,$$ and $\psi'(z)=0$
would imply that the image under $\psi$ of a small half-disk $\{w:\Im
w>0,|w-z|<\varepsilon\}$ contains numbers in $-\mathbb C^+$.
An immediate consequence of this fact is that the function $\psi_\sigma$ 
defined above never takes the value $-\sigma([0,+\infty)).$

Let now $\mu$ be a probability measure on $[0,+\infty).$ It will be convenient 
to consider the function 
$$\eta_\mu(z)=\frac{\psi_\mu(z)}{1+\psi_\mu(z)},\quad z\in\mathbb C\setminus
[0,+\infty).$$
An analogue of Theorem \ref{Cauchyalt} for the function $\eta_\mu$ is the 
following proposition (see \cite{BBercoviciMult}, Proposition 2.1).

\begin{proposition}\label{eta}
Let $\eta\colon\mathbb C\setminus[0,+\infty)\longrightarrow\mathbb C\setminus
\{0\}$ be an  analytic function such that $\eta(\overline{z})=
\overline{\eta(z)}$ for all $z\in\mathbb C\setminus[0,+\infty)$. The following 
two conditions are equivalent.
\begin{enumerate}
\item[{\rm (1)}] There exists a probability measure $\mu\ne\delta_0$ on
$[0,+\infty)$ such that
$\eta=\eta_\mu$.
\item[{\rm(2)}] $\eta(0-)=0$ and $\arg\eta(z)\in[\arg z,\pi)$ for all
$z\in\mathbb C^+$.
\end{enumerate}
Moreover, $\mu\neq\delta_0$ is a point mass if and only if there exists a $z\in
\mathbb C^+$ such that $\arg\eta_\mu(z)=\arg z.$
\end{proposition}

\begin{proof}
Observe first that, for $z\in\mathbb C^+$, the set
$$C_z=\left\{\frac{tz}{1-tz}:t>0\right\}$$
is a circular arc contained in $\mathbb C^+$, with endpoints 0 and $-1$, and
tangent to the
line $\{tz:t\in\mathbb R\}$. Thus, given $\mu\ne\delta_0$, $\psi_\mu(z)$ will
belong to
the convex hull $\text{co}(C_z)$, and therefore
$$\eta_\mu(z)\in\left\{\frac{w}{w+1}:w\in\text{co}(C_z)\right\}.$$
This last set is precisely $\{\lambda\in\mathbb C^+:\arg\lambda\ge\arg z\}$, and
we conclude
that $\eta_\mu$ satisfies the conditions in (2).  Conversely, assume that $\eta
$ satisfies (2),
and set $\psi=\eta/(1-\eta)$, so that $\eta=\psi/(1+\psi)$. For $z\in\mathbb
C^+$ we have
then $\eta(z)\in\{w/(1+w):w\in\text{co}(C_z)\}$, and thus
$\psi(z)\in\text{co}(C_z)$. Thus
$1+\psi(z)\in1+\text{co}(C_z)$ which implies that $1+\psi(z)\in\mathbb C^+$ and
$\arg(1+\psi(z))>-\arg z$ (by the symmetry of $C_z$ about the line $x=-1/2$). 
We conclude that
$z(1+\psi(z))\in\mathbb C^+$, and the result follows from the preceding
proposition.

If $\mu=\delta_a$ for some $a\in(0,+\infty),$ then $\eta_\mu(z)=az$, so
$\arg\eta_\mu(z)=\arg z$ for all $z\in \mathbb C^+.$ Assume that 
$\arg\eta_\mu(z_0)=\arg z_0$ for some $z_0$ in the upper half-plane. Then
$\arg\psi_\mu(z_0)/(1+\psi_\mu(z_0))=\arg z_0,$ so we conclude that 
there exists a number $\ell\in(0,+\infty)$ such that $\psi_\mu(z_0)/(1+\psi_\mu
(z_0))=\ell z_0,$ or, equivalently,
$$\psi_\mu(z_0)=\frac{\ell z_0}{1-\ell z_0}.$$ We conclude that
$\psi_\mu(z_0)\in C_{z_0},$ and thus $\mu=\delta_\ell$.
\end{proof}
\begin{remark}\label{Inj}{\rm
An immediate consequence of the definition of $\eta_\mu$ and of Proposition 
\ref{Moment} is that
for any $\mu\neq\delta_0,$ ${\eta_{\mu}}((-\infty,0))\subseteq(-\infty,0),$ 
and thus, as shown before, $\eta_\mu|_{(-\infty,0)}$ is injective.}
\end{remark}

Next we prove a subordination result for analytic self-maps of the 
slit complex plane.
\begin{theorem}\label{Secondsubordinationresult}
Let $\eta_j\colon\mathbb C\setminus[0,+\infty)\longrightarrow
\mathbb C\setminus[0,+\infty)$, $j\in\{1,2\}$ be two analytic functions
satisfying conditions
\begin{trivlist}
\item[{\rm \ \ (a)}] $\lim_{x\to0,x<0}\eta_j(x)=0,$ $j\in\{1,2\}$;
\item[{\rm \ \ (b)}] $\pi>\arg\eta_j(z)>\arg z$, $z\in\mathbb C^+,$ $j\in\{1,2\}$;
\item[{\rm \ \ (c)}] $\eta_j(\overline{z})=\overline{\eta_j(z)},$ $z\in\mathbb C
\setminus[0,+\infty)$, $j\in\{1,2\}$.
\end{trivlist}
Define $h_j\colon\mathbb C\setminus[0,+\infty)\longrightarrow\mathbb C$ by
$h_j(z)=\frac{\eta_j(z)}{z},$
and let $f\colon(\mathbb C^+\cup\mathbb R\setminus\{0\})\times\mathbb C^+
\longrightarrow\mathbb C$,
$$f(z,w)=h_2(zh_1(zw)).$$
Then
\begin{enumerate}
\item[{\rm (1)}] The function $f$ is well-defined and takes values in
$\mathbb C^+.$
\item[{\rm (2)}] There exists a unique analytic self-map $h$ of
$\mathbb C^+$ such that $f(z,h(z))=h(z)$, $z\in\mathbb C^+.$
\end{enumerate}
\end{theorem} 

\begin{proof}
We start by proving (1).
It follows from (b) and (c) that the functions $h_j$ map the upper half-plane 
into itself,
$h_j(\overline{z})=\overline{h_j(z)},$ and $h_j((-\infty,0))\subset(0,+\infty),
$ $j\in\{1,2\}.$ Thus, if $z\in(0,+\infty)$ it is obvious that $f(z,w)$ is well
defined and takes values in $\mathbb C^+$. If $z\in(-\infty,0),$ then
$h_1(zw)\in-\mathbb C^+,$ so that $zh_1(zw)\in\mathbb C^+,$ and we conclude
again that $f(z,w)\in\mathbb C^+.$

Let now $z\in\mathbb C^+.$ If $w\in\mathbb C^+$ is such that $\arg(zw)\leq\pi,$
then
$$\arg zh_1(zw)=\arg z+\arg\eta_1(zw)-\arg(zw)\geq\arg z>0,$$
and $$\arg zh_1(zw)=\arg z+\arg\eta_1(zw)-\arg(zw)<\pi-\arg w<\pi.$$
If $w\in\mathbb C^+$ is such that $\arg(zw)>\pi,$ then, using (b) and (c), we 
have
 $$\arg\eta_1(zw)=2\pi-\arg\overline{\eta_1(zw)}=2\pi-\arg\eta_1(\overline{zw})
<2\pi-\arg(\overline{zw})=\arg(zw),$$
so
$$\arg zh_1(zw)=\arg\eta_1(wz)-\arg w<\arg(wz)-\arg w=\arg z<\pi,$$
and
$$\arg zh_1(zw)=\arg\eta_1(zw)-\arg w>\pi-\arg w>0,$$
by (c).
We conclude that $zh_1(zw)\in\mathbb C^+$ for all $(z,w)\in(\mathbb C^+\cup
\mathbb R\setminus\{0\})\times\mathbb C^+$, and thus $f$ is well defined and
its range is included in the upper half-plane.

We prove now (2).
For each $z\in\mathbb C^+$ denote
$f_z\colon\mathbb C^+\longrightarrow\mathbb C^+$, $f_z(w)=f(z,w).$ We claim
that $f_z$ must have its Denjoy-Wolff point in $\mathbb C^+.$

Indeed, for any $r\in(0,+\infty),$ $\lim_{w\to r}f(z,w)=h_2(zh_1(rz))\in
\mathbb C^+,$ so $r$ cannot be a fixed point of $f_z$, so the Denjoy-Wolff
point of $f_z$ does not belong to the interval $(0,+\infty)$.

Assume now that zero is the Denjoy-Wolff point of $f_z$. Then, by 
the definition of the Denjoy-Wolff point, we have
 $\lim_{y\to0}
f_z(iy)=0$, and $\lim_{y\to0}\frac{f_z(iy)}{iy}$ must belong to $(0,1].$
But
$$
\lim_{y\to0}\frac{f_z(iy)}{iy}  =  \lim_{y\to0}\frac{h_2(zh_1(ziy))}{h_1(ziy)
}
=  \lim_{y\to0}\frac{\eta_2\left(z\frac{\eta_1(ziy)}{ziy}\right)}{\eta_1(ziy)}
.$$
By Theorem \ref{JuliaCaratheodory}, 
$r=\lim_{y\downarrow0}\frac{\eta_1(ziy)}{ziy}$ exists in $(0,+\infty)\cup
\{\infty\}.$

If $r\in(0,+\infty),$ then
$\lim_{y\downarrow0}\eta_2\left(z\frac{\eta_1(ziy)}{ziy}\right)=\eta_2(rz)$, 
so, using
(a) together with Theorem \ref{Lindelof}, we obtain
$\lim_{y\to0}\frac{f_z(iy)}{iy} =\frac{\eta_2(rz)}{0}=\infty,$ which 
contradicts the fact that $\lim_{y\to0}\frac{f_z(iy)}{iy}\in(0,1]$.

We conclude that $r$ must be infinity. Since, by property (b), we have 
$\arg\left(z\frac{\eta_1(zw)}{zw}
\right)\in(\arg z,\pi+\arg z)$ for any $w\in\mathbb C^+,$ we obtain 
$$\lim_{y\downarrow0}\eta_2\left(z\frac{\eta_1(ziy)}{ziy}\right)=
\lim_{x\to-\infty}\eta_2(x)\in[-\infty,0),$$ by Theorem \ref{Lindelof} and
the monotonicity of $\eta_2|_{(-\infty,0)}$.
Property (a) yields again
$\lim_{y\to0}\frac{f_z(iy)}{iy}  = \infty,$ which contradicts the fact 
that $\lim_{y\to0}\frac{f_z(iy)}{iy}\in(0,1].$
We conclude that zero is not the Denjoy-Wolff point of $f_z$.

Suppose that $f_z$ has infinity as its Denjoy-Wolff point.
Then we must have
$$\lim_{y\to+\infty}f_z(iy)=\infty,\ {\rm and}\ \lim_{y\to+\infty}
\frac{f_z(iy)}{iy}\geq1.$$
However, $$\lim_{y\to+\infty}
\frac{f_z(iy)}{iy}=\lim_{y\to+\infty}
\frac{\eta_2\left(z\frac{\eta_1(ziy)}{ziy}\right)}{\eta_1(ziy)}.$$
Theorem \ref{JuliaCaratheodory} assures us that 
$l=\lim_{y\to+\infty}\frac{\eta_1(ziy)}{ziy}
\in[0,+\infty).$ If $l>0$, then we have  $\lim_{y\to+\infty}\eta_2\left(z
\frac{\eta_1(ziy)}{ziy}\right)=\eta_2(zl),$ so
$$\lim_{y\to+\infty}
\frac{f_z(iy)}{iy}=\lim_{y\to+\infty}
\frac{\eta_2\left(z\frac{\eta_1(ziy)}{ziy}\right)}{\eta_1(ziy)}=\frac{\eta_2(z
l)}{\lim_{y\to+\infty}\eta_1(ziy)}=0,$$
a contradiction. If $l=0,$ then, by Theorem \ref{Lindelof}, the limit
$$\lim_{y\to+\infty}\eta_2\left(z\frac{\eta_1(ziy)}{ziy}\right),$$ 
if it exists, must equal zero.
Since
$$\lim_{y\to+\infty}
\frac{f_z(iy)}{iy}=\lim_{y\to+\infty}
\frac{\eta_2\left(z\frac{\eta_1(ziy)}{ziy}\right)}{\eta_1(ziy)}$$
and, by Theorem \ref{Lindelof}, $\lim_{y\to+\infty}\eta_1(ziy)=
\lim_{z\to-\infty}\eta_1(z)<0,$ 
we conclude that $$\lim_{y\to+\infty}\frac{f_z(iy)}{iy}=0.$$ 
Thus, infinity cannot be the Denjoy-Wolff point of $f_z$.

We still need to show that no point on the negative half-line is the 
Denjoy-Wolff point of $f_z.$ Suppose to the contrary that 
$$\lim_{\stackrel{w\longrightarrow r}{\sphericalangle}}f_z(w)=r$$
for some $r<0.$ As we have already seen, 
$$\lim_{\stackrel{w\longrightarrow r}{\sphericalangle}}h_1(zw)=
\lim_{\stackrel{w\longrightarrow r}{\sphericalangle}}\frac{\eta_1(zw)}{zw}
=\frac{\eta_1(zr)}{zr}\in-\mathbb C^+,$$
so 
$$\lim_{\stackrel{w\longrightarrow r}{\sphericalangle}}f_z(w)=
h_2\left(\frac{1}{r}\eta_1(zr)\right)\in\mathbb C^+.$$
Thus $r$ cannot be a fixed point of $f_z$. This proved our claim.

Let $h(z)\in\mathbb C^+$ denote the Denjoy-Wolff point of $f_z.$
Since we have shown in part (1) that on the one hand $f_z(r)$ takes values in 
the upper half-plane for any $r>0$, and on the other $f_z((-\infty,0))
\subseteq(-\infty,0),$ we conclude that $f_z$ is not a conformal automorphism 
of $\mathbb C^+.$
Thus, Theorem \ref{DenjoyWolff} assures us that the derivative of $f_z$ in the 
point $h(z)$ has absolute value strictly less than one, $|f_z'(h(z))|<1.$
By applying the implicit function theorem to the function
$(z,w)\mapsto f(z,w)-w$, we conclude that $h(z)$ depends
analytically on $z$. Thus, we have determined the existence of $\omega$. The
uniqueness follows immediatly from the uniqueness of the Denjoy-Wolff point.
This proves part (2) of the theorem.
\end{proof}

\begin{corollary}\label{subordinationwithomega}
Let $\eta_j$,
 $j\in\{1,2\}$, be two
analytic functions as in Theorem \ref{Secondsubordinationresult}.
Then there exists a unique analytic function $$\omega\colon\mathbb C\setminus
[0,+\infty)\longrightarrow\mathbb C\setminus[0,+\infty)$$ such that
$$\omega(z)=\frac{\omega(z)\eta_2\left(z\frac{\eta_1(\omega(z))}{\omega(z)}
\right)}{\eta_1(\omega(z))},\quad z\in\mathbb C\setminus[0,+\infty).$$
Moreover, $\arg z<\arg\omega(z)<\pi$ for all $z\in\mathbb C^+$.
\end{corollary}
\begin{proof}
With the notations from Theorem \ref{Secondsubordinationresult}, let 
$\omega(z)=zh(z),$ $z\in
\mathbb C^+.$ It is clear that $\arg\omega(z)>\arg z$ for all numbers $z$ in
the upper half-plane. Using the relation $f(z,h(z))=h(z)$, we have:
\begin{eqnarray*}
\omega(z) & = & zf(z,h(z))\\
& = & zh_2(zh_1(zh(z)))\\
& = & z\frac{\eta_2\left(z\frac{\eta_1(h(z))}{\omega(z)}\right)}{z\frac{
\eta_1(\omega(z))}{\omega(z)}}\\
& = & \frac{\omega(z)\eta_2\left(z\frac{\eta_1(\omega(z))}{\omega(z)}\right)}{
\eta_1(\omega(z))}.
\end{eqnarray*}
The above equality holds for all $z\in\mathbb C^+.$
We conclude that
$\omega(z)\in\mathbb C^+,$ and $\arg\omega(z)>\arg z.$

Denote $A_z=\{w\in\mathbb C^+\colon\arg w>\arg z\},$
and observe that $\omega(z)$ is a fixed point for the function
$$g_z\colon A_z\longrightarrow A_z,\quad g_z(w)=w
\frac{\eta_2\left(z\frac{\eta_1(w)}{w}\right)}{\eta_1(w)}.$$
Since the functions $\eta_j$ satisfy the condition (b) from Theorem
\ref{Secondsubordinationresult}, the function $g_z$ cannot be the identity 
function. By the Denjoy-Wolff theorem (Theorem \ref{DenjoyWolff}), we obtain 
the uniqueness of $\omega.$

Since $h$ extends continuously to the negative half-line, it is now obvious 
that so does $\omega.$
By the Schwarz reflection principle, we can extend $\omega$ to the whole
$\mathbb C\setminus[0,+\infty),$ and its range will be included in 
$\mathbb C\setminus[0,+\infty).$ This concludes the proof.
\end{proof}

As in the first section, let us observe that Theorem 
\ref{Secondsubordinationresult} and Corollary \ref{subordinationwithomega} 
provide an analytic function 
$\tilde{\omega}\colon\mathbb C\setminus[0,+\infty)
\longrightarrow\mathbb C\setminus[0,+\infty),$ $\tilde{\omega}(z)=z\eta_1(\omega(z))/\omega(z),$
satisfying the conditions 
$$\tilde{\omega}(z)=\frac{\tilde{\omega}(z)\eta_1\left(z\frac{\eta_2(
\tilde{\omega}(z))}{\tilde{\omega}(z)}
\right)}{\eta_2(\tilde{\omega}(z))},\quad z\in\mathbb C\setminus[0,+\infty),$$
and $\arg\tilde{\omega}(z)>\arg z,$ $z\in\mathbb C^+.$
The pair of functions $(\omega,\tilde{\omega})$ is uniquely 
determined by the equations
\begin{equation}
\eta_1(\omega(z))=\eta_2(\tilde{\omega}(z))
=\frac{1}{z}\omega(z)\tilde{\omega}(z)
\label{eq1.4}
\end{equation}
We record here for later use the form of equations \ref{eq1.4} in terms
of the functions $\psi=\eta/(1-\eta)$:
\begin{equation}
\psi_1(\omega(z))=\psi_2(\tilde{\omega}(z))
=\frac{\omega(z)\tilde{\omega}(z)}{z-\omega(z)\tilde{\omega}(z)}
\label{eq1.5}
\end{equation}
In the following theorem we shall describe the boundary behavior of the 
functions $\omega$ and $\tilde{\omega}.$

\begin{theorem}\label{BoundaryX1}
Let $\eta_1,\eta_2,f$ and $\omega$ be as in 
Corollary \ref{subordinationwithomega}, and consider their restrictions to
$\mathbb C^+.$
Fix a point $a\in(0,+\infty)$ and define the following two self-maps of the 
upper
half-plane: $u_j(z)=a\eta_j(z)/z$, $j\in\{1,2\},$ $z\in\mathbb C^+.$ Assume that 
$u_2\circ u_1$ is not a conformal automorphism of $\mathbb C^+.$ Then:
\begin{enumerate}
\item[{\rm (1)}] If $C(\omega,a)\cap\mathbb C^+\neq\varnothing$, then the 
function $\omega$ extends analytically in a neighbourhood of $a$.
\item[{\rm (2)}] The limit $\lim_{z\to a}\omega(z)$ exists in $\mathbb C\cup
\{\infty\}.$
\end{enumerate}
\end{theorem}
\begin{proof}
We shall prove (1) by showing that the function $h$ extends analytically 
through $a$. Let $z_n$ be a sequence converging to $a$ such that 
$a\ell=\lim_{n\to\infty}\omega(z_n)$ belongs to the upper half-plane. We have:
$$\ell=\lim_{n\to\infty}h(z_n)=\lim_{n\to\infty}f(z_n,h(z_n))=
f(a,\ell)=f_a(\ell),$$
so $\ell$ is a fixed point of $f_a$ belonging to $\mathbb C^+.$
If $f_a$ is not a conformal automorphism of the upper half-plane, then 
the implicit function argument used in the proof of part (2) of Theorem
\ref{Secondsubordinationresult} implies that there
exists a neighbourhood of $a$ in $\mathbb C$ on which $h$, and hence $\omega$, 
extends analytically.

If $f_a(\cdot)$ is a conformal automorphism of the upper half-plane, denote
by $k$ its inverse. Consider $u_j,\ j=1,2,$ as in the statement of the theorem.
The equality
$$w=f_a(k(w))=h_2(ah_1(ak(w)))$$
is equivalent to $$u_2(u_1(ak(w)))=aw,\quad w\in\mathbb C^+.$$
We conclude that $u_1\circ k$ (and hence $u_1$) must be injective. 
Let $l(w)=ak(a^{-1}w),$ $w\in\mathbb C^+,$ so that $u_2(u_1(l(w)))=w.$ By 
applying $u_1\circ l$ to the above equality, we obtain
$$u_1(l(u_2(z)))=z,\ {\rm for\ all}\ z\in (u_1\circ l)(\mathbb C^+).$$
Since $(u_1\circ l)(\mathbb C^+)$ is open, the above equality must be true for
all $z\in\mathbb C^+,$ allowing us to conclude that $u_1$ is also surjective,
so a conformal automorphism. But then $u_2\circ u_1$ must also be a conformal
automorphism of $\mathbb C^+.$ 
This contradicts the hypothesis.

We now prove part (2) of the theorem.
Assume that $C(\omega,a)$ contains more than one point. 
By part (1), we must have $C(\omega,a)\subseteq\mathbb R\cup\{\infty\},$ and, 
by Lemma \ref{ClusterConex}, either $C(\omega,a)\setminus\{\infty\}$ is a 
closed interval in $\mathbb R$, or $\mathbb R\setminus C(\omega,a)$ is an open
interval in $\mathbb R.$

As we have shown in the proof of Theorem \ref{Boundary+} (2), for any
$c\in C(\omega,a)\setminus\{\infty\},$ with the possible exception of three
points, there exists a sequence $\{z_n^{(c)}\}_{n\in\mathbb N}$ converging to
$a$ such that $\lim_{n\to\infty}\omega(z_n^{(c)})=c$ and $\Re\omega(z^{(c)}_n)
=c$ for all $n$. By Theorem \ref{Fatou}, the limit $\lim_{n\to\infty}
\eta_1(\omega(z_n^{(c)}))$ exists for almost all $c\in C(\omega,a).$ Denote
it by $\eta_1(c).$ We claim that for any $c\in C(\omega,a),$ with at most 
three exceptions, $\eta_1(c)\not\in\mathbb C^+.$ Indeed, suppose that $\eta_1(
c)\in\mathbb C^+$ for some $c\in C(\omega,a), $ and assume that there 
exists a sequence $z_n^{(c)}$ such that $\omega(z_n^{(c)})$ converges 
nontangentially to $c$. Then, using (\ref{eq1.4}) and the fact that $\arg
\eta_{j}(z)>\arg z$ for all $z\in\mathbb C^+,$ $j\in\{1,2\},$ we obtain
\begin{eqnarray*}
\arg\eta_1(c) & = & \lim_{n\to\infty}\arg\eta_1(\omega(z_n^{(c)}))\\
& = & \lim_{n\to\infty}\arg\eta_2\left(\frac{z_n^{(c)}\eta_1(\omega(z_n^{(c)}))
}{\omega(z_n^{(c)})}\right)\\
& = & \arg\eta_2\left(\frac{a\eta_1(c)}{c}\right)\\
& > & \arg\left(\frac{a\eta_1(c)}{c}\right)\\
& = & \arg\eta_1(c),
\end{eqnarray*}
which is a contradiction.

Assume that $c_0\in{\rm Int}C(\omega,a)$ is a point where $\eta_1$ does not 
continue meromorphically. Proposition \ref{SeidelCaratheodory} (b) shows that 
the set
$$E=\{c\in C(\omega,a)\colon a\eta_1(c)/c\in(-\infty,0)\}$$
has nonzero Lebesgue measure.
In particular, for all $c\in E$, 
\begin{eqnarray*}
\eta_1(c)
& = & \lim_{n\to\infty}\eta_1(\omega(z_n^{(c)}))\\
& = &\lim_{n\to\infty}\eta_2\left(\frac{z_n^{(c)}\eta_1(\omega(z_n^{(c)})
)}{\omega(z_n^{(c)})}\right)\\
& = & \eta_2\left(\frac{a\eta_1(c)}{c}\right).
\end{eqnarray*}
By analytic continuation, this implies that $\eta_1(z)=\eta_2(a\eta_1(z)/z)$
for all $z\in\mathbb C^+.$
Rewritting this equality gives
$$a\eta_2\left(\frac{a\eta_1(z)}{z}\right)\frac{z}{a\eta_1(z)}=z,\quad z\in
\mathbb C^+,$$
or, equivalently,
$$u_2(u_1(z))=z,\quad z\in\mathbb C^+.$$ In particular, $u_2\circ u_1$ is
an automorphism of the upper half-plane, contradicting the hypothesis.

We conclude that $\eta_1$ must extend meromorphically through the set 
${\rm Int}C(\omega,a).$ Let $\tilde{\omega}(z)=z\eta_1(\omega(z))/\omega(z)$. 
We shall  argue that the set $C(\tilde{\omega},a)\subseteq\mathbb R\cup\{\infty
\}$ must also be infinte. Indeed, suppose this were not the case. Then, for any
$c\in{\rm Int}C(\omega,a)$, with at most one exception,
$$\lim_{z\to a}\tilde{\omega}(z)=\lim_{n\to\infty}\frac{z_n^{(c)}\eta_1(
\omega(z_n^{(c)}))}{\omega(z_n^{(c)})}=\frac{a\eta_1(c)}{c},$$
so, by analytic continuation,
$\eta_1(z)=z\lim_{z\to a}\tilde{\omega}(z)/a,$ $z\in\mathbb C^+.$ This 
contradicts the hypothesis that $\arg\eta_1(z)>\arg z$ for all $z$ in the 
upper half-plane. So $C(\tilde{\omega},a)$ is an infinite set. As before,
$\eta_2$ must continue meromorphically through all of ${\rm Int}C(\tilde{\omega
},a).$

By passing to a subsequence, we may assume that $\lim_{n\to\infty}\tilde{\omega
}(z_n^{(c)})$ exists for every $c\in{\rm Int}C(\omega,a)$, with at most one 
exception.
Suppose there were a point $d\in C(\tilde{\omega},a)$ and a set $V_{d}
\subset{\rm Int}C(\omega,a)$ of nonzero Lebesgue measure
such that $\lim_{n\to\infty}\tilde{\omega}(z_n^{(c)})=d$ for all $c\in V_{d}.$
Taking limit as $n\to\infty$ in the equality
$$z_n^{(c)}\eta_1(\omega(z_n^{(c)})=\omega(z_n^{(c)})\tilde{\omega}(z_n^{(c)})$$
gives 
$$a\eta_1(c)=cd,\quad {\rm for\ all}\ c\in V_d.$$
We conclude by Theorem \ref{Privalov} that $\eta_1(z)=(d/a)z$ for all $z\in\mathbb C^+$, and hence $\arg\eta_1(z)=\arg z$ for all $z\in\mathbb C^+,$ 
contradicting the hypothesis. Thus, there exists a set $E\subseteq C(\omega,a)$
of nonzero Lebesgue measure such that $\{\tilde{c}=\lim_{n\to\infty}\tilde{
\omega}(z_n^{(c)})\colon c\in E\}\subseteq{\rm Int}C(\tilde{\omega},a).$ 
Then, as $\eta_2$ extends analytically through ${\rm Int}C(\tilde{\omega},a)$, 
we have, by (\ref{eq1.4})
\begin{eqnarray*}
\eta_1(c) & = & \lim_{n\to\infty} \eta_1(\omega(z_n^{(c)}))\\
& = & \lim_{n\to\infty} \eta_2(\tilde{\omega}(z_n^{(c)}))\\
& = & \lim_{n\to\infty}\eta_2\left(\frac{z_n^{(c)}\eta_1(\omega(z_n^{(c)}))}{
\omega(z_n^{(c)})}\right)\\
& = & \eta_2\left(\frac{a\eta_1(c)}{c}\right)
\end{eqnarray*}
for all $z\in E.$ By analytic continuation, we conclude that
$$\eta_1(z)=\eta_2(a\eta_1(z)/z),\quad {\rm for\ all\ } z\in\mathbb C^+.$$ 
As we have already
proved, this implies that $u_1$ and $u_2$ are conformal automorphisms of 
the upper half-plane, contradicting the hypothesis. This concludes the proof. 
\end{proof}

The results above will be used for analyzing properties of free 
multiplicative convolution. As in Section 1, we start by stating several
well-known results.

In \cite{Voiculescu2}, Voiculescu introduces the $S-$transform, an analogue of the 
$R-$transform for free multiplicative convolution. We give below the main 
facts about the $S-$transform, and for proofs we refer to
\cite{Voiculescu2} and
\cite{BercoviciVoiculescuIUMJ}. 
Let $\mu$ be a probability measure on $[0,+\infty)$, $\mu\neq\delta_0.$ 
As observed in Remark \ref{Inj}, the restriction of $\eta_\mu$ to $(0,+\infty)$
is injective. Denote by $\eta_\mu^{-1}(z)$ its right inverse, defined
on some neighbourhood of $\eta_\mu((-\infty,0))$ in $\mathbb C^+.$ Let 
$\Sigma_\mu(z)=\eta_\mu^{-1}(z)/z.$
\begin{theorem}\label{Stransform}
Let $\mu,\nu$ be two probability measures supported on $[0,+\infty),$ none
of them concentrated at zero. Then $\Sigma_{\mu\boxtimes\nu}(z)=\Sigma_\mu(z)
\Sigma_\nu(z)$ for all $z$ in the common domain of the three functions.
\end{theorem}
The function $S_\mu(z)=\Sigma_\mu\left(\frac{z}{1+z}\right)$ is called the 
$S-$transform of $\mu$.

\bigskip

Let $\mu,\nu$ be two Borel probability measures on $[0,+\infty)$, none of them 
a point mass. Proposition \ref{eta} assures us that the functions $\eta_1
=\eta_\mu$ and $\eta_2=\eta_\nu$ satisfy the conditions of Corollary
\ref{subordinationwithomega}. With the notations from (\ref{eq1.4}), let
$\eta_3(z)=\eta_1(\omega(z))=\eta_2(\tilde{\omega}(z)).$ Observe that,
by Corollary \ref{subordinationwithomega}, $\eta_3$ satisfies the conditions
of Proposition \ref{eta}. By Remark \ref{Inj} and Proposition \ref{eta}, there 
exists an interval $(-s,0)$ such 
that $\eta_1,\eta_2,$ and $\eta_3$ have right inverses on $(-s,0).$
Let us replace $z$ in the relation
$$z\eta_3(z)=\omega(z)\tilde{\omega}(z),$$
(equation (\ref{eq1.4})) by $\eta_3^{-1}(z).$
We obtain 
$$z\eta_3^{-1}=\omega(\eta_3^{-1}(z))\tilde{\omega}(\eta_3^{-1}(z))=
\eta_1(z)\eta_2(z),$$
or, equivalently,
$$\frac{1}{z}\eta_3^{-1}(z)=\Sigma_\mu(z)\Sigma_\nu(z),\quad z\in(-s,0).$$
By Theorem \ref{Stransform}, we conclude that $\eta_3=\eta_{\mu\boxtimes\nu}.$

We can give now a new proof to Biane's subordination result (see \cite{Biane1},
Theorem 3.6):
\begin{corollary}\label{Biane2}
Given two probability measures $\mu,\nu$ on $[0,+\infty),$ there exists a 
pair of analytic functions $\omega_1,\omega_2\colon\mathbb C\setminus
[0,+\infty)\longrightarrow\mathbb C\setminus
[0,+\infty)$ such that
$\eta_\mu(\omega_1(z))=\eta_\nu(\omega_2(z))=\eta_{\mu\boxtimes\nu}(z)$, $
z\in\mathbb C\setminus
[0,+\infty).$ If neither of $\mu,\nu$ is concentrated at zero, the
 functions $\omega_1,\omega_2$ are unique and satisfy the following four
properties
\begin{enumerate}
\item[{\rm (a)}] $z\eta_{\mu\boxtimes\nu}(z)=\omega_1(z)\omega_2(z),$ $z\in
\mathbb C\setminus[0,+\infty).$
\item[{\rm (b)}] $\omega_j(0-)=0,$ $j\in\{1,2\},$
\item[{\rm (c)}] $\arg\omega_j(z)\in[\arg z,\pi)$ for all $z\in\mathbb C^+,$
$j\in\{1,2\},$
\item[{\rm (d)}] $\omega_j(\overline{z})=\overline{\omega_j(z)}$ for all $z\in
\mathbb C\setminus
[0,+\infty),$ $j\in\{1,2\}.$
\end{enumerate}
\end{corollary}

\begin{proof}
If one of the two measures, say $\mu$, equals $\delta_0,$ then 
$\mu\boxtimes\nu=\delta_0$,
and $\omega_1(z)=0,$ $\omega_2(z)=0,$ $z\in\mathbb C\setminus
[0,+\infty)$ will satisfy the required conditions. If one of the two measures,
say $\mu$, is concentrated at a point $a\in(0,+\infty),$ then the corollary is 
obviously true, with $\omega_1(z)=\eta_2(az)/a$ and $\omega_2(z)=az.$
If none of the two measures is a point mass, then, as we have already seen, if
we let $\eta_1=\eta_\mu$ and $\eta_2=\eta_\nu$, then the functions 
$\omega_1=\omega$ and $\omega_2=\tilde{\omega}$ provided by Corollary
\ref{subordinationwithomega} will satisfy the requirements from the statement
of the Corollary. 

Property (a) has been proved before for the functions $\omega,\tilde{\omega}$ 
(see relation (\ref{eq1.4})). To prove property (b),
observe that, by Proposition \ref{eta} and Remark \ref{Inj}, we have 
$$\omega_1(0-)=
\lim_{y\uparrow0}\omega_1(y)=\lim_{y\uparrow0}\eta_\mu^{-1}(\eta_{\mu\boxplus
\nu}(y))=\lim_{z\uparrow0}\eta_\mu^{-1}(z)=0.$$ Property (c) has been 
proved in Corollary \ref{subordinationwithomega}, and property (d) follows
trivially from the similar property of the functions $\eta.$
\end{proof}

Next, we use the boundary properties of the subordination functions
$\omega_1,\omega_2$ to describe the atomic and absolutely continuous parts with
respect to the Lebesgue measure of the free convolution
$\mu\boxtimes\nu$ of two probability measures on the positive half-line,
none of them  a point mass. 
          We start with some observations about the boundary behaviour of
the function $\psi_\mu$, particulary in the neighborhood of an
atom. Fix a probability measure $\mu$ on $[0,+\infty)$ and a real number $a>0$.
The following two
lemmas indicate how the presence or absence of an atom of $\mu$
can be detected from the boundary behaviour of $\psi_\mu$.

\begin{lemma}\label{atoms&tang}
 Suppose that $\{z_n=x_n+iy_n\}_{n=0}^\infty$
is a bounded sequence in $\mathbb{C}\setminus[0,+\infty)$ which converges
tangentially to $[0,+\infty)$, i.e.
$$\lim_{n\to\infty}\frac{y_n}{x_n-t}=0$$ for all $t\in[0,+\infty)$.
Then $\lim_{n\to\infty}y_n\psi_\mu(z_n)=0.$

\end{lemma}

\begin{proof}  Indeed, $$|y_n\psi_\mu(z_n)|\le
y_n\int_{[0,+\infty)}\left|\frac{t(x_n+iy_n)}{1-t(x_n+iy_n)}\right|\,d\mu(t)$$

$$\le\int_{[0,+\infty)}\frac{|y_nx_n|}{(y_n^2+(\frac{1}{t}-x_n)^2)^{1/2}}\,d\mu
(t)$$
$$+\int_{[0,+\infty)}\frac{y_n^2}{((\frac{1}{t}-x_n)^2+y_n^2)^{1/2}}\,d\mu(t)$$

The hypothesis implies that
$$\lim_{n\to\infty}\frac{y_nx_n}{((\frac{1}{t}-x_n)^2+y_n^2)^{1/2}}=\lim_{n\to\infty}\frac{y_n^2}{((\frac{1}{t}-x_n)^2+y_n^2)^{1/2}}=0$$
for all $t\in\mathbb{R}$.

 Moreover
$$\left|y_n\frac{x_n+iy_n}{((\frac{1}{t}-x_n)^2+y_n^2)^{1/2}}\right|=\left|y_n\frac{tz_n}{1-tz_n}\right|\le2\sup_n|z_n|,$$
and the desired result follows by dominated convergence.
\end{proof}

\begin{lemma}\label{atoms&nontang}
Suppose that
$\{z_n=x_n+iy_n\}_{n=0}^\infty$, is a bounded sequence in
$\mathbb{C}\setminus[0,+\infty)$, $a\in(0,+\infty)$ and
$z_n\to\frac{1}{a}$ nontangentially, i.e.
$$\inf_{n\in\mathbb{N}}\left|\frac{y_n}{x_n-(1/a)}\right|>0.$$ Then
$\lim_{n\to\infty}(1-az_n)\psi_\mu(z_n)=\mu(\{a\})$.
\end{lemma}

\begin{proof}  Fix $m>\sup_{n\in\mathbb{N}}\left|\frac{x_n-(1/a)}{y_n}\right|$, and observe that
$$\frac{|(1-z_na)z_nt|^2}{|1-z_nt|^2}=a^2\frac{((1/a)-x_n)^2+y_n^2}{((1/t)-x_n)^2+y_n^2}=$$
$$=a^2(x_n^2+y_n^2)\frac{1+\left(\frac{(1/a)-x_n}{y_n}\right)^2}{1+\left(\frac{(1/t)-x_n}{y_n}\right)^2}<\sup_{n\in\mathbb{N}}|z_n|(m^2+1)$$f\
\
or
$t\neq0$. The inequality remains trivially true for $t=0$ as well.
Now $\frac{(1-z_na)z_nt}{1-z_nt}\to0$ if $t\neq\frac{1}{a}$, and
$\frac{(1-z_na)z_nt}{1-z_nt}\to1$ if $t=1/a$. The conclusion of
the lemma follows then from the dominated convergence theorem.
\end{proof}

\begin{remark}\label{remarkaboutpsi}{\rm
Relation $G_\mu(1/z)=z(\psi_\mu(z)+1),$ $z\in\mathbb C^+$, implies that 
whenever $\psi_\mu|_{\mathbb C^+}$ extends analytically through a point $x\in
(0,+\infty)$, $G_\mu$ will extend analytically through $1/x$, and hence
the density of ${\mu}^{(ac)}$ with respect to the Lebesgue measure will 
be analytic in the point $1/x$.
}
\end{remark}

We take care first of a particular case, not covered by Theorem \ref{BoundaryX1}.
\begin{proposition}\label{SpecialcaseX1}
Let $\mu,\nu$ be two Borel probability measures supported on the interval
$[0,+\infty),$ none of them a point mass. Let $a\in(0,+\infty),$ and define
$u_\mu(z)=a\eta_\mu(z)/z,$ $u_\nu(z)=a\eta_\nu(z)/z,$ $z\in\mathbb C^+.$
The function $u_\mu\circ u_\nu$ is a conformal automorphism of the upper
half-plane if and only if $\mu$ and $\nu$ are convex combinations of two point 
masses. Moreover, in this case, $(\mu\boxtimes\nu)^{sc}=0$, and the density of 
$(\mu\boxtimes\nu)^{ac}$ with respect to the Lebesgue measure is analytic 
everywhere, with the exception of at most four points.
\end{proposition}
\begin{proof} Observe that $u_\mu$ and $u_\nu$ are nonconstant functions
taking values in the upper half-plane. 
We claim that both $u_\mu$ and $u_\nu$ are conformal automorphisms 
of $\mathbb C^+.$ Indeed, let $k$ be the inverse with respect to composition of
$u_\mu\circ u_\nu$, so that 
                 $$u_\mu(u_\nu(k(z)))=z,\quad z\in\mathbb C^+.$$
This implies that $u_\nu$ is injective. 
Applying $u_\nu$ to both sides of the above equality gives $u_\nu(u_\mu(w))=w$
for all $w$ in the open set $(u_\nu\circ k)(\mathbb C^+)$, so, by analytic 
continuation, for all $w\in\mathbb C^+.$ This proves that $u_\nu$ must be 
surjective, and thus an automorphism of $\mathbb C^+.$ Now it is obvious that 
$u_\mu$ must also be an automorphism of $\mathbb C^+$. This proves our claim.

As in the proof of Proposition \ref{Specialcase+}, there exist real
numbers  $p,q,r,s\in\mathbb R$ such that
$${\rm det}\left[\begin{array}{cc}
p & q \\
s & r
\end{array}\right]>0,\quad{\rm and}\quad u_\mu(z)=\frac{pz+q}{sz+r},\quad z
\in\mathbb C^+.$$
Since $zu_\mu(z)=a\eta_\mu(z)\in\mathbb C^+$ for all $z\in\mathbb C^+,$ we 
conclude that $s\neq0.$ Thus, we may assume without loss of generality 
that $s=1.$ Modulo a translation of $\mu$ and $\nu$, we may assume that 
$a=1$. By direct computation,
$$\eta_\mu(z)=zu_\mu(z)=\frac{pz^2+qz}{z+r},$$ so that 
   $$\psi_\mu(z)=\frac{z(pz+q)}{-pz^2-z(q-1)+r},\quad z\in\mathbb C^+.$$
Thus, $\mu$ is a convex combination of two atoms, at

$$\frac{q-1-\sqrt{(q-1)^2+4pr}}{2r}\quad{\rm and}\quad
\frac{q-1+\sqrt{(q-1)^2+4pr}}{2r},$$
with weights

$$\frac{1}{2}+\frac{q+1}{2r}\quad{\rm and}\quad\frac{1}{2}-\frac{q+1}{2r}.$$
A similar statement hold for $\nu.$

Conversely, if $\nu=t\delta_u+(1-t)\delta_v,$ then a direct computation shows 
that 
$$\psi_\nu(z)=z\frac{tu+(1-t)v-zuv}{1-z(u+v)+z^2uv},$$
so
$$\eta_\nu(z)=\frac{\psi_\nu(z)}{1+\psi_\nu(z)}=z\frac{tu+(1-t)v-zuv}{1-z((1-t)
u+tv)},$$ 
and thus
$$u_\nu(z)=\frac{-zuv+tu+(1-t)v}{-z((1-t)u+tv)+1},\quad z\in\mathbb C^+.$$
Observe that 
$$\det\left[\begin{array}{cc}
-uv & tu+(1-t)v\\
-(1-t)u-tv & 1
\end{array}\right]=t(1-t)(u-v)^2>0.$$
This proves the first statement of the proposition.

Assume now that $\mu=s\delta_\alpha+(1-s)\delta_u$ and $\nu=t\delta_\beta+
(1-t)\delta_v.$ Modulo a translation, we may assume that $\alpha=\beta=1$,
$u,v>0.$ We shall compute the function $\eta$ corresponding 
to the measure $(s\delta_1+(1-s)\delta_u)\boxtimes(t\delta_1+(1-t)\delta_v)$
by the same method as in the Proposition \ref{Specialcase+}. Unfortunately,
the expression of the function is much less manageable than in the 
additive case. 
As we have seen
before, 
$$\eta_\mu(z)=zu_\mu(z)=z\frac{(s+(1-s)u)-uz}{1-z(su+(1-s))},$$
and $$
\eta_\nu(z)=zu_\nu(z)=z\frac{(t+(1-t)v)-vz}{1-z(tv+(1-t))},$$
for all $z\in\mathbb C^+.$
From equation (\ref{eq1.4}), we conclude that the subordination function
$\omega_1(z)$ must satisfy
$$\eta_\mu(\omega_1(z))=\eta_\nu\left(\frac{z\eta_\mu(\omega_1(z))}{\omega_1(z)
}\right),$$ or, equivalently,

$${\frac {\omega_1(z) \left( s+ \left( 1-s \right) u \right) -u{\omega_1(z)}^{2}}{1-\omega_1(z)
 \left( su+1-s \right) }}=\ \ \ \ \ \ \ \ $$
$$
\frac{\left( {\frac {z \left( \omega_1(z) \left( s+ \left( 1-s \right) u \right) -u{\omega_1(z)
}^{2} \right)  \left( t+ \left( 1-t \right) v \right) }{\omega_1(z) \left( 1-\omega_1(z)
 \left( su+1-s \right)  \right) }}-{\frac {v{z}^{2} \left( \omega_1(z) \left( s+
 \left( 1-s \right) u \right) -u{\omega_1(z)}^{2} \right) ^{2}}{{\omega_1(z)}^{2} \left( 1
-\omega_1(z) \left( su+1-s \right)  \right) ^{2}}} \right)}{  \left( 1-{\frac {z
 \left( \omega_1(z) \left( s+ \left( 1-s \right) u \right) -u{\omega_1(z)}^{2} \right) 
 \left( tv+1-t \right) }{\omega_1(z) \left( 1-\omega_1(z) \left( su+1-s \right)  \right) }
} \right)}.
$$
With the help of Maple, one can obtain the solution 
\begin{eqnarray*}
\omega_1(z) & = & \frac{z^2uv+z((v-u)(s+t-1)+(t-s)(vu-1))-1}{2(zu(1-t-tv)-(1-
s+su))}\\
& & \mbox{}-\frac{\sqrt{Az^4+Bz^3+Cz^2+Dz+1}}{2(zu(1-t-tv)-(1-
s+su))},
\end{eqnarray*}
where the coefficients $A,B,C,$ and $D$ are given below:
$$A=u^2v^2,$$
$$B=2uv\left((v+u)((t-1)(s-1)-ts)+(1+uv)(t(s-1)+s(t-1))\right),$$
\begin{eqnarray}
C & = & (u^2+v^2)(t+s-1)^2\nonumber\\
& & \mbox{}+2(s(1-s)+t(1-t))(u+v-2uv+u^2v^2+uv^2)\nonumber\\
& & \mbox{}+4uv+(u^2v^2+1)(s-t)^2,\nonumber
\end{eqnarray}
$$D=-2\left((u+v)+(u-1)(v-1)(t(1-s)+s(1-t))\right).$$
From the expression of $\omega_1$ one can obtain $\omega_2$ by simply 
interchanging $(s,u)$ and $(t,v).$
By (\ref{eq1.4}), $\eta_{\mu\boxtimes\nu}(z)=\eta_\mu(\omega_1(z)),\ z\in
\mathbb C\setminus[0,+\infty),$ so 
$$\eta_{\mu\boxtimes\nu}(z)=\omega_1(z)\frac{(s+(1-s)u)-u\omega_1(z)}{1-
\omega_1(z)(su+1-s)},\quad z\in\mathbb C\setminus[0,+\infty).$$
Let us notice that, since $G_{\mu\boxtimes\nu}(1/z)=z(1-\eta_{\mu\boxtimes\nu}
(z))^{-1},$ the only points in $(0,+\infty)$ where the limit of
$G_{\mu\boxtimes\nu}$ is infinite are the points where the limit of $
\eta_{\mu\boxtimes\nu}$ equals 1, and $G_{\mu\boxtimes\nu}$ can be continued
analytically to $x\in[0,+\infty)$ if and only if $\eta_{\mu\boxtimes\nu}$ can 
be continued meromorphically to $1/x$, and $\eta_{\mu\boxtimes\nu}(1/x)\neq1.$
The only points where $\eta_{\mu\boxtimes\nu}$ may fail to extend 
meromorphically 
to $[0,+\infty)$ are the roots of the equation $Az^4+Bz^3+Cz^2+Dz-1=0$, while
$\eta_{\mu\boxtimes\nu}(1/x)=1$ if and only if 
$$\omega_1(1/x)\frac{(s+(1-s)u)-u\omega_1(1/x)}{1-\omega_1(1/x)(su+1-s)}=1/x,$$ i.e. 
if $\omega_1(1/x)$ satisfies an equation
of degree two with coefficients in $\mathbb R$; the algebraic expression of 
$\omega_1$ indicates clearly
that there can be at most eight points $x$ such that this happens. 
We conclude that there are at most twelve points in the interval $[0,+\infty)$ 
where $G_{\mu\boxtimes\nu}$ may fail to continue analytically. We conclude that
$(\mu\boxtimes\nu)^{sc}=0,$ and the density of $(\mu\boxtimes\nu)^{ac}$ is 
analytic everywhere, with at most twelve exceptions. This concludes the proof.
\end{proof}

\bigskip

We are now ready to prove the main regularity result for free multiplicative 
convolutions of probability measures on the positive half-line (parts (1) and
(2) of the theorem below appear also in \cite{B}).

\begin{theorem}\label{RegX1}
Let $\mu,\nu$ be two Borel probability measures supported on the interval
$[0,+\infty).$ Then:
\begin{enumerate}
    \item [\rm{(1)}] The following are equivalent:

    {\rm{(i)}} $\mu\boxtimes\nu$ has an atom at
    $a\in(0,+\infty)$;

    {\rm{(ii)}} there exist $u,v\in(0,+\infty)$ so that $uv=a$ and
$\mu(\{u\})+\nu(\{v\})>1$.
Moreover,
$\mu(\{u\})+\nu(\{v\})-1=(\mu\boxtimes\nu)(\{a\})$.
    \item [\rm{(2)}]$(\mu\boxtimes\nu)(\{0\})=\max\{\mu(\{0\}),
\nu(\{0\})\}.$
    \item[{\rm(3)}] Assume that neither $\mu,$ nor $\nu$ is a point mass. Then
$(\mu\boxtimes\nu)^{ac}$ is always nonzero, and its density with respect to 
the Lebesgue measure is analytic outside a closed set of Lebesgue measure 
zero. The support of $(\mu\boxtimes\nu)^{sc}$ is closed, of zero Lebesgue 
measure, and included in the support of $(\mu\boxtimes\nu)^{ac}$.
\end{enumerate}
\end{theorem}
\begin{proof}
We start by proving (1), (i)$\Rightarrow$(ii). So let
$\mu,\nu$ satisfy condition (i). Set
$z_n=\frac{1}{a}+i\frac{1}{n}$, $n\ge1$, and note that
$$\lim_{n\to\infty}(1-az_n)\psi_{\mu\boxtimes\nu}(z_n)=(\mu\boxtimes\nu)(\{a\})>0,$$
by Lemma \ref{atoms&nontang}. This implies that
$\lim_{n\to\infty}\psi_{\mu\boxtimes\nu}(z_n)=\infty$.
We shall prove the existence of two numbers $u,
v\in(0,+\infty)$ such that, after possibly dropping to a
subsequence,
\begin{trivlist}

\item[(j)]$uv=a$;
\item[(jj)]$\lim_{n\to\infty}\omega_\mu(z_n)=1/u$;
\item[(jjj)]$\lim_{n\to\infty}\omega_\nu(z_n)=1/v$;
\item[(jv)]the sequence $\Im\omega_\mu(z_n)(\Re\omega_\mu(z_n)-(1/u))^{-1}$ does not converge to
zero; and
\item[(v)]the sequence $\Im\omega_\nu(z_n)(\Re\omega_\nu(z_n)-(1/v))^{-1}$
does not converge to zero.

\end{trivlist}

First, we claim that the sequence $\{\omega_\mu(z_n)\}_{n=0}^\infty$ is
bounded. Assume to the contrary that there is a subsequence of
$\omega_\mu(z_n)$ which tends to infinity.
Observe that $\arg(\omega_\mu(z))\ge\arg(z)$ by Corollary \ref{Biane2},
so
$\frac{\Im\omega_\mu(z_n)}{\Re\omega_\mu(z_n)}\ge\frac{\Im(z_n)}{\Re(z_n)}=\frac{a}{n}$. Therefore
$$\left|(1-az_n)\frac{\omega_\mu(z_n)t}{1-\omega_\mu(z_n)t}\right|^2=\left(\frac{a}{n}\right)^2\frac{(\Re\omega_\mu(z_n))^2+(\Im\omega_\mu(z_n))^2}{(\frac{1}{t}-\Re\omega_\mu(z_n))^2+(\Im\omega_\mu(z_n))^2}$$
$$\le\left(\frac{a}{n}\right)^2\frac{1+\left(\frac{\Im\omega_\mu(z_n)}{\Re\omega_\mu(z_n)}\right)^2}{\left(\frac{\Im\omega_\mu(z_n)}{\Re\omega_\mu(z_n)}\right)^2}=\left(\frac{a}{n}\right)^2+\left(\frac{a\Re\omega_\mu(z_n)}{n\Im\omega_\mu(z_n)}\right)^2\le1+\left(\frac{a}{n}\right)^2,$$
which shows that the family of functions $\left\{\mathbb{R}_+\ni
t\to\-i\frac{a}{n}\frac{\omega_\mu(z_n)t}{1-\omega_\mu(z_n)t}\right\}_{n=0}^\infty$
is uniformly bounded and comverges pointwise to zero as
$n\to\infty.$
So we conclude that
$\lim_{n\to\infty}i\frac{a}{n}\psi_\mu(\omega_\mu(z_n))=0$. But by Lemma
\ref{atoms&tang}, we obtain

$$0<(\mu\boxtimes\nu)(\{a\})=\lim_{n\to\infty}(1-az_n)\psi_{\mu\boxtimes\nu}(z_n)=\lim_{n\to\infty}\left[-i\frac{a}{n}\psi_\mu(\omega_\mu(z_n))\right].$$
This contradicts the previous equality, and therefore
$\{\omega_\mu(z_n)\}_{n\in\mathbb N}$ must be bounded, as claimed. Passing
to a subsequence, we infer the existence of $u\in\mathbb{C}$ such that
$\omega_\mu(z_n)\to\frac{1}{u} $ as $n\to\infty$.

Setting $v=a/u$, we have
$$\lim_{n\to\infty}\omega_\nu(z_n)=\frac{1}{a\lim_{n\to\infty}\omega_\mu(z_n)}=\frac{1}{a/u}=\frac{1}{v}$$
by (\ref{eq1.5}) and Lemma \ref{atoms&nontang}. Note that $u, 
v\in(0,+\infty)$; indeed, if $u\not\in(0,+\infty)$, we would have
$\psi_{\mu\boxtimes\nu}(z_n)=\psi_\mu(\omega_\mu(z_n))\to\psi_\mu(u)$, which is
a finite number. This proves (j), (jj) and (jjj).

 To prove (jv), note that
$$\lim_{n\to\infty}|\Im\omega_\mu(z_n)\psi_\mu(\omega_\mu(z_n))|=\lim_{n\to\infty}\left|\frac{\Im\omega_\mu(z_n)}{(-ia/n)}(1-az_n)\psi_{\mu\boxtimes\nu}(z_n)
\right|$$
$$=(\mu\boxtimes\nu)(\{a\})\lim_{n\to\infty}\left|\frac{\Im\omega_\mu(z_n)}{(-ia/n)}\right|\ge(\mu\boxtimes\nu)(\{a\})\lim_{n\to\infty}|\Re\omega_\mu(z_n)|$$
$$=(\mu\boxtimes\nu)(\{a\})u>0,$$ by Corollary \ref{Biane2}
and Lemma \ref{atoms&nontang}.
Also, the sequence
$$\left\{\left|\frac{\Im\omega_\mu(z_n)}{\Re\omega_\mu(z_n)-(1/u)}\right|\right\}_{n=0}^\infty$$
does not converge to zero, by Lemma \ref{atoms&tang}. Dropping if necessary to
a subsequence, we can assume that it is bounded away from zero.
This proves (jv). Statement (v) follows by symmetry.

Conditions (jj) and (jv) allow us to apply Lemma \ref{atoms&nontang} to 
conclude that
$$\lim_{n\to\infty}(1-\omega_\mu(z_n)u)\psi_\mu(\omega_\mu(z_n))=\mu(\{u\})$$
and analogously
$$\lim_{n\to\infty}(1-\omega_\nu(z_n)v)\psi_\nu(\omega_\nu(z_n))=\nu(\{v\}).$$
We conclude that
   $$\mu(\{u\})+\nu(\{v\})=$$
   $$=\lim_{n\to\infty}[(1-\omega_\mu(z_n)u)\psi_\mu(\omega_\mu(z_n))+(1-\omega_\nu(z_n)v)\psi_\nu(\omega_\nu(z_n))]$$
   $$=\lim_{n\to\infty}[(1-\omega_\mu(z_n)u)\psi_\mu(\omega_\mu(z_n))+\omega_\mu(z_n)u(1-\omega_\nu(z_n)v)\psi_\nu(\omega_\nu(z_n))]$$
   $$=\lim_{n\to\infty}[(1-\omega_\mu(z_n)\omega_\nu(z_n)uv)\psi_{\mu\boxtimes\nu}(z_n)]$$
   $$=\lim_{n\to\infty}[(1-az_n)\psi_{\mu\boxtimes\nu}(z_n)+a(z_n-\omega_\mu(z_n)\omega_\nu(z_n))\psi_{\mu\boxtimes\nu}(z_n)]$$
   $$=\lim_{n\to\infty}[(1-az_n)\psi_{\mu\boxtimes\nu}(z_n)]+\lim_{n\to\infty}[a\omega_\mu(z_n)\omega_\nu(z_n)]$$
$$=(\mu\boxtimes\nu)(\{a\})+a\overline{a}=1+(\mu\boxtimes\nu)(\{a\}),$$
   where we used (\ref{eq1.5}) in the next to the last equality.
This
completes the proof of (i)$\Rightarrow$(ii).

 To prove
(ii)$\Rightarrow$(i), consider two free selfadjoint random
variables $X$ and $Y$ in the tracial noncommutative probability
space $({\mathcal A}, \tau)$ having distributions $\mu$ and $\nu$
respectively. The fact that $\mu$ has an atom at $u$, $\nu$ has an
atom at $v$ and $\mu(\{u\})+\nu(\{v\})>1$ is equivalent to the
fact that there are projections $p$ and $q$ with $\tau(p)+\tau(
q)>1$ such that $Xp=up$ and $Yq=vq$. Then $p\wedge q\neq0$ and
$X^{1/2}YX^{1/2}(p\wedge q)=uv(p\wedge q)$, so $\mu\boxtimes\nu$ has an
atom at $uv$. This concludes the proof of part (1).

To prove (2),  let $X$ and $Y$ be two selfadjoint random variables
as in the proof of (a). Note that for any projection $p$ such that
$Xp=0$, we also have that $X^{1/2}YX^{1/2}p=0$, and we deduce that
$(\mu\boxtimes\nu)(\{0\})\geq\mu(\{0\}).$ Since $\boxtimes$ is commutative, we
have $(\mu\boxtimes\nu)(\{0\})\geq\max\{\mu(\{0\}), \nu(\{0\})\}.$

We first prove the opposite inequality for $(\mu\boxtimes\nu)(\{0\})<1$.
To do this, we analyze the behaviour of $\psi_\mu(z)$ as
$z\to\infty$ nontangentially to $[0,+\infty)$. Consider the truncated cone
$\Gamma=\{x+iy:y>0, x\le y, |x|+|y|\ge 1\}$ and note that for $z=x+iy\in\Gamma$
we have

$$\left|\frac{zt}{1-zt}\right|^2=\left|\frac{(x+iy)t}{1-(x+iy)t}\right|^2<3,
\ \ t\in\mathbb{R}_+.$$
   Since $\lim_{z\to\infty, z\in\Gamma}\frac{zt}{1-zt}$ equals
   $-1$ if $t\neq 0$ and $0$ if $t=0$, the dominated convergence
   theorem yields $\lim_{z\to\infty,
   z\in\Gamma}\psi_\mu(z)=\mu(\{0\})-1.$
   Applying this fact to $\mu\boxtimes\nu$ combined with (\ref{eq1.5}), we 
obtain
$$\lim_{z\to\infty,
z\in\Gamma}\frac{\omega_\mu(z)\omega_\nu(z)}{\omega_\mu(z)\omega_\nu(z)-z}=1-(\mu\boxtimes\nu)(\{0\}).$$
We deduce that $\lim_{z\to\infty}\omega_\mu(z)\omega_\nu(z)=\infty$ and
therefore there is a sequence $z_n\in\Gamma$ such that
$z_n\to\infty$ and either $\omega_\mu(z_n)\to\infty$, or
$\omega_\nu(z_n)\to\infty$.
Assume for simplicity that
$\omega_\mu(z_n)\to\infty$, and note that $\omega_\mu(z_n)\in\Gamma$ by
Corollary \ref{Biane2}. By the preceding calculations, we must 
have
 $$(\mu\boxtimes\nu)(\{0\})-1=\lim_{n\to\infty}\psi_{\mu\boxtimes\nu}(z_n)=\lim_{n\to\infty}\psi_\mu(\omega_\mu(z_n))=\mu(\{0\})-1.$$
 Thus $(\mu\boxtimes\nu)(\{0\})\leq\max\{\mu(\{0\}), \nu(\{0\})\}$ whenever
 $(\mu\boxtimes\nu)(\{0\})<1$.
To conclude the proof, we must show that $\mu(\{0\})<1$ and
$\nu(\{0\})<1$ implies $(\mu\boxtimes\nu)(\{0\})<1$. Indeed, in this case
Corollary \ref{Biane2} applies to show that
$\psi_{\mu\boxtimes\nu}=\psi_\mu\circ\omega_\mu$ is not identically zero.

The proof of (3) is similar to the proof of part (2) of Theorem \ref{Reg+}.
Suppose that $(\mu\boxtimes\nu)^{ac}=0.$
Part (1) of the theorem shows that $\mu\boxtimes\nu$ cannot be purely atomic,
and thus, $(\mu\boxtimes\nu)^{sc}\neq0.$ Since $\psi_{\mu\boxtimes\nu}=
(1/z)G_{\mu\boxtimes\nu}(1/z)-1,$ we conclude as in the proof of part (2) of
Theorem \ref{Reg+} that $\lim_{y\downarrow0}\Im\psi_{\mu\boxtimes\nu}(x+iy)=0$
for almost all $x\in\mathbb R$, and that there exists a point $x_0\in(0,+\infty
)$ such that $\psi_{\mu\boxtimes\nu}$ does not continue analytically through
$x_0$. Theorem \ref{Seidel} can now be applied to the function $
\psi_{\mu\boxtimes\nu}$ and the point $x_0$ to conclude that $C(\psi_{\mu
\boxtimes\nu},x_0)=\overline{\mathbb C^+},$ and hence $C(\eta_{\mu
\boxtimes\nu},x_0)=\overline{\mathbb C^+}.$ 
Relation (\ref{eq1.4}) implies that $C(\omega_j,x_0)\cap\mathbb C^+$, $j\in\{1,
2\}$ are also infinite sets (we recall that $\omega_1,\ \omega_2$ are the 
subordination functions provided by Corollary \ref{Biane2}).
Applying Theorem \ref{BoundaryX1} and Proposition \ref{SpecialcaseX1} yields
a contradiction. Thus, $\mu\boxtimes\nu$ cannot be purely singular.

Next we show that there exists a closed subset of $[0,+\infty)$ of zero 
Lebesgue measure on whose complement the density $f(x)=\frac{d
(\mu\boxtimes\nu)^{ac}}{dx}$ is analytic. By Theorem \ref{Fatou} and 
Proposition \ref{Moment}, there 
exists a set $E\subset[0,+\infty)$ of zero Lebesgue measure such that for
all $x\in\mathbb R\setminus E$ the limits $\lim_{y\to0}\psi_{\mu\boxtimes\nu}
(x+iy),$ $\lim_{y\to0}\omega_{j}(x+iy)$, $j\in\{1,2\}$ exist and are finite.
Also, for almost all $x\in(0,+\infty)$ such that $(1/x)\in{\rm supp}((\mu
\boxtimes\nu)^{ac})\setminus E$, with respect to $(\mu\boxtimes\nu)^{ac}$, we 
have $\lim_{y\downarrow0}
G_{\mu\boxtimes\nu}\left(\frac{1}{x+iy}\right)\in\mathbb C^+,$ so that 
$$\lim_{y\downarrow
0}\psi_{\mu\boxtimes\nu}(x+iy)=\left(\frac{1}{x}\lim_{y\downarrow0}G_{\mu
\boxtimes\nu}\left(\frac{1}{x+iy}\right)-1\right)\in\mathbb C^+.$$
By equation (\ref{eq1.5}), at least one of $\lim_{y\to0}\omega_{j}(x+iy)$, $j
\in\{1,2\}$ must also be in $\mathbb C^+.$ For definiteness, assume that 
$\omega_1(x)=\lim_{y\to0}\omega_{1}(x+iy)\in\mathbb C^+.$ Part (1) of Theorem 
\ref{BoundaryX1} and Proposition \ref{SpecialcaseX1} assure us that $\omega_1
$ extends analytically through $x$. Since $\omega_1(x)\in\mathbb C^+,$
the function $\psi_{\mu\boxtimes\nu}=\psi_\mu\circ\omega_1$ will also 
extend analytically through $x$. By Remark \ref{remarkaboutpsi}, we conclude
that the density of $\mu\boxtimes\nu$ with respect to the Lebesgue measure must
be analytic in $1/x$. The set $A$ of points $x\in{\rm supp}((\mu\boxtimes\nu
)^{ac})$ where $f(x)=\frac{d(\mu\boxtimes\nu)^{ac}(x)}{dx}$ is analytic 
is open in $\mathbb R$, and its complement in ${\rm supp}((\mu\boxplus\nu)^{ac})$ 
is of zero Lebesgue measure. On the other hand, if there exists an interval 
$J\subseteq[0,+\infty)\setminus{\rm supp}((\mu\boxplus\nu)^{ac})$, then 
$\lim_{y\downarrow0}\Im\psi_{\mu\boxtimes\nu}(x+iy)=0$ for 
almost all $x\in\mathbb R$ such that $1/x\in J$, and the argument used in the
proof of the existence of the absolutely continuous part of $\mu\boxtimes\nu$
shows that $\psi_{\mu\boxtimes\nu}$ extends meromorphically through $\{x\colon
1/x\in J\}.$
In particular, $J\cap{\rm supp}((\mu\boxtimes\nu)^{sc})=\varnothing.$
\end{proof}


\section[]{Regularity for the free multiplicative convolution of probability 
measures on the unit circle}

The tools and methods from this section are very similar to the ones used in 
the previous section. Let $\mu$ be a Borel probability measure on the unit
circle $\mathbb T$, and define the moment generating function
$$\psi_\sigma(z)=\int_{\mathbb T}\frac{zt}{1-zt}d\sigma(t),
\quad z\in\mathbb D.$$ 
(Recall that $\mathbb D$ denotes the unit disk.) Observe that
$$\psi_\mu(z)=\int_0^{2\pi}\frac z{e^{it}-z}\,d\mu(e^{-it})=-\frac12+\frac12
\int_0^{2\pi}\frac {e^{it}+z}{e^{it}-z}\,d\mu(e^{-it}),\quad z\in\mathbb D.$$
This immediately yields the following
result.

\begin{proposition}\label{Momentcircle}
Let $\psi:\mathbb D\to\mathbb C$ be an analytic function. The following
conditions
are equivalent.
\begin{enumerate}
\item[{\rm (1)}] There exists a probability measure $\mu$ on $\mathbb T$ such
that
$\psi=\psi_\mu$.
\item[{\rm (2)}] $\psi(0)=0$ and $\Re\psi(z)\ge -1/2$ for all $z\in\mathbb D$.
\end{enumerate}
\end{proposition}

The corresponding result for the function $\eta_\mu=\psi_\mu/(1+\psi_\mu)$ is as
follows.

\begin{proposition}\label{etacircle}
Let $\eta:\mathbb D\to\mathbb C$ be an analytic function. The following
conditions
are equivalent.
\begin{enumerate}
\item[{\rm (1)}] There exists a probability measure $\mu$ on $\mathbb T$ such
that
$\eta=\eta_\mu$.
\item[{\rm (2)}] $\eta(0)=0$ and $|\eta(z)|<1$ for all $z\in\mathbb D$.
\item[{\rm (3)}] $|\eta(z)|\le|z|$ for all $z\in\mathbb D$.
\end{enumerate}
Moreover, $\mu$ is a point mass if and only if there exists $z_0\in\mathbb D
\setminus\{0\}$ such that $|\eta_\mu(z_0)|=|z_0|.$
\end{proposition}
\begin{proof}
The equivalence of (2) and (3) is, of course, the Schwarz lemma, while the 
equivalence between (1) and (2) is a trivial consequence of Proposition
\ref{Momentcircle}. 
Suppose now that there exists $z_0\in\mathbb D\setminus\{0\}$ such 
that $|\eta_\mu(z_0)|=|z_0|.$ Then $\eta_\mu(z_0)=e^{i\theta}z_0$ for some
$\theta\in[0,2\pi).$ The function $g(z)=\eta_\mu(z)/z$ maps $\mathbb D$ into
$\overline{\mathbb D}$, so, since $|g(z_0)|=1,$ by the maximum principle we 
conclude that $g$ must be constant, and thus $\eta_\mu(z)=e^{i\theta}z$ for 
all $z\in\mathbb D$. Thus, $\psi_\mu(z)=ze^{i\theta}(1-ze^{i\theta})^{-1},$ 
so that $\mu=\delta_{e^{i\theta}}.$ The converse is trivial.
\end{proof}

Next we prove the analogue of Theorem \ref{Secondsubordinationresult} for
self-maps of the unit disk which fix the origin.

\begin{theorem}\label{Thirdsubordinationresult}
Let $\eta_j\colon\mathbb D\longrightarrow\mathbb D$, $j\in\{1,2\}$, be two
analytic functions, none of them a conformal automorphism of the unit
disk, satisfying $\eta_j(0)=0,\ \eta_j'(0)\neq0$, $j\in\{1,2\}.$ Define
$f\colon(\mathbb D\cup\mathbb T)\times\mathbb D\longrightarrow\mathbb D$ by
$$f(z,w)=\frac{w\eta_2\left(z\frac{\eta_1(w)}{w}\right)}{\eta_1(w)},$$
where we understand $\frac{\eta(0)}{0}$ as $\eta'(0).$
Then
\begin{enumerate}
\item[{\rm (1)}] The function $f$ is analytic on $\mathbb D\times\mathbb D$
and takes values in $\mathbb D.$
\item[{\rm (2)}] There exists a unique analytic function
$\omega\colon\mathbb D\longrightarrow\mathbb D$ such that
$$f(z,\omega(z))=\omega(z),\quad z\in\mathbb D.$$ Moreover,
$\omega(0)=0.$
\end{enumerate}
\end{theorem}
\begin{proof}
It is obvious that
$$\lim_{z\to z_0,w\to 0}f(z,w)=\frac{\eta_2(z_0\eta_1'(0))}{\eta_1'(0)},$$
for any $z_0\in\mathbb D$. The only problem that might occur is when there
exists a $w_0\neq0$ such that $\eta_1(w_0)=0$. But then
$$\lim_{z\to z_0,w\to w_0}f(z,w)=\lim_{z\to z_0,w\to w_0}z\frac{\eta_2\left(z\frac{
\eta_1(w)}{w}\right)}{z\frac{\eta_1(w)}{w}}=z_0\eta_2'(0).$$
So $f$ is a well-defined complex-valued analytic function. To see that it
takes values in the unit disc, we apply the Schwarz lemma to conclude
that for any $z\in\mathbb D$, $w\in\mathbb D\setminus\{0\}$ such that $\eta_1(
w)\neq0$ we have
$$|f(z,w)|=\left|\frac{\eta_2\left(z\frac{\eta_1(w)}{w}\right)}{\frac{
\eta_1(w)}{w}}\right|\leq \left|z\frac{z\frac{\eta_1(w)}{w}}{z
\frac{\eta_1(w)}{w}}\right|=|z|<1.$$ This inequality extends by continuity to
all points $z,w\in\mathbb D$. This proves (1).

We now prove (2).
Fix a point $z\in\mathbb D\setminus\{0\}.$ As in the proof of Theorem
\ref{Firstsubordinationresult},
the function $f_z(w)=f(z,w)$ is a self-map of the unit disc, satisfying, as
shown above, the inequality $|f_z(w)|\leq|z|<1.$ Thus, according to Theorem
\ref{DenjoyWolff}, it must have a (unique) Denjoy-Wolff point. Denote it by
$\omega(z)$. Since (by the same Theorem \ref{DenjoyWolff}), we have $|f_z'
(\omega(z))|<1,$ we can apply the Implicit Function Theorem to conclude that 
the dependence $z\mapsto\omega(z)$ is analytic. Obviously, 
$\lim_{z\to0}\omega(z)=1.$ This proves (2) and concludes the proof.
\end{proof}

As before, we observe that there exists an analytic function
$\tilde{\omega}\colon\mathbb D\longrightarrow\mathbb D$, 
$\tilde{\omega}(z)=z\eta_1(\omega(z))/\omega(z),$ such that 
$$\tilde{\omega}(z)=\frac{\tilde{\omega}(z)\eta_1\left(z\frac{\eta_2(
\tilde{\omega}(z))}{\tilde{\omega}(z)}
\right)}{\eta_2(\tilde{\omega}(z))},\quad z\in\mathbb D.$$
The pair of functions $(\omega,\tilde{\omega})$ is uniquely
determined by the equations
\begin{equation}
\eta_1(\omega(z))=\eta_2(\tilde{\omega}(z))
=\frac{1}{z}\omega(z)\tilde{\omega}(z)
\label{eq1.6}
\end{equation}
As in Section 2, equatin (\ref{eq1.6}) can be expressed in terms
of the functions $\psi=\eta/(1-\eta)$:
\begin{equation}
\psi_1(\omega(z))=\psi_2(\tilde{\omega}(z))
=\frac{\omega(z)\tilde{\omega}(z)}{z-\omega(z)\tilde{\omega}(z)}
\label{eq1.7}
\end{equation}
\begin{remark}\label{remarkabour1stmoment0}{\rm
It is obvious that the condition $\eta_j'(0)\neq0$, $j\in\{1,2\}$ cannot be 
removed, since in that case we would (by Theorem \ref{DenjoyWolff}) have 
$\omega(z)=\tilde{\omega}(z)=0$ for all $z\in\mathbb D.$
However, one can still consider the case when only one of the two functions
$\eta,$ say $\eta_1$, has its derivative at the origin equal to zero.
Indeed, in this case the proof of Theorem \ref{Thirdsubordinationresult} 
applies to provide a function $\omega$ with the property that 
$\omega(z)$ is the fixed point of $f_z.$ In this case we can still
define $\tilde{\omega}$ by $\tilde{\omega}(z)=z\eta_1(\omega(z))/\omega(z),$ 
but $\tilde{\omega}(z)$ will not be necessarily the fixed point
of $g_z=g(z,\cdot),$ where  
$$ g(z,w)=\frac{w\eta_1\left(z\frac{\eta_2(w)}{w}\right)}{\eta_2(w)};$$
in fact, if there exists a point $w_0\neq0$ in the unit disk such that 
$\eta_2(w_0)=0$, then $g_z(0)=0$ for all $z\in\mathbb D$. 
}
\end{remark}

In the following theorem we describe some boundary properties of the 
functions $\omega$ and $\tilde{\omega}.$ Since the proof is similar 
to the proofs of Theorems \ref{Boundary+} and \ref{BoundaryX1}, we
omit it.
\begin{theorem}\label{BoundaryX2}
Let $\eta_1,\eta_2,f$ and $\omega$ be as in
Theorem \ref{Thirdsubordinationresult}.
Fix a point $a\in\mathbb T$ and define the following two self-maps of the
upper
half-plane: $u_j(z)=a\eta_j(z)/z$, $j\in\{1,2\},$ $z\in\mathbb D.$ Assume that
$u_2\circ u_1$ is not a conformal automorphism of $\mathbb D.$ Then:
\begin{enumerate}
\item[{\rm (1)}] If $C(\omega,a)\cap\mathbb D\neq\varnothing$, then the
function $\omega$ extends analytically in a neighbourhood of $a$.
\item[{\rm (2)}]  Assume that there exist open intervals $I_1$, $I_2$ such that
$\eta_j$ continues analytically through $\{e^{i\theta}\colon\theta\in I_j\},$ 
$j\in\{1,2\}.$
Then the limit $\lim_{z\to a}\omega(z)$ exists in $\overline{\mathbb D}.$
\end{enumerate}

\end{theorem}

We shall apply the previous results to the study of free convolutions of 
probability measures on the unit circle (for details, see \cite{Voiculescu2}). 
The definition and properties of the 
$S-$transform for measures on $\mathbb T$ are similar to the ones from 
Section 2. However,
in order to
invert the
function $\eta_\mu$ around zero one needs the condition $\eta_\mu'(0)\ne0$.
This is equivalent
to $\int_0^{2\pi}e^{it}\,d\mu(e^{it})\ne0$ because
$$\eta'_\mu(0)=\psi'_\mu(0)=\int_0^{2\pi}e^{it}\,d\mu(e^{it}).$$
Under this condition, the inverse $\eta_\mu^{-1}$ is defined in a neighborhood
of zero, and so is $\Sigma_\mu(z)=\eta^{-1}_\mu(z)/z$. The identity
$$\Sigma_{\mu\boxtimes\nu}(z)=\Sigma_\mu(z)\Sigma_\nu(z)$$
holds in a neighborhood of zero provided that $\eta_\mu'(0)\ne0\ne\eta'_\nu(0)$.
We can now prove a subordination theorem for free multiplicative convolution 
on the unit circle (see \cite{Biane1}, Theorem 3.5).
\begin{corollary}\label{Biane3}
Given two Borel probability measures $\mu,\nu$ on the unit circle, there exists
a pair of analytic functions $\omega_1,\omega_2\colon\mathbb D\longrightarrow
\mathbb D$ such that $\eta_\mu(\omega_1(z))=\eta_\nu(\omega_2(z))=
\eta_{\mu\boxtimes\nu}(z),$ $z\in\mathbb D.$
The pair of functions $\omega_1,\omega_2$ is uniquely determined by
the following properties:
\begin{enumerate}
\item[{\rm(a)}]  $z\eta_{\mu\boxtimes\nu}(z)=\omega_1(z)\omega_2(z),$ $z\in
\mathbb D;$
\item[{\rm (b)}] $\omega_j(0)=0,$ $j\in\{1,2\};$
\item[{\rm (c)}] $|\omega_j(z)|\leq|z|$ for all $z\in\mathbb D,$
$j\in\{1,2\}.$
\end{enumerate}
\end{corollary}
\begin{proof}
If both $\mu$ and $\nu$ have first moment equal to zero, then we 
conclude by Definitions \ref{free} and \ref{Free} that $\mu\boxtimes\nu$
has all moments equal to zero, so $\mu\boxtimes\nu$ is the uniform 
distribution (the Haar measure) on $\mathbb T.$ We conclude that
$\eta_{\mu\boxtimes\nu}(z)=0$ for all $z\in\mathbb D,$ so the 
functions $\omega_1(z)=\omega_2(z)=0$, $z\in\mathbb D$ will do. 
If only one of the two measures has first moment zero, then, choosing in Remark
\ref{remarkabour1stmoment0}, $\omega_1=\omega$ and $\omega_2=\tilde{\omega}$ will
provide the required pair. Observe that in this case the function $\omega_1$
is uniquely determined.
If both $\mu$ and $\nu$ have nonzero first moment, then choosing $\omega_1=\omega$
in Theorem \ref{Thirdsubordinationresult} and $\omega_2=\tilde{\omega}$ from 
equation (\ref{eq1.6}) will do. The uniqueness is guaranteed by Theorem
\ref{Thirdsubordinationresult}, and the properties (a) through (c) are 
elementary consequences of Theorem
\ref{Thirdsubordinationresult}.
\end{proof}

For measures on $\mathbb T$ we shall consider absolute continuity with respect 
to the uniform distribution on $\mathbb T$. The notations will be the same
as in Section 1.
The following regularity result holds for free convolutions of probability measures
supported on $\mathbb T$. The proof (see also \cite{B}) is similar to the proofs 
of Theorems \ref{Reg+} and \ref{RegX1} and will be omitted.
\begin{theorem}\label{RegX2}
Let $\mu,\nu$ be two Borel probability measures supported on $\mathbb T.$ Then:
\begin{enumerate}
    \item [\rm{(1)}] The following are equivalent:

    {\rm{(i)}} $\mu\boxtimes\nu$ has an atom at
    $a\in(0,+\infty)$;

    {\rm{(ii)}} there exist $u,v\in(0,+\infty)$ so that $uv=a$ and
$\mu(\{u\})+\nu(\{v\})>1$.
Moreover,
$\mu(\{u\})+\nu(\{v\})-1=(\mu\boxtimes\nu)(\{a\})$.
    \item[{\rm(2)}] Assume that neither $\mu,$ nor $\nu$ is a point mass. Then
$(\mu\boxtimes\nu)^{ac}$ is always nonzero, and its density with respect to
the Lebesgue measure is analytic outside a closed set of Lebesgue measure
zero. The support of $(\mu\boxtimes\nu)^{sc}$ is closed, of zero Lebesgue
measure, and included in the support of $(\mu\boxtimes\nu)^{ac}$.
\end{enumerate}
\end{theorem}


\chapter[Semigroups relative to free convolutions]{Semigroups relative to free 
convolutions}

In this chapter we study from an analytic point of view the existence and
some properties of partially defined free convolution semigroups.
The existence of these objects,
which have no classical counterpart, has been first
observed in \cite{BVSemigroup}: it is shown there that for any given compactly
supported probability measure $\mu$ on the real line, 
there exists a number $T\geq1$ and a family of compactly supported 
probability measures $\{\mu_t\colon t\geq T\}$ on $\mathbb R$ such that 
$\mu_{t+s}=\mu_t\boxplus\mu_s$ for all $s,t\geq T$ and $\mu_t=\underbrace{\mu
\boxplus\cdots\boxplus\mu}_{t\ \rm times}$ for $t\in\mathbb N,$ $t\geq T.$
Using combinatorial tools, Nica and Speicher have proved in \cite{NS} that all
the above statements are true for $T=1$. 

An analytic approach to the existence of the family $\{\mu_t\colon t\geq T\}$ 
can be found in \cite{BB}. As a consequence of results in Chapter 1 we know
that for any $t\in\mathbb N,$ there exists an analytic self-map $\omega_t$ of 
$\mathbb C^+$ such that $G_{\mu_t}=G_\mu\circ\omega_t$. In \cite{BB} it is 
shown that this result extends to noninteger values of $t>1$. The proof
of the existence of the subordination function turns out to provide a 
sufficient argument for the existence of the measures $\mu_t$ for all $t>1.$
Several regularity results are deduced from properties of the function 
$\omega_t$.

The approach presented in this chapter is based on Theorem \ref{DenjoyWolff}.
We show existence and regularity properties for measures in partial 
semigroups with respect to both additive and multiplicative convolutions. 
We also investigate connections with infinite divisibility with respect
to free and boolean convoluitons, and point out some new properties of the
Cauchy transforms of measures which are infinitely divisible with respect to
free convolutions. Most of the results presented here appear in 
\cite{BBercoviciMult}.

\section[]{Partial semigroups with respect to free additive convolution}

We present first a global inversion theorem for certain analytic 
functions defined in the upper half-plane. A somewhat weaker result has been 
proved in \cite{BB}, Theorem 2.2.

\begin{theorem}\label{InversionC+}
Let $H:\mathbb C^+\to\mathbb C$ be an analytic function such that $\Im
H(z)\le\Im z$ for $z\in\mathbb C^+$,
and the limit $$a=\lim_{y\to+\infty}\frac{H(iy)}{iy}$$
is strictly positive.\begin{enumerate}
\item[{\rm (1)}] For every $\alpha\in\mathbb C^+$ there exists a unique
$z\in\mathbb C^+$ such that
$H(z)=\alpha$.
\item[{\rm (2)}] There exists a continuous function $\omega:\mathbb
C^+\cup\mathbb R\longrightarrow\mathbb C^+
\cup\mathbb R$ such that
$\omega(\mathbb C^+)\subset\mathbb C^+$, $\omega|_\mathbb C^+$ is analytic, and
$H(\omega(z))=z$ for $z\in\mathbb C^+$.
\item[{\rm (3)}] For each $\alpha\in\mathbb C^+\cup\mathbb R$, $\omega(\alpha)$
is the Denjoy-Wolff point
of the map $g_\alpha:\mathbb C^+\longrightarrow\mathbb C^+\cup\mathbb R$ 
defined by$$g_\alpha(z)=z+\alpha-H(z),\quad z\in\mathbb C^+.$$
\item[{\rm (4)}] The function $\omega$
satisfies $\Re\omega'(z)>1/2,$ $z\in\mathbb C^+$, and
 $$|\omega(z_1)-\omega(z_2)|\ge\frac12|z_1-z_2|,\quad z_1,z_2\in\mathbb
C^+\cup\mathbb R.$$
In particular, $\omega$ is one-to-one.
\item[{\rm (5)}] If $\alpha\in\mathbb R$ is such that $\Im\omega(\alpha)>0$,
then $\omega$ can be
continued analytically to a neighborhood of $\alpha$.
\end{enumerate}
\end{theorem}
\begin{proof}
Let us first consider the function $\varphi(z)=z-H(z)$ which maps $\mathbb C^+$
to its closure, and therefore
can be written in Nevanlinna integral form
$$\varphi(z)=c+bz+\int_{-\infty}^\infty\frac{1+zt}{t-z}\,d\sigma(t),\quad
z\in\mathbb C^+,$$
with $c\in\mathbb R$, $b\ge0$, and $\sigma$ a Borel, finite, positive measure on
the real line. Since clearly
$$\lim_{y\to+\infty}\frac{\varphi(iy)}{iy}=b,$$
and $\Im\varphi(z)\ge b\Im z$, we easily deduce that $a=1-b$ and $a\in(0,1]$.

In order to prove (1), fix $\alpha\in\mathbb C^+$.
We shall argue that the Denjoy-Wolff
point of $g_\alpha$ belongs to $\mathbb C^+$. First,
$$\Im g_\alpha(z)=\Im z+\Im\alpha-\Im H(z)\ge\Im\alpha>0,$$
so that the Denjoy-Wolff point cannot be a real number. On the other hand it
cannot be infinity either
because $$\lim_{y\to\infty}\frac{iy}{iy+\alpha-H(iy)}=\frac1{1-a}>1,$$
where the last fraction must be understood as $+\infty$ if $a=1$. Denote then by
$z_\alpha\in\mathbb C^+$
the Denjoy-Wolff point of $g_\alpha$, and note that it is then the unique fixed
point of $g_\alpha$.
Fixed points $z$ for $g_\alpha$ are precisely points with $H(z)=\alpha$, and the
function
$\omega(\alpha)=z_\alpha$ is analytic, by the implicit function theorem.
This proves (1), and the part of (2) and (3)
pertaining to points $\alpha\in\mathbb C^+$.

Let us note at this point that $g_\alpha$ is not a conformal automorphism of
$\mathbb C^+$, and
therefore $|g'_\alpha(z_\alpha)|<1$. Equivalently,
$$\left|1-\frac1{\omega'(\alpha)}\right|=|1-H'(\omega(\alpha))|=|1-H'(z_\alpha)|<1,$$
so that $\Re\omega'(\alpha)>1/2$ for all $\alpha\in\mathbb C^+$. This implies
the inequality
$$|\omega(z_1)-\omega(z_2)|\ge\frac12|z_1-z_2|,\quad z_1,z_2\in\mathbb C^+.$$
We shall show now that $\omega$ extends by continuity to a 
function $:\overline{\mathbb C^+}\cup\{\infty\}\longrightarrow
\overline{\mathbb C^+}\cup\{\infty\}$, and thus the preceding inequality will 
persist for
$z_1,z_2\in\overline{\mathbb C^+}$ (the proof appears also in 
\cite{BB}, Proposition 2.4).
First assume that there exists a sequence $z_n\to x$ such that the limit
$\lambda=\lim_{n\to\infty}\omega(z_n)$ exists and belongs to ${\mathbb C}^+.$
In this case we have $H(\lambda)=x.$
Denote by $n\geq1$ the order of the zero of
$H(z)-x$ at $z=\lambda$.
We can find analytic functions $\omega_1,\omega_2,\dots,\omega_n$ defined in a
set of the form $\Omega=\{w:0<|w-x|<\delta, w\not\in x-i{\mathbb R}_+\}$ such
that
$H(\omega_j(w))=w$ for $w\in \Omega$ and $j=1,2,\dots,n$. Clearly $\omega$ must
coincide with one of the functions $\omega_j$ on $\Omega\cap{\mathbb C}^+$ and
it follows that $\omega$ extends continuously to the interval $(x-\delta,
x+\delta)$.

Assume to the contrary that there is no sequence $z_n$ as in the first part of
the argument. In other words, if $z_n\to x$ and $\lim_{n\to\infty} \omega(z_n)$
exists, this limit is either infinite or real. Assume now that two sequences
$z_n,w_n\in{\mathbb C}^+$ have limit equal to $x$ and the limits
$\lim_{n\to\infty} \omega(z_n)$, $\lim_{n\to\infty} \omega(w_n)$ exist and
are different. Consider a continuous path $\gamma:(0,1)\longrightarrow
{\mathbb C}^+$ passing through all the points $z_n$ and $w_n$, and such that
$\lim_{t\to1}\gamma(t)=x.$ There exists then an open interval $(\alpha,\beta)
\subset\mathbb R$
such that for every $s\in(\alpha,\beta)$ there is a sequence $t_n\to1$ such
that $\omega(\gamma(t_n))\to s.$ In fact $t_n$ can be chosen so that
$\omega(\gamma(t_n))\to s$ nontangentially as $n\to\infty.$ Since
$H(\omega(\gamma(t_n)))=\gamma(t_n),$ we deduce that the nontangential limit
$H(s)$ of $H$ at $s$ is equal to $x$ almost everywhere. Theorem \ref{Privalov} 
shows now that $H$ must be constant, and this is a contradiction.
Therefore $\lim_{z\to x} \omega(z)$ exists. The case $x=\infty$ is treated
similarily.

To conclude the proof of (2) we still need
to show that $\omega(x)$
is finite if $x\in\mathbb R$. This will follow once we establish (3) for points
on the boundary.

We proceed by proving that (3) and (5) hold for a point $\alpha\in\mathbb R$
with the property that
$\omega(\alpha)\in\mathbb C^+$. Fix such a point $\alpha$, and set
$w=\omega(\alpha)$. With the trivial
exception of $H(z)=z+c, c\in\mathbb R$, the function $g_\alpha$ still maps
$\mathbb C^+$
to $\mathbb C^+$, and $g_\alpha(w)=w$. Indeed,
$$H(w)=H(\omega(\alpha))=\lim_{y\downarrow0}H(\omega(\alpha+iy))=
\lim_{y\downarrow0}(\alpha+iy)=\alpha.$$ Thus $w$ is indeed the Denjoy-Wolff
point of $g_\alpha$.
We conclude that $|g'_\alpha(w)|\le1$, with equality only when $g_\alpha$
is a conformal automorphism of $\mathbb C^+$.
The last case is trivial. Indeed, when $a<1$, we
have $g_\alpha(\infty)=\infty$,
so that $g_\alpha$
can only be an automorphism of the form $g_\alpha(z)=(1-a)z+b$ for some
$b\in\mathbb R$. In this case
$H(z)=az+\alpha-b$, and the conclusions of the theorem are obvious since
$$\omega(z)=\frac1a(z-\alpha+b).$$
If $a=1$ and $g_\alpha$ is an automorphism, we must have
$$g_\alpha(z)=b-\frac c{d-z}$$ for some real numbers
$b,c,d$ with $c>0$. In this case $$H(z)=z+\alpha-b+\frac c{z-d}, \quad
H'(z)=1-\frac c{(z-d)^2},$$
and the zeros of $H'$ are $d\pm\sqrt c\in\mathbb R$. We conclude that
$H'(w)\ne0$. The same conclusion
holds when $g_\alpha$ is not an automorphism because
$$|1-H'(w)|=|g'_\alpha(w)|<1.$$ In both cases, the local inverse of $H$ at $w$ will be
an analytic continuation
of $\omega$ in a neighborhood of $\alpha$.

To conclude, we must complete the proof of (3) by showing that
$\omega(\alpha)$ cannot be infinite
if $\alpha\in\mathbb R$, and that $\omega(\alpha)$ is the Denjoy-Wolff point of
$g_\alpha$.
Fix then $\alpha\in\mathbb R$, and exclude the case $\omega(\alpha)\in\mathbb
C^+$ which has
just been treated. We can also assume that the map $g_\alpha$ is not a conformal
automorphism of $\mathbb C^+$,
since in this case the function $\omega$ can be calculated explicitly.
As in the first part of the proof,
$$\lim_{y\to\infty}\frac{iy}{g_\alpha(iy)}=\lim_{y\to\infty}\frac{iy}{iy+\alpha-H(iy)}=\frac1{1-a}>1,$$
which shows that $\infty$ is not the Denjoy-Wolff point of $g_\alpha$.
 We conclude then that
this Denjoy-Wolff point $z_\alpha$ is a real number. We will show that
$\omega(\alpha)=z_\alpha$.
To do this we observe that
(since the Julia-Carath\'eodory
derivative of $g_\alpha$ at $z_\alpha$ is no more than one) $g_\alpha$ maps the
horodisk
$D_y =\{z:|z-(z_\alpha+iy/2)|<y/2\}$ into itself for $y>0$. In particular, since
$g_\alpha$ has no fixed points in
$\mathbb C^+$,
$$\Im g_\alpha(z_\alpha+iy)<y,\quad y>0.$$Therefore
$$\Im H(z_\alpha+iy)=\Im(z_\alpha+iy+\alpha-g_\alpha(z_\alpha+iy))=y-\Im
g_\alpha(z_\alpha+iy)>0.$$
We also have
$$\lim_{y\downarrow0}H(z_\alpha+iy)=\lim_{y\downarrow0}(z_\alpha+iy+\alpha-g_\alpha(z_\alpha+iy))=\alpha,$$
so that, at long last,
$$\omega(\alpha)=\lim_{y\downarrow0}\omega(H(z_\alpha+iy))=\lim_{y\downarrow0}(z_\alpha+iy)=z_\alpha,$$
as claimed.
\end{proof}

Let now $\mu$ be a Borel probability measure on $\mathbb R.$
We shall use Theorem \ref{InversionC+} to prove the existence of a family
$\{\mu_t\colon t\geq1\}$ such that $\mu_1=\mu$ and $\mu_{t+s}=\mu_t\boxplus
\mu_s.$ We recall that, by Theorem \ref{R-transform}, it will be enough
to show that for each $t>1$ there exists a probability measure $\mu_t$
on the real line such that $\varphi_{\mu_t}(z)=t\varphi_\mu(z)$ for
$z$ in the common domain of the two functions. A slightly weaker version of 
this result can be found in \cite{BB}.

\begin{theorem}\label{mut+}
Consider a Borel probability measure $\mu$ on $\mathbb R$, and a real number
$t\geq1.$
\begin{enumerate}
\item[{\rm (1)}] There exists a probability measure $\mu_t$ satisfying
$\varphi_{\mu_t}(z)=t\varphi_\mu(z)$ for $z$ in the common domain of the two 
functions.
\item[{\rm (2)}] There exists an injective analytic map
$\omega_t:{\mathbb C}^+\rightarrow{\mathbb C}^+$ such that $F_{\mu_t}(z)=
F_\mu(\omega_t(z)),$ for $z\in{\mathbb C}^+.$
\item[{\rm (3)}] We have $\omega_t(z)=\frac{1}{t}z+\left(1-\frac{1}{t}\right)
F_{\mu_t}(z),$ and $H_t(\omega_t(z))=z,$ where $H_t(z)=tz+(1-t)F_\mu(z),$ for 
$z\in{\mathbb C}^+.$
\end{enumerate}
\end{theorem}
\begin{proof} If $t=1$, clearly $\mu_1=\mu$ and $\omega_1(z)=z$ will
satisfy the conclusions of the theorem. Assume therefore that $t>1.$
We clearly have
\begin{eqnarray*}
\Im H_t(z) & = & t\Im z-(t-1)\Im F_\mu(z) \\
       & \leq & t\Im z-(t-1)\Im z \\
          & = & \Im z,
\end{eqnarray*}
and
$$\lim_{y\to+\infty}\frac{H_t(iy)}{iy}
=t-(t-1)\lim_{y\to+\infty}\frac{F_\mu(iy)}{iy}=1.$$
Therefore, Theorem \ref{InversionC+} implies the existence of an analytic 
function $\omega_t:{\mathbb C}^+\rightarrow{\mathbb C}^+$ satisfying
$H_t(\omega_t(z))=z,\ z\in{\mathbb C}^+.$ We also have
$\Im\omega_t(z)\geq\Im z$ and $\lim_{y\to+\infty}\omega_t(iy)/iy=1.$ It follows
that the function
$$F_t(z)=\frac{t\omega_t(z)-z}{t-1},\quad z\in {\mathbb C}^+$$
satisfies the conditions $\Im F_t(z)\geq\Im z,$ $z\in{\mathbb C}^+$, and
$\lim_{y\to+\infty}F_t(iy)/(iy)=1.$ By Proposition \ref{Cauchyalt}, these 
conditions imply the existence of a
Borel probability measure $\mu_t$ on $\mathbb R$ satisfying  $F_{\mu_t}=F_t.$
Note that the definition of $F_t$ yields the first formula in (3). To prove
(2) we observe that
$$F_\mu(z)=\frac{tz-H_t(z)}{t-1},\quad z\in{\mathbb C}^+,$$
so that
$$F_\mu(\omega_t(z))=\frac{t\omega_t(z)-z}{t-1}=F_t(z),\quad z\in{\mathbb C}^+.$$
Finally, let us observe that, for $z$ in the domain of definition of $
\varphi_{\mu_t},$ we have
\begin{eqnarray*}
z & = & F_{\mu_t}\left(z+\varphi_{\mu_t}(z)\right) \\
  & = & F_\mu\left(\omega_t\left(z+\varphi_{\mu_t}(z)\right)\right) \\
  & = & F_\mu\left(\frac{1}{t}\left[z+\varphi_{\mu_t}(z)\right]+
\left(1-\frac{1}{t}\right)F_{\mu_t}(z+\varphi_{\mu_t}(z))\right) \\
  & = & F_\mu\left(\frac{1}{t}\left[z+\varphi_{\mu_t}(z)\right]
+\left(1-\frac{1}{t}\right)\frac{1}{z}\right) \\
  & = & F_\mu\left(z+\frac{1}{t}\varphi_{\mu_t}(z)\right),
\end{eqnarray*}
where we used (2) in the second equality, and (3) in the third equality.
We conclude that the function $\rho(z)=\frac{1}{t}\varphi_{\mu_t}(z)$ satisfies
$\lim_{y\to0}y\rho(iy)=0$ and $F_\mu\left(z+\rho(z)\right)=z.$
Therefore, $\rho(z)=\varphi_\mu(z),$ which proves (1). 
\end{proof}
Due to their properties, it is natural to denote the probability measures
$\mu_t$ from the theorem above by $\mu^{\boxplus t}.$
In the following, we describe regularity properties of probability measures
$\mu^{\boxplus t}$ for $t>1.$ We shall use Theorem \ref{InversionC+} and
Theorem \ref{mut+} to describe the boundary behavior of $F_{\mu^{\boxplus t}}$
and thus characterize the absolutely continuous, singular continuous, and
atomic, parts of $\mu^{\boxplus t}.$ We recall that $-\frac1\pi\Im G_\mu(x+iy)$
is the Poisson integral of the measure $\mu, $ and hence 
$\mu$ can be recovered as the weak*-limit of the measures
$$d\nu_y(x)=-\frac1\pi\Im G_\mu(x+iy)\, dx$$
as $y\downarrow0$. In particular, the density of $\mu$ (relative to Lebesgue
measure) is the
$dx$-a.e. limit of $(-1/\pi)\Im G_\mu(x+iy)$ as $y\downarrow0$, and this limit
is, as seen also in Lemma \ref{DeLaValee}, a.e. infinite relative to
the singular part of the measure $\mu$. More specifically, $x$ is an atom of
$\mu$ if and only of
$F_\mu(x)=0$ and the Julia-Carath\'eodory derivative $F'_\mu(x)$ is finite. 
The value of this derivative is then $$F'_\mu(x)=\frac1{\mu(\{x\})}.$$
These observations can now be used to prove a somewhat stronger version of the
regularity
results proved in \cite{BB}.

\begin{theorem}\label{RegPtmut+}
Let $\mu$ be a probability measure on $\mathbb R$, and let $t>1$.
\begin{enumerate}
\item[{\rm(1)}]
A point $x\in\mathbb R$ satisfies $F_{\mu^{\boxplus t}}(x)=0$ if and only if
$x/t$ is an atom of $\mu$ with
mass $\mu(\{x/t\})\ge(t-1)/t$. If $\mu(\{x/t\})>(t-1)/t$, then $x$ is an atom of
$\mu^{\boxplus t}$, and
$$\mu^{\boxplus t}(\{x\})=t\mu\left(\left\{\frac xt\right\}\right)-(t-1).$$
\item[{\rm(2)}]
The nonatomic part of $\mu^{\boxplus t}$ is absolutely continuous, and its
density is continuous except
at the (finitely many) points $x$ such that $F_{\mu^{\boxplus t}}(x)=0$.
\item[{\rm(3)}]
The density of $\mu^{\boxplus t}$ is analytic at all points where it is
different from zero.
\end{enumerate}
\end{theorem}
\begin{proof}
Assume first that $F_{\mu^{\boxplus t}}(x)=0$. This is equivalent to
$\omega_t(x)=x/t$, or, equivalently, $x/t$ is the Denjoy-Wolff point of the
function $g_x(z)=z-x+H_t(z)$ defined in Theorem \ref{InversionC+} (3). 
By Theorem \ref{DenjoyWolff}, and the definition of $g_x$, this is equivalent to
$H_t(x/t)=x$ and $H'_t(x/t)\ge0$. We have then
$$x=H_t\left(\frac xt\right)=t\frac xt-(t-1)F_\mu\left(\frac
xt\right)=x-(t-1)F_\mu\left(\frac xt\right),$$
so that $F_\mu(x/t)=0$. Moreover, $H_t'(z)=t-(t-1)F'_\mu(z)$, so that
$$t-(t-1)F'_\mu\left(\frac xt\right)\ge 0.$$ This yields $F'_\mu(x/t)\le
t/(t-1)$, indicating that $x/t$ is an atom of
$\mu$ with the required mass. Conversely, if $\mu(\{x/t\})\ge(t-1)/t$, then the
calculations above show
that $H_t(x/t)=x$ and $H'_t(x/t)\ge0$. Thus, $g_x(x/t)=x/t$ and $g'(x/t)\leq1.$
We conclude by Theorem \ref{DenjoyWolff} that $x/t$ is the Denjoy-Wolff point of $g_x$.
Thus, $x/t\in\omega_t(\mathbb R)$, and $\omega_t(x)=x/t$. Therefore
$$F_{\mu^{\boxplus t}}(x)=\frac{t\omega_t(x)-x}{t-1}=0.$$ In case
$\mu(\{x/t\})>(t-1)/t$, we have
$H'_t(x/t)=t-(t-1)F'_\mu(x/t)>0$. This implies that $g_x'(x/t)<1$ and therefore 
$$\lim_{\stackrel{z\longrightarrow x/t}{\sphericalangle}}\frac{H_t(z)-x}{z-x/t}$$
exists and belongs to $(0,1).$ We conclude that $H_t$ maps any nontangential path $
\gamma$ in$ \mathbb C^+$ ending at $x/t$ into a nontangential path $H_t(\gamma)$
ending at $x$, and having the part sufficiently close to $x$ included in $\mathbb C^+.$
Thus, $\omega'_t(x)=1/H'_t(x/t)$. We deduce that
$$F'_{\mu^{\boxplus
t}}(x)=\frac{t\omega_t'(x)-1}{t-1}=\frac{t\frac{\mu(\{x/t\})}{t\mu(\{x/t\})-(t-1)}-1}{t-1}
=\frac1{t\mu(\{x/t\})-(t-1)}.$$ This shows that $x$ is an atom of $\mu^{\boxplus
t}$ with the required mass.

As noted above, $F_{\mu^{\boxplus t}}(x)=0$ almost everywhere relative to the
singular part of $\mu^{\boxplus t}$.
We deduce that this singular part has finite support, and therefore it is purely
atomic. The density of
$\mu^{\boxplus t}$ is simply $\Im(1/F_{\mu^{\boxplus t}}(x))$, and is therefore
continuous at all points
where the denominator is not zero. The analyticity of this density follows from
the analyticity of $\omega_t$
at points $x\in\mathbb R$ where $\omega_t(x)\notin\mathbb R$
(see Theorem \ref{InversionC+} (5)).
\end{proof}

In the following, we apply Theorems \ref{mut+} and \ref{InversionC+} to point 
out connections between free and boolean infinite divisibility and describe 
some properties of reciprocals of Cauchy transforms of infinitely divisible 
measures with respect to free additive convolution. 

\begin{definition}\label{infdiv}
Given an arbitrary convolution $\star,$ we say that a probability
measure $\mu$ is infinitely divisible with respect to $\star$ if
for any $n\in\mathbb N$ there exists a probability measure $\mu_n$ such
that $\mu=\underbrace{\mu_n\star\cdots\star\mu_n}_{n\ \rm
times}.$
\end{definition}
We shall denote by $\mathcal I\mathcal D(\star)$ the set of
all probability measures which are infinitely divisible with
respect to the convolution $\star$.

In \cite{BercoviciVoiculescuIUMJ} Bercovici and Voiculescu have completely 
described infinitely divisible probability measures with respect to free 
additive convolution in terms of their $R-$transforms. We state the result
below:
\begin{theorem}\label{infdiv+}
\
\begin{trivlist}
\item[\ {\rm (i)}] A probability measure $\mu$ on $\mathbb R$ is
$\boxplus$-infinitely divisible if and only if $\varphi_\mu$ has an analytic
extension defined on $\mathbb C^+$ with values in $-\mathbb C^+\cup\mathbb R$.
\item[\ {\rm (ii)}] Let $\varphi\colon\mathbb C^+\longrightarrow-\overline{
\mathbb C^+}$
be an analytic function. Then $\varphi$ is a continuation of $\varphi_\mu$ for
some $\boxplus$-infinitely divisible measure $\mu$ if and only if
$$\lim_{\stackrel{z\longrightarrow\infty}{\sphericalangle}}\frac{\varphi(z)}{z}
=0.$$
\end{trivlist}
\end{theorem}

Next, we give 
an analytic characterization of boolean additive convolution (Definition 
\ref{Boolean}, (a)), due to Speicher and Woroudi:
\begin{theorem}\label{booleananalitic}
Consider two probability measures $\mu$ and $\nu$ on $\mathbb R.$ Then 
$$F_{\mu\uplus\nu}(z)=F_\mu(z)+F_\nu(z)-z,\quad z\in\mathbb C^+.$$
\end{theorem}
(It is easy to observe that the operation $\uplus$ is commutative and 
associative. However, unlike for classical and free additive convolutions, 
$\mu\uplus\delta_a$ is not usually equal to the shift of $\mu$ by the amount
of $a$.)

Theorem \ref{booleananalitic} makes the following result of Speicher and Woroudi
less surprising:
\begin{theorem}\label{infdivboolean}
All probability measures on $\mathbb R$ are infinitely divisible with respect to 
boolean additive convolution.
\end{theorem}
We can now prove the following result:
\begin{proposition}\label{?}
Fix a number $t>1$. With the notations from Theorem \ref{mut+}, the 
following hold: 
\begin{enumerate}
\item[{\rm(a)}] For any probability measure 
$\mu\in\mathcal I\mathcal D(\boxplus)$, we have $\Re F_\mu(z)>1/2$ 
for all $z\in\mathbb C^+$, and $F_\mu$ extends continuously to $\mathbb R$. 
Moreover, $\mu^{sc}=0,$ $\mu$ has at most one atom, and the density of $\mu^{ac}
$ with respect to Lebesgue measure is analytic wherever positive; 
\item[{\rm(b)}] The correspondence $H_t\longleftrightarrow\omega_t$ induces 
a bijective correspondence 
   $$\Psi_t\colon\mathcal I\mathcal D(\uplus)\longrightarrow\mathcal I\mathcal D
   (\boxplus),$$
with the property that $\Psi_t(\mu\uplus\nu)=\Psi_t(\mu)\boxplus\Psi_t(\nu)$
for all probability measures $\mu,\nu.$
\end{enumerate}
\end{proposition}
\begin{proof}
Assume $\mu\in\mathcal I\mathcal D(\boxplus)$. By Theorem \ref{infdiv+}, $
\varphi_\mu$ can be continued analytically to $\mathbb C^+.$ 
Let us observe, by part (ii) of 
Theorem \ref{infdiv+}, the function $H(z)=z+\varphi_\mu(z)$ satisfies the 
hypothesis of Theorem \ref{InversionC+}, with $a=1.$ (Recall that 
$a=\lim_{y\to+\infty}H(iy)/(iy).$) Thus, there exists an analytic self-map 
$\omega$ of the upper half-plane such that $H(\omega(z))=z,$ $z\in\mathbb C^+.$
Applying $\omega$ to both sides of this equality gives $\omega(H(w))=w$ 
for all $w\in\omega(\mathbb C^+).$ Part (4) of Theorem \ref{InversionC+}
implies that $\lim_{y\to+\infty}\omega(iy)/(iy)\geq1/2$, so we conclude that
$\lim_{y\to+\infty}\omega(iy)/(iy)=1.$ But $F_\mu(H(z))=z$ for $z$ in some
large enough truncated cone at infinity. We conclude that $\omega=F_\mu$, and
by part (4) of Theorem \ref{InversionC+}, we obtain that $\Re F_\mu(z)>1/2$ for 
all $z\in\mathbb C^+.$ The continuous extension of $F_\mu$ to $\mathbb R$ follows
from part (2) of the same theorem. 

By definition of infinite divisibility, we know that
$\mu=\mu_{1/2}\boxplus\mu_{1/2}$ for some probability measure $\mu_{1/2}$ on
$\mathbb R$. Thus, the statements refering to the absolutely continuous and 
singular continuous parts of $\mu$ follow from parts (2) and (3) of Theorem 
\ref{RegPtmut+}. Finally, by part (1) of Theorem \ref{RegPtmut+} we know that
in order for $x$ to be an atom of $\mu$, we must have $F_\mu(x)=0$. Since
$\Re F_\mu(z)>1/2$ for all $z\in\mathbb C^+,$ we conclude that this can 
happen for at most one $x\in\mathbb R$. This completes the proof of (a).

To prove (b), consider two arbitrary probability measures $\mu,\nu$ on $\mathbb
R$. Define $H_t^\mu(z)=tz+(1-t)F_\mu(z),$ $H_t^\nu(z)=tz+(1-t)F_\nu(z).$ As we 
have seen in the proof of (a),
the right inverses $\omega_t^\mu$ and $\omega_t^\nu$ of $H_t^\mu$ and $H_t^\nu,
$ respectively, are reciprocals of Cauchy 
transforms of probability measures which are infinitely divisible with respect
to $\boxplus$; denote them by $\Psi_t(\mu)$ and $\Psi_t(\nu).$
Observe that, by its definition, 
$$\varphi_{\Psi_t(\mu)}(z)=H_t^\mu(z)-z=(1-t)(F_\mu(z)-z),\quad z\in\mathbb C^+.
$$
Thus, 
\begin{eqnarray*}
\varphi_{\Psi_t(\mu\uplus\nu)}(z) & = & H_t^{\mu\uplus\nu}(z)-z\\
& = & (1-t)(F_{\mu\uplus\nu}(z)-z)\\
& = & (1-t)(F_\mu(z)-z+F_\nu(z)-z)\\
& = & (H_t^\mu(z)-z)+(H_t^\nu(z)-z)\\
& = & \varphi_{\Psi_t(\mu)}(z)+\varphi_{\Psi_t(\nu)}(z),\quad z\in\mathbb C^+.
\end{eqnarray*}
(We have used Theorem \ref{booleananalitic} in the third equality.)
We conclude by Theorem \ref{R-transform} that $\Psi_t(\mu\uplus\nu)=
\Psi_t(\mu)\boxplus\Psi_t(\nu).$ This shows that $\Psi_t$ is an injective 
homomorphism. To show surjectivity, consider $\lambda\in\mathcal I\mathcal D
(\boxplus)$ and let $H(z)=z+\varphi_\lambda(z),$ $z\in\mathbb C^+.$ 
By Theorem \ref{infdiv+}, we have $\Im H(z)\leq\Im z$ for all $z$ in the
upper half-plane, and $\lim_{y\to+\infty}H(iy)/(iy)=1.$ 
Let 
 $$F(z)=\left(1+\frac{1}{t-1}\right)z-\frac{1}{t-1}H(z),\quad z\in\mathbb C^+.
$$ 
Observe that 
$$\Im F(z) = \left(1+\frac{1}{t-1}\right)\Im z - \frac{1}{t-1}\Im H(z)\geq
\left(1+\frac{1}{t-1}\right)\Im z - \frac{1}{t-1}\Im z =\Im z,$$
for all $z\in\mathbb C^+,$ and
$$\lim_{y\to+\infty}\frac{F(iy)}{iy}=1+\frac{1}{t-1}-\frac{1}{t-1}
\lim_{y\to+\infty}\frac{H(iy)}{iy}=1.$$
We conclude by Proposition \ref{Cauchyalt} that there exists a probability
measure $\mu$ on $\mathbb R$ such that $F=F_\mu$. A computation similar to
the one above shows that $\Psi_t(\mu)=\lambda.$
\end{proof}

A somewhat weaker version of part (a) of the above theorem can be found in
\cite{BVSemigroup}, while for $t=2,$ the morphism described in part (b) is the 
Bercovici-Pata bijection (see \cite{BPAnnals} for details).

\section[]{Partially defined semigroups with respect to 
free multiplicative convolution on the unit circle}

A result similar to the one from Theorem \ref{mut+} holds for 
free multiplicative convolutions of measures supported on the
unit circle. However, as we shall see below, it is not always possible to 
start the semigroup at $t=1$, and the semigroup is not unique.
We begin with an inversion theorem for maps defined in the unit
disk.
\begin{theorem}\label{InversionX1}
Let $\Phi:\mathbb D\to\mathbb C\cup\{\infty\}$ be a meromorphic function 
such that $\Phi(0)=0$ and
$|\Phi(z)|\ge |z|$ for all $z\in\mathbb D$. Then:
\begin{enumerate}
\item[{\rm (1)}] For every $\alpha\in\mathbb D$ there exists a unique
$z\in\mathbb D$ such that
$\Phi(z)=\alpha$.
\item[{\rm (2)}] There exists a continuous function $\omega:\overline{\mathbb
D}\longrightarrow\overline{\mathbb D}$
such that $\omega(\mathbb D)\subset\mathbb D$, $\omega|_\mathbb D$ is analytic,
and $\Phi(\omega(z))=z$ for $z\in\mathbb D$.
\item[{\rm (3)}] For each $\alpha\in\overline{\mathbb D}$, $\omega(\alpha)$ 
is the Denjoy-Wolff point of the map $g_\alpha:\mathbb D\to
\overline{\mathbb D}$ defined by
$$g_\alpha(z)=\frac{\alpha z}{\Phi(z)},
\quad  z\in\mathbb D.$$
\item[{\rm (4)}] The function $\omega$ is one-to-one.
\item[{\rm (5)}] If $\zeta\in\mathbb T $ is such that $|\omega(\zeta)|<1$, 
then $\omega$ can be continued analytically to a neighborhood of $\zeta$.
\end{enumerate}
\end{theorem}
\begin{proof}
If $\alpha=0$, then $z=0$ is indeed the unique point for which $\Phi(z)=0$.
Assume therefore
that $\alpha\in\mathbb D\setminus\{0\}$.
The function $g_\alpha$ is analytic on $\mathbb D$, and
$$|g_\alpha(z)|\le|\alpha|<1,\quad z\in\mathbb D.$$ We conclude that the
Denjoy-Wolff point
of $g_\alpha$ is in $\mathbb D$, and therefore there exists a unique
$z_\alpha\in\mathbb D$ such that
$g_\alpha(z_\alpha)=z_\alpha$, namely this Denjoy-Wolff point.
Since
$$g_\alpha(0)=\frac\alpha{\Phi'(0)}\ne0,$$
the point $z_\alpha$ cannot be zero, and therefore the equation
$g_\alpha(z_\alpha)=z_\alpha$ is equivalent to $\Phi(z)=\alpha$. The function
$\omega(\alpha)=z_\alpha$
is easily seen to be analytic.

To complete the proof of (2), we must show that for every $\zeta\in\mathbb T $, the limit
$\lim_{z\to\zeta}\omega(z)$ exists. We distinguish two cases:\begin{enumerate}
\item[(a)] There exists a sequence $\alpha_n\in\mathbb D$ such that
$\lim_{n\to\infty}\alpha_n=\zeta$
and the limit $\lim_{n\to\infty}\omega(\alpha_n)$ exists and belongs to $\mathbb D$;
\item[(b)] For any sequence $\alpha_n\in\mathbb D$ such that
$\lim_{n\to\infty}\alpha_n=\zeta$
and the limit $\lim_{n\to\infty}\omega(\alpha_n)$ exists, this limit is in $\mathbb T
$.
\end{enumerate}
Assume first that (a) holds, and set $w=\lim_{n\to\infty}\omega(\alpha_n)$. We
claim that
$g_\zeta(w)=w$. Indeed,
$$\Phi(w)=\lim_{n\to\infty}\Phi(\omega(\alpha_n))=\lim_{n\to\infty}\alpha_n=\zeta.$$
We deduce that $w$ is the Denjoy-Wolff point of $g_\zeta$, and hence
$|g'_\zeta(w)|\le1$, with strict inequality unless $g_\zeta$ is a conformal
automorphism of $\mathbb D$.
This last case means that$$g_\zeta(z)=\beta\frac{z-\gamma}{1-\overline\gamma
z}$$ for some
$\beta\in\mathbb T $ and $\gamma\in\mathbb D$. In this case
$$\Phi(z)=\frac\zeta\beta z\frac{1-\overline\gamma z}{z-\gamma},\quad
\Phi'(z)=\frac\zeta\beta
\frac{-\overline\gamma z^2+2\gamma\overline\gamma z-\gamma}{(z-\gamma)^2}.$$
If $\gamma=re^{it}$, then the zeros of $\Phi'$ are equal to $e^{it}(r\pm
i\sqrt{1-r^2})$.
Both of these zeros are on the unit
circle, so that $\Phi'(w)\ne0$. We deduce that $\Phi$ is locally invertible at
$w$, and the local inverse is an analytic
continuation of $\omega$ to a neighborhood of $\zeta$.

When $g_\zeta$ is not a conformal automorphism, we must have $|g'_\zeta(w)|<1$.
Since
$$g'_\zeta(w)=\zeta\frac{\Phi(w)-w\Phi'(w)}{\Phi(w)^2}=1-w\frac{\Phi'(w)}\zeta,$$ we deduce that $\Phi'(w)\ne0$,
and the analytic continuation of $\omega$ at $\zeta$ is obtained as in the
previous case.

This argument shows that, in case (a) holds, $\omega$ has an analytic, in
particular continuous,
continuation at $\zeta$, which also proves assertion (5).

Next we deal with case (b). Assume that there exist two sequences
$\alpha_n,\beta_n\in\mathbb D$
such that $\lim_{n\to\infty}\alpha_n=\lim_{n\to\infty}\beta_n=\zeta$, the limits
$\lim_{n\to\infty}\omega(\alpha_n)$, $\lim_{n\to\infty}\omega(\beta_n)$ exist,
and they are different.
Choose a continuous path $\gamma:[0,1)\to\mathbb D$ passing through all the
points $\alpha_n,\beta_n$,
and such that $\lim_{t\to1}\gamma(t)=\zeta$. It follows from our assumption that
there exists an open interval
$I$ in the unit circle $\mathbb T $ such that every point $\xi\in I$ is of the
form $\xi=\lim_{n\to\infty}
\omega(\gamma(t_{\xi,n}))$, where $t_{\xi,n}\in[0,1)$ for all $n$. Moreover,
the sequence $\{t_{\xi,n}\}_{n\in\mathbb N}$ can be chosen
so that $\omega(\gamma(t_{\xi,n}))/\xi>0$, and thus $\omega(\gamma(t_{\xi,n}))$
approaches $\xi$
nontangentially. Since
$$\lim_{n\to\infty}\Phi(\omega(\gamma(t_{\xi,n})))=
\lim_{n\to\infty}\gamma(t_{\xi,n})=\zeta,$$
we conclude that the nontangential limit of $\Phi$ at $\xi$ is $\zeta$ for
almost every $\xi\in I$. 
Privelov's theorem
(Theorem \ref{Privalov}) implies that $\Phi$ is identically
equal to $\zeta$.
This however is not true because $\Phi(0)=0$. This contradiction shows that the
limit $\omega(\zeta)=
\lim_{z\to\alpha,z\in\mathbb D}\omega(z)$ exists in case (b) as well.

We have already shown that $\omega(\alpha)$ is the Denjoy-Wolff point of
$g_\alpha$ provided that
$|\alpha|<1$. We will show that this is still true when $\alpha\in\mathbb T $.
Fix $\alpha\in\mathbb T $, and denote by $z_\alpha$ the Denjoy-Wolff point of
$g_\alpha$. If
$|z_\alpha|<1$, then $g_\alpha(z_\alpha)=z_\alpha$ so that
$\Phi(z_\alpha)=\alpha$. Since $\Phi$ is analytic
at $z_\alpha$, hence open, there exist points $z_r$ such that
$\lim_{r\uparrow1}z_r=z_\alpha$ and
$\Phi(z_r)=r\alpha$ for $r$ close to 1. However, property (1) shows that
$\omega(r\alpha)$ is
the only point in $\mathbb D$ such that $\Phi(\omega(r\alpha))=r\alpha$. We
deduce that
$$\omega(\alpha)=\lim_{r\uparrow1}\omega(r\alpha)=\lim_{r\uparrow1}z_r=z_\alpha,$$
as claimed.

On the other hand, if $|z_\alpha|=1$, Julia's theorem shows that, for each
$\varepsilon\in(0,1/2)$,
$g_\alpha$ maps the horodisk $\{w\in\mathbb
C:|w-(1-\varepsilon)z_\alpha|<\varepsilon\}$ into itself.
In particular, as $g_\alpha$ is not the identity map,
$|g_\alpha((1-2\varepsilon)z_\alpha)|>1-2\varepsilon$ and
$|z_\alpha-g_\alpha((1-2\varepsilon)z_\alpha)|<2\varepsilon$. Therefore
$$|\Phi((1-2\varepsilon)z_\alpha)|=\left|\frac{\alpha(1-2\varepsilon)z_\alpha}{g_\alpha((1-2\varepsilon)z_\alpha)}
\right|<1$$
and $\lim_{\varepsilon\downarrow0}\Phi((1-2\varepsilon)z_\alpha)=\alpha$. We
conclude that
$$\omega(\alpha)=\lim_{\varepsilon\downarrow0}\omega(\Phi((1-2\varepsilon)z_\alpha))
=\lim_{\varepsilon\downarrow0}(1-2\varepsilon)z_\alpha=z_\alpha.$$

It remains to prove that $\omega$ is one-to-one on $\overline{\mathbb D}$. To do
this, we observe first
that for each $\alpha\in\mathbb D$, $\omega(\alpha)$ is the unique fixed point
of $g_\alpha$, and
$g_\alpha$ is not a conformal automorphism of $\mathbb D$ as
$$g_\alpha(\mathbb D)\subset\{z:|z|\le|\alpha|\}\ne\mathbb D.$$
Thus $|g'_\alpha(\omega(\alpha))|<1$,
so that $$\left|1-\omega(\alpha)\frac{\Phi'(\omega(\alpha))}\alpha\right|<1.$$
Using the relation
$\Phi'(\omega(\alpha))=1/\omega'(\alpha)$, the inequality becomes
$$\left|1-\frac{\omega(\alpha)}{\alpha\omega'(\alpha)}\right|<1$$
which is equivalent to
$$\Re\left(\frac{\alpha\omega'(\alpha)}{\omega(\alpha)}\right)>\frac12.$$
Consider now the logarithm of $\omega$. More precisely, let
$h:(0,1]\times\mathbb R\to\mathbb C$
be a continuous function such that $e^{h(r,t)}=\omega(re^{it})$ for all $r$ and
$t$. We have
$$\frac d{dt}h(r,t)=i\frac{re^{it}\omega'(re^{it})}{\omega(re^{it})},$$ so that
$$\Im\frac d{dt}h(r,t)>\frac12$$
provided that $r<1$.
This implies the inequality
$$|h(r,t_1)-h(r,t_2)|\ge\frac12|t_1-t_2|$$
whenever $r<1$. This inequality extends by continuity to $r=1$ as well. Assume
now that
$0\le t_1<t_2<2\pi$ are such that $\omega(e^{it_1})=\omega(e^{it_2})$. In this
case there must exist
an integer $k$ such that $h(1,t_1)=h(1,t_2+2k\pi)$, and this leads to a
contradiction since
$t_1\ne t_2+2k\pi$.
\end{proof}

As in Section 1, we use the above theorem to prove the existence of
a partially defined free convolution semigroup of measures on $\mathbb T.$

\begin{theorem}\label{mutX1}
Consider a probability measure $\mu$ on $\mathbb T$ having nonzero first 
moment, and let $t>1$ be a real number.
Assume that the function $\eta_\mu$ never vanishes on $\mathbb D\setminus\{0\}.
$
\begin{enumerate}
\item[{\rm (1)}] There exists a probability measure $\mu_t$ on $\mathbb T$ 
such that 
$\int_{\mathbb T}\zeta\,d\mu_t(\zeta)\ne0$, and $\Sigma_{\mu_t}(z)=\Sigma_\mu(z)^t$ in a
neighborhood
of zero. Moreover, $\eta_{\mu_t}$ never vanishes on $\mathbb D\setminus\{0\}$.
\item[{\rm (2)}] There exists an analytic function 
$\omega_t:\mathbb D\longrightarrow\mathbb D$ such that
$|\omega_t(z)|\le|z|$ and $\eta_{\mu_t}(z)=\eta_\mu(\omega_t(z))$ for
$z\in\mathbb D$.
\item[{\rm (3)}] We have
$$\omega_t(z)=\eta_{\mu_t}(z)\left[\frac z{\eta_{\mu_t}(z)}\right]^{1/t},\quad
z\in\mathbb D,$$ and $\Phi_t(\omega_t(z))=z$ for $z\in\mathbb D$, where

$$\Phi_t(z)=z\left[\frac z{\eta_\mu(z)}\right]^{t-1},\quad z\in\mathbb D.$$
\end{enumerate}
\end{theorem}
\begin{proof}
The function $\Phi_t$, defined as in (3), satisfies the conditions of Theorem
\ref{InversionX1} (this is an immediate consequence of Proposition 
\ref{etacircle} and the hypothesis),
which proves the existence of a right inverse $\omega_t$ for $\Phi_t.$ We
define then $\eta_t=\eta_\mu\circ\omega_t$, and observe that 
$\eta_t=\eta_{\mu_t}$ for some measure $\mu_t$ on $\mathbb T$. Indeed, 
this follows immediatly from Theorem \ref{InversionX1} (2) and 
Proposition \ref{etacircle}. 
 Clearly $\eta_{\mu_t}$ only vanishes for $z=0$; this follows
from the corresponding hypothesis on $\mu$, and from the fact that $\omega_t$ 
is injective.

To conclude, we must verify the identities in (1) and (3).
The inversion relation $\Phi_t(\omega_t(z))=z$
amounts to
$$\omega_t(z)\left[\frac{\omega_t(z)}{\eta_\mu(\omega_t(z))}\right]^{t-1}=z,$$
or
$$\left[\frac{\omega_t(z)}{\eta_{\mu_t}(z)}\right]^{t-1}=\frac z{\omega_t(z)}.$$
This can also be written as
$$\left[\frac{\omega_t(z)}z\right]^{t-1}\left[\frac
z{\eta_{\mu_t}(z)}\right]^{t-1}=\frac z{\omega_t(z)},$$
by injectivity of $\omega_t$ and the hypothesis on $\eta_\mu$. 
$$\left[\frac{\omega_t(z)}z\right]^{t}=\left[\frac
z{\eta_{\mu_t}(z)}\right]^{1-t}
=\frac z{\eta_{\mu_t}(z)}\left[\frac{\eta_{\mu_t}(z)}z\right]^t$$
and hence
$$\frac{\omega_t(z)}z=\left[\frac
z{\eta_{\mu_t}(z)}\right]^{1/t}\frac{\eta_{\mu_t}(z)}z,$$
which verifies (3).
Finally, if $z$ is close to zero, we have
\begin{eqnarray*}
z\Sigma_{\mu_t}(z)&=\eta_{\mu_t}^{-1}(z)=\omega_t^{-1}(\eta_\mu^{-1}(z))
=\Phi_t(\eta^{-1}_\mu(z))
=\eta^{-1}_\mu(z)\left[\frac{\eta^{-1}_\mu(z)}{\eta_\mu(\eta^{-1}_\mu(z))}\right]^{t-1}\\
&=\eta^{-1}_\mu(z)\left[\frac{\eta^{-1}_\mu(z)}z\right]^{t-1}
=z\Sigma_\mu(z)\Sigma_\mu(z)^{t-1}.
\end{eqnarray*}
Dividing by $z$ we obtain the desired relation
$\Sigma_{\mu_t}(z)=\Sigma_\mu(z)^t$.
\end{proof}

It seems natural to denote $\mu_t$ by $\mu^{\boxtimes t}$, but note that these
measures are
only determined up to a rotation by a multiple of $2\pi t$. This is due to the
choice involved in
extracting roots. In case $\eta_\mu$ has zeroes in $\mathbb D\setminus\{0\}$,
we are assured by Theorem \ref{Biane3} that $\eta_{\mu\boxtimes\mu}$ has no
zeros in $\mathbb D\setminus\{0\}$.
Thus, one can still define
$$\mu^{\boxtimes t}=(\mu\boxtimes\mu)^{\boxtimes(t/2)},$$
at least for $t>2$. It will still be true that $\eta_{\mu^{\boxtimes
t}}=\eta_\mu\circ\omega_t$
for some conformal map $\omega_t$, but $\omega_t$ might not have a globally
defined left
inverse.  Its left inverse $z(z/\eta_\mu(z))^{t-1}$ is however the composition
of two functions
defined on $\mathbb D$ because $\omega_t$ is the composition of two functions
which do
have global left inverses.

Let us also observe that whenever $\int_{\mathbb T}\zeta\,d\mu(\zeta)=0$,
$\mu_t$ is defined for any $t\geq 2$, and equals the uniform distribution
(the Haar measure) on $\mathbb T$.

In the following, we shall discuss regularity properties of measures $\mu_t$, $
t>1.$
As we have observed at the beginning of Section 1.3 of the previous chapter,
in the case of probability measures $\mu$ on $\mathbb T$, we have
$$\frac1\pi\left(\psi_\mu(z)+\frac12\right)=\frac1{2\pi}\int_0^{2\pi}\frac{e^{it}+z}{e^{it}-z}\,d\mu(e^{-it}),
\quad z\in\mathbb D.$$
The real part of this function is the Poisson integral of the measure
$d\mu(e^{-it})$, and therefore
this measure can be obtained as the weak*-limit of the measures
$$d\nu_r(e^{it})=\Re\frac1\pi\left(\psi_\mu(re^{it})+\frac12\right)\,dt$$
as $r\uparrow 1$. In particular, the functions
$(1/\pi)\Re(\psi_\mu(re^{it})+1/2)$ converge $dt$-a.e.
to the density of $\mu(e^{-it})$, and they converge to infinity a.e. relative to
the singular part of this
measure.
A number $1/\zeta\in\mathbb T$ is an
atom of $\mu$ if and only if
$\eta_\mu(\zeta)=1$ and the Julia-Carath\'eodory derivative $\eta'_\mu(\zeta)$
is finite (this follows from Lemma \ref{atoms&tang} and Theorem 
\ref{JuliaCaratheodory}). In this case this derivative is still given by
$$\zeta\eta'_\mu(\zeta)=\frac1{\mu(\{1/\zeta\})}.$$
In the following result we will assume that a choice of arguments has been made
so that the results of Theorem \ref{mutX1} are correct. When $t$ is an integer, the
following result is true without the assumption that $\eta_\mu$ have no zeros 
different from $0$.

\begin{theorem}\label{regmutX1}
Let $\mu$ be a probability measure on $\mathbb T$ such that $\eta_\mu'(0)\ne0$,
$\eta_\mu$ has no zeros
in $\mathbb D\setminus\{0\}$, and let $t>1$. Consider the measure
$\mu^{\boxtimes t}$ provided by Theorem \ref{mutX1},
and the corresponding functions $\Phi_t$ and $\omega_t$
\begin{enumerate}
\item[{\rm(1)}]
A point $\zeta\in\mathbb T$ satisfies
$\eta_{\mu^{\boxtimes t}}(\zeta)=1$ if and only if
$\omega_t(\zeta)\in\mathbb T$ and
$1/\omega_t(\zeta)$ is an atom of $\mu$ with
mass $\mu(\{1/\omega_t(\zeta)\})\ge(t-1)/t$. If
$\mu(\{1/\omega_t(\zeta)\})>(t-1)/t$, then $1/\zeta$ is
an atom of $\mu^{\boxtimes t}$, and
$$\mu^{\boxtimes t}(\{1/\zeta\})=t\mu(\{1/\omega_t(\zeta)\})-(t-1).$$
If $\eta_{\mu^{\boxtimes t}}(\zeta)=1$, there is a real number $\theta$ such
that $\zeta=e^{i\theta}$
and $\omega_t(\zeta)=e^{i\theta/t}$.
\item[{\rm(2)}]
The nonatomic part of $\mu^{\boxtimes t}$ is absolutely continuous, and its
density is continuous except
at the (finitely many) points $\zeta$ such that $\eta_{\mu^{\boxtimes
t}}(\zeta)=1$.
\item[{\rm(3)}]
The density of $\mu^{\boxtimes t}$ is analytic at all points where it is
different from zero.
\end{enumerate}
\end{theorem}
\begin{proof}
As in the proof of Theorem \ref{RegPtmut+}, (2) and (3) follow easily from
the continuity and analyticity properties of $\omega_t$, once (1)
is established. The formulas
$$\eta_{\mu^{\boxtimes t}}(z)=\eta_\mu(\omega_t(z))=
z\left[\frac{\omega_t(z)}z\right]^{\frac t{t-1}},
\quad z\in \mathbb D,$$ show that $\eta_{\mu^{\boxtimes t}}$
extends continuously to $\overline{\mathbb D}$, and these
equalities persist when $z\in\mathbb T$ even if
$|\omega_t(z)|=1$, provided that $\eta_\mu(\omega_t(z))$ is viewed
as a nontangential limit. This being said, let $\zeta\in\mathbb T$
satisfy $\eta_{\mu^{\boxtimes t}}(\zeta)=1$. The above formulas
show that $|\omega_t(\zeta)|=1$, and in fact $\omega_t(\zeta)$ is
one of the $(1/t)$-powers of $\zeta$. Since, according
to Theorem \ref{InversionX1}, $\omega_t(\zeta)$ is the 
Denjoy-Wolff point of the function $g_\zeta(z)=\zeta z/\Phi_t(z)$,
we see, by Theorem \ref{DenjoyWolff}, that $\Phi_t(\omega_t(\zeta))=\zeta$, and
$$\frac{\omega_t(\zeta)}{\zeta}\Phi'_t(\omega_t(\zeta))\ge0.$$
We
will treat this derivative a little more carefully, on account of
the fractional powers involved. We have
$$\frac{\Phi'_t(z)}{\Phi_t(z)}=\frac
tz-(t-1)\frac{\eta'_\mu(z)}{\eta_\mu(z)}.$$Replacing $z$ by
$\omega_t(\zeta)$, and using the equalities
$\Phi_t(\omega_t(\zeta))=\zeta$ and
$\eta_\mu(\omega_t(\zeta))=\eta_{\mu^{\boxtimes t}}(\zeta)=1$, we
obtain
$$\frac{\omega_t(\zeta)}\zeta\Phi'_t(\omega_t(\zeta))=t-(t-1)\omega_t(\zeta)
\eta'_\mu(\omega_t(\zeta))=\frac{t\mu(\{1/\omega_t(\zeta)\})-(t-1)}{
\mu(\{1/\omega_t(\zeta)\})}.$$This number is nonnegative precisely
when $\mu(\{1/\omega_t(\zeta)\})\ge(t-1)/t$. 
Conversely, if $\omega_t(\zeta)\in\mathbb T$, and
$\mu(\{1/\omega_t(\zeta)\})\ge(t-1)/t$, then the above calculation
will show that $\Phi_t$ has nonnegative Julia-Carath\'eodory
derivative at $\omega_t(\zeta)$, and, since $g_\zeta(\omega_t(\zeta))=
\omega_t(\zeta)$, we conclude that $\omega_t(\zeta)$ is the Denjoy-Wolff point
of $g_\zeta,$ so that
$\Phi_t(\omega_t(\zeta))=\zeta$. This 
implies
$\eta_{\mu^{\boxtimes t}}(\zeta)=1$. The mass $\mu^{\boxtimes
t}(\{1/\zeta\})$ can then be calculated using the chain rule in
case $\mu(\{1/\omega_t(\zeta)\})>(t-1)/t$, as in the proof of
Theorem \ref{RegPtmut+}.
\end{proof}

A connection similar to the one presented in Proposition \ref{?} exists for
measures on the unit circle. We present first an analytic characterization
of multiplicative boolean convolution (Definition \ref{Boolean} (b)), due to 
Franz \cite{Franz}:

\begin{theorem}\label{utimesanalitic}
Consider  two probability measures $\mu$ and $\nu$ on $\mathbb T$. We have
$$z\eta_{\mu\cup\kern-.5383em\lower-.65ex\hbox{$_\times$}\nu}(z)=\eta_\mu(z)\eta_\nu(z),\quad z\in\mathbb D.$$
\end{theorem}
Unlike in the additive case, not all probability measures are infinitely
divisible with respect to multiplicative boolean convolution, as the following
theorem of Franz shows (see \cite{Franz}, Theorem 3.6):

\begin{theorem}\label{Franzinfdiv}
A probability measure $\mu$ on the unit circle, different from the Haar 
measure, is infinitely divisible with
respect to multiplicative boolean convolution if and only if there 
exists a real number $b\in[0,2\pi)$ and a finite measure $\rho$ on $\mathbb T$ 
such that 
$$\eta_\mu(z)=z\exp\left(ib-\int_\mathbb T\frac{e^{it}+z}{e^{it}-z}d\rho(e^{it}
)\right),\quad z\in\mathbb D$$
\end{theorem}
Observe that, if $\mu$ is the uniform distribution on $\mathbb T,$ then 
$\eta_\mu(z)=0$ for all $z\in\mathbb D$, and thus $\mu$ is infinitely divisible
with respect to $\Utimes.$ 

We state below the analogue of Theorem \ref{infdiv+} for free multiplicative 
convolution, as it appears in \cite{BercoviciVoiculescuPJM}, Theorem 6.7:
\begin{theorem}\label{infdivX1}
\
\begin{enumerate}
\item[{\rm(i)}] A probability measure $\mu$ on $\mathbb T$ with nonzero first 
moment is 
infinitely divisible if and only if there exists a function $u(z)$, analytic
in $\mathbb D$, such that $\Re u(z)\geq0$ for all $z\in\mathbb D$ and
$\Sigma_\mu(z)=\exp(u(z)).$
\item[{\rm(ii)}] Let $u\colon\mathbb D\longrightarrow-i\overline{\mathbb C^+}$
be an analytic function. Then the function $\Sigma(z)=\exp(u(z))$ has the
form $\Sigma=\Sigma_\mu$ for some $\boxtimes$-infinitely divisible measure
$\mu$ with nonzero first moment.
\end{enumerate}
\end{theorem}
As in the case of multiplicative boolean convolution, we observe that 
the uniform distribution is also infinitely divisible with respect to $\boxtimes.$

\begin{proposition}\label{??}
Fix a number $t>1$. With the notations from Theorem \ref{mutX1}, the following
hold:
\begin{enumerate}
\item[{\rm(a)}] For any $\boxtimes$-infinitely divisible probability measure 
$\mu$ on $\mathbb T$, different from the uniform distribution on the unit 
circle, $\eta_\mu$ extends continuously to $\mathbb T,$ and the extension is
injective. Moreover, $\mu^{sc}=0,$ $\mu$ has at most one atom, and the density
of $\mu^{ac}$ with respect to the uniform distribution on $\mathbb T$ is
analytic wherever positive;
\item[{\rm(b)}] The correspondence $\Phi_t\longleftrightarrow\omega_t$ induces
a bijective correspondence 
 $$\Psi_t\colon\mathcal I\mathcal D(\Utimes)\longrightarrow\mathcal I\mathcal D
(\boxtimes),$$
with the property that $\Psi_t(\mu\utimes\nu)=\Psi_t(\mu)\boxtimes\Psi_t(\nu)$
for all $\mu,\nu\in\mathcal I\mathcal D(\Utimes).$
\end{enumerate}
\end{proposition}
\begin{proof}
Assume $\mu\in\mathcal I\mathcal D(\boxtimes).$ 
By definition, $\mu=\mu_{1/2}\boxtimes\mu_{1/2}$ for some probability measure $
\mu_{1/2}$ on $\mathbb T$ with nonzero first moment. Thus, by Theorem 
\ref{regmutX1}, $\mu^{sc}=0,$ and $\mu^{ac}$ has analytic density with respect
to the uniform distribution on $\mathbb T$.
To prove injectivity of $\eta_\mu,$ observe that, by Theorem \ref{infdivX1}, 
$\Sigma_\mu=\exp\circ u$ for some analytic function $u$ on $\mathbb D$
with $\Re u(z)\geq0,$ $|z|<1.$ We claim that the function $\Phi(z)=z\Sigma_\mu
(z)$ satisfies the hypotheses of Theorem \ref{InversionX1}. Indeed, 
$\Phi(0)=0$ trivially, and 
$$|\Phi(z)|=|z||\exp(u(z))|\geq|z|e^{\Re u(z)}\geq|z|,\quad z\in\mathbb D,$$
by Theorem \ref{infdivX1} (i). Thus, there exists an analytic self-map of the
unit disk $\omega$ such that $\Phi(\omega(z))=z,$ $z\in\mathbb D$. 
Applying $\omega$ to both sides of the equality gives $\omega(\Phi(w))=w$ for
all $w\in\omega(\mathbb D).$
Since $\Phi(z)=0$ if and only if $z=0,$ we must also have $\omega(0)=0.$
But $\eta_\mu(\Phi(z))=z$ for $z$ in some neighbourhood of $0$. We conclude 
that $\eta_\mu=\omega,$ and thus, by part (4) of Theorem \ref{InversionX1},
we conclude that $\eta_\mu$ exitends as a continuous injective function to 
$\overline{\mathbb D}.$ Finally, by part (1) of Theorem \ref{regmutX1}, 
we know that in order for $x\in\mathbb T$ to be an atom of $\mu$, we must have
$\eta_\mu(x)=1.$ The injectivity of $\eta_\mu$ implies that this can happen for
at most one $x\in\mathbb T.$ This proves part (a) of the proposition.

To prove part (b), consider two arbitrary probability measures $\mu,\nu
\in\mathcal I \mathcal D(\Utimes).$ Assume first that neither of them is the 
uniform distribution on $\mathbb T$. Define 
$$\Phi_t^\mu(z)=z\left[\frac{z}{\eta_\mu(z)}\right]^{t-1},\quad
\Phi_t^\mu(z)=z\left[\frac{z}{\eta_\mu(z)}\right]^{t-1},\quad z\in\mathbb D.$$
Observe that, by Theorem \ref{Franzinfdiv}, the above functions are indeed
well-defined, and 
$|\Phi_t^\mu(z)|\ge|z|$, $|\Phi_t^\nu(z)|\geq|z|.$
As we have seen in the proof of (a), there exist probability measures, denote
them by 
$\Psi_t(\mu)$ and $\Psi_t(\nu)$
such that the $\omega_t^\mu=\eta_{\Psi_t(\mu)},$ and $\omega_t^\nu=\eta_{
\Psi_t(\nu)}$,
where $\omega_t^\mu,\omega_t^\nu$ denote the right inverses of 
$\Phi_t^\mu$ and $\Phi_t^\nu$, respectively, provided by Theorem \ref{mutX1}.
We claim that $\Psi_t(\mu),\Psi_t(\nu)\in\mathcal I \mathcal D(\boxtimes)$.
Indeed, on some neighbourhood of zero we have $\eta_{\Psi_t(\mu)}(
\Phi_t^\mu(z))=z,$ and hence $\Sigma_{\Psi_t(\mu)}(z)=\Phi_t^\mu(z)/z=
(z/\eta_\mu(z))^{t-1}.$ By analytic continuation, this equality holds for all 
$z\in\mathbb D.$
In particular, the function
$$u(z)=(t-1)\log\frac{z}{\eta_\mu(z)},\quad z\in\mathbb D$$
is well-defnied, analytic, and $\Sigma_{\Psi_t(\mu)}(z)=\exp(u(z)),$ $z\in
\mathbb D.$
Of course, $|z/\eta_\mu(z)|\geq1$ implies that $\Re u(z)=\exp(|u(z)|)
=\exp((t-1)|\log(z/\eta_\mu(z))|)\geq0.$ We conclude by Theorem \ref{infdivX1}
that $\Psi_t(\mu)\in\mathcal I \mathcal D(\boxtimes).$ Similarily,
$\Psi_t(\nu)\in\mathcal I \mathcal D(\boxtimes).$ 

Assume now that $\mu$ is the uniform distribution on $\mathbb T.$ We define 
$\Psi_t(\mu)=\mu$. Observe that in this case
$\Psi_t(\mu\cup\kern-.7987em\lower-.7ex\hbox{$_\times$}\nu)=\mu=
\Psi_t(\mu)\boxtimes\Psi_t(\nu).$
We have thus proved that $\Psi_t$ is an injective homomorphism. 
To show surjectivity, let $\lambda\in\mathcal I \mathcal D(\boxtimes).$
We may assume without loss of generality that $\lambda$ is not the uniform
distribution on $\mathbb T$. 
By Theorem \ref{infdivX1}, $\Sigma_\lambda(z)=\exp(u(z))$, $z\in\mathbb D$,
for some analytic function $u$ such that $\Re u(z)\geq0.$ 
Let $\Phi(z)=z\Sigma(z),$ $z\in\mathbb D.$ We obviously have $|\Phi(z)|
\geq|z|,$ and $\Phi(0)=0.$ 
Define
  $$\eta(z)=z\left[\frac{z}{\Phi(z)}\right]^{\frac{1}{t-1}},\quad z\in
  \mathbb D.$$
Observe that $\eta(z)=z(\exp(u(z)))^{\frac{1}{1-t}}=z\exp(u(z)/(1-t)),$ and 
thus $\eta(0)=0, $ and
  $$|\eta(z)|\leq|z|(\exp(\Re u(z)))^{\frac{1}{1-t}}\leq|z|,\quad z\in
  \mathbb D.$$
We conclude by Proposition \ref{eta} that there exists a probability measure 
$\mu$ on $\mathbb T$ such that $\eta=\eta_\mu$.
Also, the function $u(z)/(1-t)$ satisfies the condition $\Re u(z)/(1-t)\leq0
,$ and thus, by Proposition \ref{Momentcircle}, is of the form
$u(z)/(1-t)=ib-\int_\mathbb T\frac{e^{it}+z}{e^{it}-z}d\rho(e^{it})$ for some 
positive finite measure $\rho$ on $\mathbb T$, and some number $b\in[0,2\pi).$ 
We conclude by Theorem \ref{Franzinfdiv} that $\mu\in
\mathcal I \mathcal D(\Utimes).$ The computation above shows that
$\Psi_t(\mu)=\lambda.$ 
\end{proof}

\section[]{Partially defined semigroups with respect to free multiplicative 
convolution on the positive half-line}

The results that will be presented in this section are similar to the ones in
the previous sections. As before, we shall start with an inversion theorem, 
this time for maps defined on the slit complex plane.

The covering map $u:\mathbb C\longrightarrow\mathbb C\setminus\{0\}$
defined by $u(z)=-e^{z}$ provides a conformal map of the strip 
${{\mathcal S}}_\pi=\{x+iy\in\mathbb C:|y|<\pi\}$ onto
$\mathbb C\setminus[0,+\infty)$.  Call $v:\mathbb C
\setminus[0,+\infty)\longrightarrow{{\mathcal S}}_\pi$ the inverse of this 
conformal map.
For any function $\Phi:\mathbb C
\setminus[0,+\infty)\longrightarrow\mathbb C\setminus\{0\}$ there exists 
therefore a function
$f:{{\mathcal S}}_\pi\longrightarrow\mathbb C$ such that $\Phi=u\circ f\circ v$. When
$\Phi((-\infty,0))\subset(-\infty,0)$, the
function $f$ is uniquely determined if we require it to be real-valued on
$\mathbb R$ or, equivalently,
$f(\overline z)=\overline{f(z)}, z\in{{\mathcal S}}_\pi$.
Thus, inversion results for functions defined on $\mathcal S_\pi$ can 
be naturally formulated for maps on $\mathbb C\setminus[0,+\infty)$.
We will consider
analytic functions defined in the band ${\mathcal S}_\pi=\{x+iy\in\mathbb
C:|y|<\pi\}\longrightarrow\mathbb C$, which is
of course conformally equivalent to the right half-plane, the equivalence being
realized by
$\varphi:{\mathcal S}_\pi\longrightarrow -i\mathbb C^+$, $\varphi(z)=e^{z/2}$. The map
$\varphi$ sends $+\infty$ to $+\infty$
and $-\infty$ to zero, and from this one can easily deduce what the
Julia-Carath\'eodory derivatives should be
at $\pm\infty$. For a function $g:{\mathcal S}_\pi\longrightarrow
{\mathcal S}_\pi$ they are
$$\lim_{x\uparrow\infty}\frac{\varphi(g(x))}{\varphi(x)}=
\lim_{x\uparrow\infty}e^{(x-g(x))/2}$$
at $+\infty$, and
$$\lim_{x\downarrow-\infty}e^{(g(x)-x)/2}$$
at $-\infty$. In order to use Julia's theorem, we will need to look at the
horodisks associated with
a point $x+i\pi\in\partial{\mathcal S}_\pi$. If $D$ is such a horodisk, then
$\varphi(D)$ is an ordinary disk contained
in the right-half-plane, and tangent to the imaginary axis at the point $e^{
x/2}i$.
The point in (the boundary of) $D$ with the lowest imaginary
part corresponds with the point in $\varphi(D)$ with lowest argument. If
$\varphi(w)$ is this point, then the segment
$[0,\varphi(w)]$ must be tangent to $\varphi(D)$ at $\varphi(w)$, and therefore
it must have the same length as the
other tangent from the origin, namely $[0,e^{x/2}i]$. We deduce that
$|\varphi(w)|=e^{ x/2}$, so that
$\Re w=x$. It is easy to see that any point in ${\mathcal S}_\pi$, with real 
part $x$, is the lowest point in some
horodisk. The consequence we will use is as follows: if $g:{\mathcal S}_\pi
\longrightarrow{\mathcal S}_\pi$ is an analytic function
with Denjoy-Wolff point $x+i\pi$, then $\Im(g(x+iy))\ge y$ for $y\in(0,\pi)$.

\begin{theorem}\label{InversionX2}
Consider an analytic function $f:{{\mathcal S}}_\pi\longrightarrow\mathbb C$ 
with the property that
$f(\overline z)=\overline{f(z)}$ for all $z\in{{\mathcal S}}_\pi$. Assume that there
exists a number $k>1$
such that $\lim_{t\downarrow-\infty}(kt+f(t))=-\infty$ and
$\Im z\le\Im f(z)\le k\Im z$ for all $z\in{{\mathcal S}}_\pi$ with
$\Im z>0$. Then:
\begin{enumerate}
\item[{\rm (1)}] For every $\alpha\in\mathbb {\mathcal S}_\pi$ there exists a unique
$z\in{\mathcal S}_\pi$ such that
$f(z)=\alpha$.
\item[{\rm (2)}] There exists a continuous function $\omega:\overline{{\mathcal
S}_\pi}\longrightarrow\overline{{\mathcal S}_\pi}$
such that
$\omega({{\mathcal S}_\pi})\subset{\mathcal S}_\pi$, $\omega|{\mathcal S}_\pi$ is analytic,
$|\Im \omega(z)|\le |\Im z|,$ $\omega(\overline{z})=\overline{\omega(z)},$
and
$f(\omega(z))=z$ for $z\in{\mathcal S}_\pi$.
\item[{\rm (3)}] For each $\alpha\in\overline{{\mathcal S}_\pi}$, 
$\omega(\alpha)$ is the Denjoy-Wolff point
of the map $g_\alpha:{\mathcal S}_\pi\longrightarrow{\mathcal S}_\pi$ defined by
$$g_\alpha(z)=z-\tau(f(z)-\alpha),\quad z\in{\mathcal S}_\pi,$$where
$\tau=\pi/(k+\pi)$.
\item[{\rm (4)}] The function $\omega$
satisfies $$|\omega(z_1)-\omega(z_2)|\ge\frac {\tau}{2}|z_1-z_2|,\quad
z_1,z_2\in\overline{{\mathcal S}_\pi}.$$
In particular, $\omega$ is one-to-one.
\item[{\rm (5)}] If $\alpha\in\partial{\mathcal S}_\pi$ is such that
$\omega(\alpha)\in{\mathcal S}_\pi$, then $\omega$ can be
continued analytically to a neighborhood of $\alpha$.
\end{enumerate}

\end{theorem}
\begin{proof}

Let us show 
that the functions $g_\alpha$ map ${\mathcal S}_\pi$ into itself. We
clearly have, from the properties of $f,$ that
$$\frac k{k+\pi}\Im z+\frac \pi{k+\pi}\Im\alpha\ge\Im g_\alpha(z)\ge\Im
z-\frac{\pi k}{k+\pi}+\frac \pi{k+\pi}\Im\alpha,\quad
\Im z\ge0,$$ and
$$\frac k{k+\pi}\Im z+\frac \pi{k+\pi}\Im\alpha\le\Im g_\alpha(z)\le\Im
z+\frac{\pi k}{k+\pi}+\frac \pi{k+\pi}\Im\alpha,\quad
\Im z\le0.$$
The left-hand side of these inequalities is a convex combination of $\Im z$ and
$\Im\alpha$, and hence
is in $(-\pi,\pi)$.  For the right-hand sides we observe that
$$\Im z-\frac{\pi k}{k+\pi}+\frac \pi{k+\pi}\Im\alpha\ge-\frac{\pi k}{k+\pi}-\frac \pi{k+\pi}\Im\alpha\geq-\frac{\pi k+\pi^2}{k+\pi}=-\pi$$
for $\Im z\ge0$,
and
$$\Im z+\frac{\pi k}{k+\pi}+\frac \pi{k+\pi}\Im\alpha\le\frac{\pi k}{k+\pi}+\frac \pi{k+\pi}\alpha\leq\frac{\pi k+\pi^2}{k+\pi}=\pi$$
for $\Im z\le0$. 
Observe that 
$$f'(x)=\lim_{y\downarrow0}\frac{f(x+iy)-f(x)}{iy}\ge\lim_{y\downarrow0}\frac{\Im f(x+iy)}y\geq1$$
for all $x\in\mathbb R$. This
implies that
$$\lim_{x\uparrow+\infty}e^{(x-g_\alpha(x))/2
}=\lim_{x\uparrow+\infty}e^{\pi(f(x)-\alpha)/2(\pi+k)}=+\infty,$$
and similarly
$$\lim_{x\downarrow-\infty}e^{(g_\alpha(x)-x)/2
}=+\infty.$$
We conclude that the Denjoy-Wolff point of $g_\alpha$ is not $\pm\infty$, with
the possible exception when $f$ is a real multiple of the identity map.
If $\alpha\in{\mathcal S}_\pi$, then, by the above inequalities, 
$$g({\mathcal S}_\pi)\subseteq\left\{z\in\mathbb C\colon|\Im z|\leq\pi
\frac{k+\Im\alpha}{k+\pi}\right\},$$
so the Denjoy-Wolff point of $g_\alpha$ belongs to $\mathcal S_\pi.$
Denote
by $z_\alpha$ this point.
It is clear  that $f(z_\alpha)=\alpha.$ The function $\omega(\alpha)=z\alpha$ is 
clearly an analytic self-map of $\mathcal S_\pi.$
This proves (1) and parts of (2) and (3).

We wish now to show that the function $\omega$ has a continuous extension to
$\overline{{\mathcal S}_\pi}$.
Fix $\alpha\in\partial{\mathcal S}_\pi$ and, as in the case of the disk, consider 
two cases:
\begin{enumerate}
\item[(a)] There exists a sequence $\{\alpha_n\}_{n\in\mathbb N}\subseteq
{\mathcal S}_\pi$ such that $\lim_{n\to\infty}\alpha_n=\alpha$
and $\lim_{n\to\infty}\omega(z_n)$ exists and belongs to ${\mathcal S}_\pi$;
\item[(b)] For any sequence $\{\alpha_n\}_{n\in\mathbb N}\subseteq{\mathcal S}_\pi$ 
such that $\lim_{n\to\infty}\alpha_n=\alpha$
and $\lim_{n\to\infty}\omega(z_n)$ exists, this limit is in
$\partial{\mathcal S}_\pi$.
\end{enumerate}
If (a) holds and we set $w=\lim_{n\to\infty}\omega(\alpha_n)$, it is imediate
that
$$f(w)=\lim_{n\to\infty}f(\omega(\alpha_n))=\lim_{n\to\infty}\alpha_n=\alpha.$$
This implies that $w$ is a fixed point for $g_\alpha$, and therefore
$w=z_\alpha$.
The map $g_\alpha$ is never a conformal automorphism. Indeed, we have
$$g'_\alpha(x)=1-\frac \pi{k+\pi}f'(x)
\ge 1-\frac k{k+\pi},$$ and this implies that $\pm\infty$ are fixed points for
$g_\alpha$. The only
conformal automorphisms of ${\mathcal S}_\pi$ with this property are translations
$g(z)=z+c$, with
$c\in\mathbb R$, and it is obvious $g_\alpha$ does not have this form. We
conclude that
$$\left|1-\frac \pi{k+\pi}f'(w)\right|=|g'_\alpha(w)|<1,$$
and therefore $f'(w)\ne0$. Thus $f$ is locally invertible at $w$, and the local
inverse continues $\omega$
analytically in a neighborhood of $\alpha$. Thus we have proved (5),
the continuity of $\omega$ at $\alpha$, as well
as assertion (3) in case (a).

Assume now that (b) holds, and that there exist two sequences
$\{\alpha_n\}_{n\in\mathbb N}$ and $\{\beta_n\}_{n\in\mathbb N}$ in 
${\mathcal S}_\pi$ such that $\lim_{n\to\infty}\alpha_n=\lim_{n\to\infty}\beta_n
=\alpha$, the limits $\lim_{n\to\infty}\omega(\alpha_n)$, $\lim_{n\to\infty}
\omega(\beta_n)$ exist, and they are different.
Choose a continuous path $\gamma:[0,1)\to{\mathcal S}_\pi$ passing through all the
points $\alpha_n,\beta_n$,
and such that $\lim_{t\to1}\gamma(t)=\alpha$. It follows from our assumption
that there exists an open interval
$I$ in $\partial{\mathcal S}_\pi$ such that every point $\xi\in I$ is
of the form $\xi=\lim_{n\to\infty}
\omega(\gamma(t_{\xi,n}))$, where $t_{\xi,n}\in[0,1)$ for all $n$. Moreover,
$t_{\xi,n}$ can be chosen
so that $\Re\omega(\gamma(t_{\xi,n}))=\Re\xi$, and thus
$\omega(\gamma(t_{\xi,n}))$ approaches $\xi$
nontangentially. Since
$$\lim_{n\to\infty}f(\omega(\gamma(t_{\xi,n})))=\lim\gamma(t_{\xi,n})=\alpha,$$
we conclude that the nontangential limit of $f$ at $\xi$ is $\alpha$ for almost
every $\xi\in I$. As $f$ takes
values in the band ${\mathcal S}_k=(k/\pi)\mathcal S_\pi$, which is conformally 
equivalent to a disk, 
Theorem \ref{Privalov}
can be used to deduce that $f$ is identically
equal to $\alpha$.
This however is not true because %
$f(0)\in\mathbb R.$ This contradiction shows
that the limit $\omega(\alpha)=
\lim_{z\to\alpha,z\in{\mathcal S}_\pi}\omega(z)$ exists in case (b) as well.

We will show next that assertion (3) is also verified in case (b). Let us note
first that the point $z_\alpha$ cannot
belong to ${\mathcal S}_\pi$ in this case. If it did, we would have
$f(z_\alpha)=\alpha$, and we would deduce as before
that $f$ is invertible in a neighborhood of $z_\alpha$, and the local inverse
continues $\omega$, which would
place us in case (a) instead of (b). Thus $z_\alpha$ is in $\partial{\mathcal S}_\pi$.
Let us assume, for definiteness, that $\Im\alpha=\pi$, in which case we must also
have $\Im z_\alpha = \pi$.
As in the previous arguments,
$g_\alpha$ maps each horodisk at $z_\alpha$ into itself. As noted before the
statement of the theorem,
this implies that, with $x=\Re z_\alpha$ and $y\in(0,\pi)$, $g_\alpha(x+iy)$
is a point close $z_\alpha$ if $y$ is close to $\pi$, and with greater imaginary
part. It follows that
$f(x+iy)\in {\mathcal S}_\pi$, and $\lim_{y\uparrow h}f(x+iy)=\alpha$.
One concludes as before that $\omega(\alpha)=z_\alpha$.

Finally, for $\alpha\in{\mathcal S}_\pi$ we have
$$\left|1-\frac{\pi}{k+\pi}\frac1{\omega'(\alpha)}\right|=
\left|1-\frac{\pi}{k+\pi}f'(z_\alpha)\right|=|g'_\alpha(z_\alpha)|\le1,$$
which implies that
$$\Re\omega'(\alpha)\ge\frac{\pi}{2(k+\pi)}.$$
This implies assertion (4), first in ${\mathcal S}_\pi$, and then by 
continuity in the closure of the domain.\end{proof}
As in the previous sections, we shall apply the above theorem to prove for 
each probability measure $\mu$ on $[0,+\infty)$ the existence of a family 
$\{\mu_t\colon t\geq1\}$ such that $\mu_1=\mu$ and $\mu_{s+t}=\mu_s
\boxtimes\mu_t.$

\begin{theorem}\label{mutX2}
Let $\mu\ne\delta_0$ be a probability measure on $[0,+\infty)$, and let $t\ge1$
be a real
number.
\begin{enumerate}
\item[{\rm(1)}]
There exists a probability measure $\mu_t\ne\delta_0$ on $[0,+\infty)$ such that
$\Sigma_{\mu_t}(z)\newline=\Sigma_\mu(z)^t$ for $z<0$ sufficiently close to zero.
\item[{\rm(2)}]
There exists an analytic self-map of the slit complex plane 
$\omega_t:\mathbb C\setminus[0,+\infty)
\longrightarrow\mathbb C\setminus[0,+\infty)$ such that
$\omega_t((-\infty,0))\subset(-\infty,0)$, $\omega_t(0-)=0$,
$\arg\omega_t(z)\in[\arg(z),\pi)$ for all $z\in\mathbb C^+$, and
$\eta_{\mu_t}(z)=\eta_\mu(\omega_t(z))$ for all $z\in\mathbb C\setminus
[0,+\infty)$.
\item[{\rm(3)}]
The function $\omega_t$ is given by
$$\omega_t(z)=\eta_{\mu_t}(z)\left[\frac z{\eta_{\mu_t}(z)}\right]^{1/t},\quad
z\in\mathbb C\setminus[0,+\infty),$$
where the power is taken to be positive for $z<0$.
\item[{\rm(4)}] The analytic function $\Phi_t:\mathbb C\setminus[0,+\infty)
\longrightarrow\mathbb C\setminus\{0\}$
defined by
$$\Phi_t(z)=z\left[\frac z{\eta_\mu(z)}\right]^{t-1},\quad z\in\mathbb C
\setminus[0,+\infty),$$ satisfies $\Phi_t(\omega_t(z))=z$ for $z\in
\mathbb C\setminus[0,+\infty)$.
\end{enumerate}
\end{theorem}
\begin{proof}
Consider the analytic function $f$ provided by the argument before Theorem
\ref{InversionX2}, such that $\Phi_t=u\circ f\circ v.$ As $\Phi((-\infty,0))
\subseteq(-\infty,0),$ the function $f$ will be uniquely determined by the
requirement that $f(\mathbb R)\subseteq\mathbb R.$ Observe that $\lim_{x\to
-\infty}(2t-1)x+f(x)=\lim_{x\to0}v(x^{(2t-1)}\Phi_t(-x)^{-1})=-\infty$, since
$\lim_{x\uparrow0}x^{2t-2}(\eta_\mu(x)/x)^{t-1}=0,$ as it follows from the 
definition of $\eta_\mu.$ 
Moreover, since $\arg\Phi_t(z)\in(t\arg z-(t-1)\pi,\arg z],$ we conclude that
$\Im z\leq\Im f(z)\leq(2t-1)\Im z.$
Thus, we can apply Theorem \ref{InversionX2} to $f$ to 
obtain a function $\omega$ such that $f(\omega(z))=z$, $z\in{\mathcal S}_\pi.$ 
An argument similar to the one in the proof of Theorem \ref{mutX1} shows that
the function $\omega_t=u\circ\omega\circ v$ will satisfy the required conditions. 
\end{proof}
A regularity result similar to the ones in Theorems \ref{RegPtmut+} and
\ref{regmutX1} holds for semigroups of probability measures on $[0,+\infty).$
Since the proof no different from the proofs of the above mentioned theorems, 
we omit it.
\begin{theorem}\label{regmuX2}
Let $\mu$ be a probability measure on $[0,+\infty)$, and let $t>1$.
\begin{enumerate}
\item[{\rm(1)}]
A point $x\in(0,+\infty)$ satisfies $\eta_{\mu^{\boxtimes t}}(x)=1$ if and only
if $x^{-1/t}$ is an atom of $\mu$ with
mass $\mu(\{x^{-1/t}\})\ge(t-1)/t$. If $\mu(\{x^{-1/t}\})>(t-1)/t$, then $1/x$
is an atom of $\mu^{\boxtimes t}$, and
$$\mu^{\boxtimes t}(\{1/x\})=t\mu(\{ x^{-1/t}\})-(t-1).$$
\item[{\rm(2)}]
The nonatomic part of $\mu^{\boxtimes t}$ is absolutely continuous, and its
density is continuous except
at the (finitely many) points $x$ such that $\eta_{\mu^{\boxtimes t}}(x)=1$.
\item[{\rm(3)}]
The density of $\mu^{\boxtimes t}$ is analytic at all points where it is
different from zero.
\end{enumerate}
\end{theorem}

\chapter[On monotonic infinite divisibility]{On monotonic infinite divisibility
}
In this chapter we improve results of Muraki
related to 
probability measures which infinitely divisible with respect to 
monotonic convolutions. Our main tools will be two theorems of
Ch. Pommerenke and N. Baker, which provide solutions to the 
Abel equation for self-maps of the upper half-plane: given an 
analytic function $f$, we say that $\phi$ satisfies the Abel equation if 
there exists a constant $c\in\mathbb C$ such that $\phi(f(z))=\phi(z)+c$
for all $z$ in the domain of $f$. 

Muraki has given in \cite{Muraki} an analytic method to compute 
additive monotonic convolutions (recall Definition \ref{Monot}) 
of probability measures on the real line in terms of the reciprocals
of their Cauchy transforms:
\begin{theorem}\label{MurakiMonot}
Let $\mu,$ $\nu$ be two probability measures on $\mathbb R$. 
Then 
$$F_{\mu\rhd\nu}(z)=F_\mu(F_\nu(z)),\quad z\in\mathbb C^+.$$
\end{theorem}
(For another proof of this theorem, see \cite{Bercovici?}.)
As definition \ref{infdiv} makes clear, a probability measure $\mu$
will be infinitely divisible with respect to $\rhd$ if and only if 
for any $n\in\mathbb N$ there exists a probability measure $\mu_n$ 
such that $F_\mu(z)=F_{\mu_n}^{\circ n}(z)$ for all $z\in\mathbb C^+.$
(Recall that for a given self-map $f$ of a domain $D$, $f^{\circ n}$
denotes the $n$-fold composition of $f$ with itself: $f^{\circ n}(z)=
\underbrace{f(f(\dots f(z)\dots))}_{n\ {\rm times}}.$)
Muraki has shown that for compactly supported measures  infinite 
divisibility of $\mu$ is equivalent to the embedability of $F_\mu$ 
into a composition semigroup, i.e. a family $\{F_t\colon t\in[0,
+\infty)\}$ of reciprocals of Cauchy transforms of probability measures 
such that $F_0=F_{\delta_0},$ $F_1=F_\mu,$ and $F_{s+t}=F_s\circ F_t$ for
all $s,t\geq0.$

Denote by $H$ the right half-plane, $H=-i\mathbb C^+,$ and consider an 
analytic function $f\colon H
\longrightarrow H$ such that $\lim_{x\to+\infty}f(x)/x=1.$ 
Following Pommerenke, we denote 
  $$f^{\circ n}(1)=z_n=x_n+iy_n,\quad n\in\mathbb N,$$
and 
$$q_n=\frac{z_{n+1}-z_n}{z_{n+1}+\overline{z_n}},\quad n\in\mathbb N.$$
By Pick's inequality, $|q_{n+1}|\leq|q_n|,$ so that $\lim_{n\to\infty}
|q_n|$ exists and is finite. One can find two different types of solutions 
to the Abel equation for $f,$ depending on whether $q_n$ 
tends to zero or not. The following two theorems from \cite{Pommerenke} and
\cite{PommerenkeBaker} make this statement more precise.

\begin{theorem}\label{Pom}
With the above notations, let $k_n\colon H\longrightarrow H,$
$$k_n(z)=\frac{f^{\circ n}(z)-iy_n}{x_n},\quad z\in H,n\in\mathbb N.$$
\begin{enumerate}
\item[{\rm (i)}] The limit $k(z)=\lim_{n\to\infty}k_n(z)$, $z\in H$ 
exists locally uniformly, $k(z)\in H,$ and there exists a number 
$b\in\mathbb R$  such that $k(f(z))=k(z)+ib,$ $z\in H.$
\item[{\rm (ii)}] 
$k(f(1))=f(1)$ if and only if $\lim_{n\to\infty}q_n=0$. In this case 
$k(z)=1$ for all $z\in H.$
\item[{\rm (iii)}] 
Assume that $\lim_{n\to\infty}|q_n|>0$ and fix $d\in(0,1)$. Then $|y_n|\to
\infty,$ $x_n/x_{n+1}\to1$ as $n\to\infty$ and there exists $m\in\mathbb N$
such that $k$ is injective on
$$U_{f,m}^d=\{z\in H\colon\exists n\ge m,
dx_n\leq \Re z\leq x_n/d,\Im z\ {\rm between}\ y_n\ {\rm and}\ y_{n+1}\}.$$ 
Furthermore, $\Re k(x+iy)\to+\infty$ as $x/x_n\to+\infty,$
and $y$ is between $y_n$ and $y_{n+1},$ $n\to\infty.$
\end{enumerate}
\end{theorem}
It is remarkable that the convergence of $\{y_n\}_{n\in\mathbb N}$ is 
eventually monotonic, since $b=\lim_{n\to\infty}(y_{n+1}-y_n)/x_n.$ For 
details we refer to \cite{Pommerenke}, page 443.
\begin{theorem}\label{PomBak}
With the above notations, let 
$$h_n(z)=\frac{f^{\circ n}(z)-z_n}{z_{n+1}-z_n},\quad z\in H,n\in
\mathbb N.$$ Assume that $\lim_{n\to\infty}|q_n|=0.$
\begin{enumerate}
\item[{\rm (i)}] The limit $h(z)=\lim_{n\to\infty}h_n(z)$, $z\in H$
exists locally uniformly in $H,$ and satisfies
$h(f(z))=h(z)+1,$ $z\in H.$
\item[{\rm (ii)}] The function $h$ is injective in
$$G_m=\bigcup_{k=m}^{\infty}\{x+iy\colon x>x_k/2,\ |y-y_k|<x_k\}.$$
\end{enumerate}
\end{theorem}
For proofs, we refer to \cite{Pommerenke} and \cite{PommerenkeBaker}, 
respectively.

One can obviously apply Theorems \ref{Pom} and \ref{PomBak} to self-maps
of the upper half-plane, via the transformation $z\mapsto if(-iz).$
We shall use these two theorems to prove that analytic self-maps of
the upper half-plane that have the Julia-Carath\'eodory derivative at 
infinity equal to 1 can have at most one ``square root'' with respect
to composition.
\begin{proposition}\label{sqrt}
Let $f,g\colon\mathbb C^+\longrightarrow\mathbb C^+$ be two analytic 
self-maps of the upper half-plane such that $\lim_{y\to+\infty}f(iy)/(iy)
=\lim_{y\to+\infty}g(iy)/(iy)=1.$ If $f\circ f=g\circ g,$ then $f=g.$
\end{proposition}
\begin{proof}
Consider the following two self-maps of the right half-plane: 
$\tilde{f}(z)=-if(iz),$ $\tilde{g}(z)=-ig(iz),$ $z\in H.$ Since
$\tilde{f}=\tilde{g}$ if and only if $f=g,$ we can reformulate
the proposition in terms of self-maps of the right half-plane.
So, consider $f$ and $g$ to be self-maps of the left half-plane $H$.
Denote $$f^{\circ n}(1)=z_n=x_n+iy_n,\quad g^{\circ n}(1)=w_n=u_n
+iv_n,\quad n\in\mathbb N,$$
and let
$$q_n=(z_{n+1}-z_n)/(z_{n+1}+\overline{z_n}),\quad n\in\mathbb N.$$
We shall consider separately the cases when $\lim_{n\to\infty}|q_n|>0$ 
and when $\lim_{n\to\infty}|q_n|=0$, corresponding to Theorems \ref{Pom}
and \ref{PomBak}, respectively.

Assume first that $\{q_n\}_{n\in\mathbb N}$ is
bounded away from zero. 
In this case
$$k_{2n}(z)=\frac{f^{\circ 2n}(z)-iy_{2n}}{x_{2n}}=\frac{g^{\circ 2n}(z)-iy_{2n}
}{x_{2n}},\quad z\in H,n\in\mathbb N,$$
converges to a nonconstant function $k$ with the property that 
$k(f(z))=k(g(z))=k(z)+ib$ for some real number $b\neq0.$
Without loss of generality, assume that %
$y_n\to\infty$ as 
$n\to\infty.$ 
As noted after Theorem \ref{Pom}, the sequences $\{y_n\}_{n\in\mathbb N}$
and $\{v_n\}_{n\in\mathbb N}$ tend either to $+\infty,$ or to $-\infty$. Since
$y_{2n}=v_{2n},$ the limit must be the same; we assume for definiteness
that it is $+\infty.$
By Theorem \ref{Pom}, the sequences 
$\{x_{n+1}/x_n\}_{n\in\mathbb N}$ and $\{u_{n+1}/u_n\}_{n\in\mathbb N}$
tend to one as $n\to\infty,$ so
there exists $m>0$ such that both these sequences have all their terms 
at distance less than $1/4$ from $1$ for $n\geq m.$
We shall choose $m_1\in\mathbb N$ large enough so that $y_n<y_{n+1}<\dots,$
$v_n<v_{n+1}<\dots,$ and $u_{n}\in((1/2)u_{n-1},2u_{n-1}),$ 
$x_{n}\in((1/2)x_{n-1},2x_{n-1})$ for $n\ge m_1.$

Observe that, since $x_{2n}=u_{2n}$ and $y_{2n}=v_{2n}$, for $n\ge m_1,$ we have
\begin{eqnarray}
\lefteqn{
\left\{x+iy\colon\frac12x_{2n}\leq x\leq2x_{2n},\ y_{2n}\leq y\leq y_{2n+1}\right\}\cap}\nonumber\\
& & \left\{x+iy\colon\frac12u_{2n}\leq x\leq2u_{2n},\ v_{2n}\leq y\leq v_{2n+1}\right\}\neq
\varnothing.\nonumber
\end{eqnarray}
Moreover, for $n\ge m_1,$
either $z_{2n+1}\in\{x+iy\colon(1/2)u_{2n}<x<2u_{2n},
\ v_{2n}<y< v_{2n+1}\},$ or
$w_{2n+1}\in\{x+iy\colon(1/2)x_{2n}<x<2x_{2n},
\ y_{2n}< y< y_{2n+1}\}.$
Choose now $n>m_1$ to be large enough so that $k$ is injective on 
$U_{f,n}^{1/2}$. Since $k(w_{2n+1})=k(z_{2n+1})$ and both $w_{2n+1}$ and 
$z_{2n+1}$ belong to $U_{f,n}^{1/2}$, we conclude that $w_{2n+1}=z_{2n+1}$.
This is true for all $n$ large enough. %
Consider a neighbourhood $V$ of $z_{2n}$. It is clear that the 
nonempty open set $f(V)\cap g(V)$ will be included in 
$\{x+iy\colon(1/2)u_{2n}< x<2u_{2n}, \ v_{2n}< y< v_{2n+1}\}\cup
\{x+iy\colon(1/2)u_{2n+1}< x<2u_{2n+1}, \ v_{2n+1}< y< v_{2n+2}\}$,
if $V$ is small enough. But this will imply that for any $z\in V$ we have $k
(f(z))=k(z)+ib=k(g(z)).$ 
Since $\{x+iy\colon(1/2)u_{2n}< x<2u_{2n}, \ v_{2n}< y< v_{2n+1}\}\cup
\{x+iy\colon(1/2)u_{2n+1}<x<2u_{2n+1}, \ v_{2n+1}< y< v_{2n+2}\}\subset
U_{f,n}^{1/2}$, we conclude that $f(z)=g(z)$ for all $z\in V,$ and hence
$f=g.$

Consider now the second alternative, that is $\lim_{n\to\infty}q_n=0.$
We claim that 
$$\lim_{n\to\infty}\frac{f^{\circ n}(z)-z_n}{z_{n+1}-z_n}=\lim_{n\to\infty}
\frac{g^{\circ n}(z)-w_n}{w_{n+1}-w_n},\quad z\in H.$$
Indeed, 
$$h_n(z)=\frac{f^{\circ n}(f(f(1)))-z_n}{z_{n+1}-z_n}\cdot
\frac{f^{\circ n}(z)-z_n}{z_{n+2}-z_n}
=h_n(f(f(1)))\cdot\frac{f^{\circ n}(z)-z_n}{z_{n+2}-z_n};$$
passing to limit in the above equality gives
$$h(z)=h(f(f(1)))\cdot\lim_{n\to\infty}\frac{f^{\circ n}(z)-z_n}{z_{n+2}-z_n}
=2\lim_{n\to\infty}\frac{f^{\circ n}(z)-z_n}{z_{n+2}-z_n}.$$
The claim follows now from the fact that 
$$\frac{f^{\circ 2n}(z)-z_{2n}}{z_{2n+2}-z_{2n}}=\frac{g^{\circ 2n}(z)-z_{2n}
}{z_{2n+2}-z_{2n}},\quad z\in H,n\in\mathbb N.$$

With the notation from Theorem \ref{PomBak} (ii), we know from 
\cite{PommerenkeBaker} that $f$ is injective on $G_m=G_m^f$ for $m$ large 
enough.
Consider a disk $D\subset H$ such that $1,f(1),g(1),f(f(1))=g(g(1))\in
D.$ Since $k_n(z)=(f^{\circ n}(z)-iy_n)/x_n\to 1$ uniformly on $\overline{D},
$ for $N$ large enough, $B_f=\cup_{n=N}^\infty f^{\circ n}(D)$ is a connected
subset of $G_m.$ Since $f(B_f)\subseteq B_f$, the functions
$$h_n^*(z)=\frac{f^{\circ(n-N)}(z)-z_n}{z_{n+1}-z_n},\quad z\in H,\ n
\in\mathbb N,\ n>N,$$
are injective on $B_f$. It is easy to observe that $h_n^*(z_N)=0$ and $h_n^*(
z_{N+1})=1.$ Thus, the sequence $\{h_n^*\}_{n>N}$ converges to a function $h^*
$, defined on $H,$ and injective on $B_f$. (For details, see 
\cite{PommerenkeBaker}, page 256.) 
We have $$h^*(f^{\circ N}(z))=\lim_{n\to\infty}h_n^*(f^{\circ N}(z))=
\lim_{n\to\infty}h_n(z)=h(z),\quad z\in H.$$
The same argument used for $h$ shows that if we replace in the above $f$ by $g$
we obtain a function defined on $B_g=\cup_{n=N}^\infty g^{\circ n}(D),$ which
coincides with $h^*$ on $B_f\cap B_g$ (it is trivial to 
observe that $B_f\cap B_g\neq\varnothing$). Moreover, by an argument 
identical to the one used in the previous case (the case when $q_n$ was bounded
away from 0), there exists 
an open set $V\subseteq\mathbb C^+$ such that $f^{\circ2N+1}(z),g^{\circ2N+1}(z
)\in B_f\cap B_g$ for all $z\in V.$ Indeed, $f^{\circ 2N+1}(i),g^{\circ 2N+1}(i
)\in f^{\circ 2N}(D)=g^{\circ 2N}(D),$ which is an open set. But then 
$h^*(f^{\circ2N+1}(z))=
h^*(g^{\circ2N+1}(z)),$ so that, by injectivity of $h^*,$ $f^{\circ2N+1}(z)=
g^{\circ2N+1}(z),$ $z\in V,$ and thus, $f=g.$
\end{proof}

In his paper \cite{Cowen81}, Carl Cowen has a different approach towards 
solving Abel's equation. This approach will be useful for proving
that every $\rhd$-infinitely divisible probability measure $\mu$ on 
$\mathbb R$ belongs to a unique semigroup $\{\mu_t\colon t\geq0\}$
with respect to monotonic additive convolution. In the following, we shall
introduce the notions and results relevant for our proof. For details, 
we refer to \cite{Cowen81}.

\begin{definition}\label{fundset}
Given a domain $\Delta\subseteq\mathbb C$ and an analytic map 
$\psi\colon\Delta\longrightarrow\Delta$, we say that $V$ is a fundamental
set for $\psi$ on $\Delta$ if $V$ is an open, connected, simply connected 
subset of $\Delta$ such that $\psi(V)\subseteq V$ and for each compact 
set $K$ in $\Delta$, there exists $N\in\mathbb N$ so that $\psi^{\circ N}
(K)\subseteq V.$
\end{definition}
The main result of \cite{Cowen81} is the following theorem:

\begin{theorem}\label{Cowen}
Let $\varphi\colon\mathbb D\longrightarrow\mathbb D$ be analytic, not
constant and not a conformal automorphism of $\mathbb D,$ and let $a
\in\overline{\mathbb D}$ be the Denjoy-Wolff point of $\varphi.$ 
Assume that $\varphi'(a)\neq0.$ Then there exists a fundamental set
$V$ for $\varphi$ on $\mathbb D$, a domain $\Omega$, either the complex
plane or the unit disk, a linear fractional transformation $\Phi$
mapping $\Omega$ onto itself, and an analytic mapping $\sigma\colon
\mathbb D\longrightarrow\Omega$ such that $\varphi$ and $\sigma$ are
injective on $V$, $\sigma(V)$ is a fundamental set for $\Phi$ on $
\Omega$, and $\Phi\circ\sigma=\sigma\circ\varphi.$ Moreover,
$\Phi$ is unique up to conjugation by a linear fractional transformation
mapping $\Omega$ onto itself, and $\Phi$ and $\sigma$ depend only on $
\varphi,$ not on the particular fundamental set $V$.
\end{theorem}
It is shown that, depending on $\varphi,$ the linear fractional 
transformation $\Phi$ is of one of the following types:
\begin{enumerate}
\item[(1)] $\Phi(z)=sz$ for some $s\in\mathbb D,$ $z\in\mathbb C,$
\item[(2)]  $\Phi(z)=z+1,$ $z\in\mathbb C,$
\item[(3)] $\Phi(z)=\frac{(1+s)z+(1-s)}{(1-s)z+(1+s)}$ for some $s\in(0,1),$
$z\in\mathbb D,$
\item[(4)] $\Phi(z)=\frac{(1\pm 2i)z-1}{z-1\mp2i}$, $z\in\mathbb D.$
\end{enumerate}

Observe that any linear fractional transformation $\Phi$ is contained
in a real analytic group with respect to composition: there exists 
$H\colon\Omega\times\mathbb R\longrightarrow\Omega$ such 
that $H(H(z,t),s)=H(z,t+s)$ and $H(z,n)=\Phi^{\circ n}(z)$ for all
$z\in\Omega$ and $n\in\mathbb N.$ Denote by $G$ the infinitesimal 
generator of the semigroup $H$. 
Consider now a function $\varphi$ as in the above theorem.
Fix a
fundamental set $V$ of $\varphi$ on $\mathbb D$, and let $\nu(z)=
\min\{n\in\mathbb N\colon\varphi^{\circ n}(z)\in V\}.$
Define
$$\tau(z)=\inf\{t\ge0\colon\nu(z)\leq t, H(\sigma(z),s)\in\sigma(V)\ 
{\rm for\ all}\ s>t\}.$$

\begin{theorem}\label{Cowensemigroup}
Let $\varphi$ satisfy the hypotheses of Theorem \ref{Cowen}, and
let $\tau$ be defined as above. There exists a function $\varphi(z,t)$
defined for $z\in\mathbb D$ and $t>\tau(z)$, complex analytic 
in the first argument, real analytic in the second, such that
$$\varphi(\varphi(z,t),s)=\varphi(z,t+s),\quad t>\tau(z),s>0,$$
and such that $\varphi(z,n)=\varphi^{\circ n}(z)$ for all $n\in
\mathbb N$, $n>\tau(z).$
The function $\varphi(z,t)$ is defined by 
$$\varphi(z,t)=(\sigma|_V)^{-1}(H(\sigma(z),t)),\quad z\in\mathbb D,t>\tau(z).
$$
Moreover, the function $g$, given by $g(z)=G(\sigma(z))(\sigma'(z))^{-1},$ 
$z\in\mathbb D,$ is meromorphic on $\mathbb D,$ holomorphic on $V$, and 
agrees with the infinitesimal generator of the semigroup $\varphi(z,t)$ on 
the set $\{z\in\mathbb D\colon\tau(z)=0\}.$
\end{theorem}
Assume that $V$ satisfies the conditions of Theorem \ref{Cowen} for 
$\varphi$ and consider a compact set with nonempty interior $K\subset V$.
Then, by Theorem \ref{Cowen}, there exists $n_0\in\mathbb N$ so that 
$\Phi^{\circ n}(H(\sigma(z),t))\in\sigma(V)$ for all $0\leq t\leq 1,$ $n\ge 
n_0,$ $z\in K.$ But by definition $\Phi^{\circ n}(H(\sigma(z),t))=
H(\sigma(z),t+n).$ Let $w=\Phi^{\circ n}(\sigma(z))=H(\sigma(z),n)\in
\sigma(V),$ so that there exists $v\in V$ such that $\sigma(v)=z.$
We claim that $\tau(v)=0.$ Indeed, since $v\in V,$ we have $\nu(v)=0,$
and 
$$H(\sigma(v),t)=H(H(\sigma(z),n),t)=H(\sigma(z),t+n)\in\sigma(V)\quad
{\rm for\ all}\ t\ge0.$$
This proves our claim. We conclude that the interior of the set
$V\cap\{z\in\mathbb D\colon\tau(z)=0\}$ is nonempty.

All the above results can easily be reformulated for self-maps of 
the upper half-plane.

Recall that the probability measure $\mu$ on $\mathbb R$ is said to 
be $\rhd$-infinitely divisible if for any $n\in\mathbb N$ there exists 
a probability measure $\mu_{1/n}$ such that $\mu=\underbrace{\mu_{1/n}
\rhd\cdots\rhd\mu_{1/n}}_{n\ \rm times}.$ The measure $\mu_{1/n}$ will
be called an $n^{\rm th}$ root of $\mu.$ The following proposition 
generalizes the results of \cite{Muraki}, Proposition 5.4, to measures
without compact support.

\begin{proposition}\label{compsemigr}
Let $\mu$ be a $\rhd$-infinitely divisible probability measure on $
\mathbb R.$ Then 
\begin{enumerate}
\item[{\rm(1)}] The family of probability measures 
$\{\mu_{1/2^n}\colon n\in\mathbb N\}$
can be uniquely extended to a weak*-continuous semigroup
$\{\mu_t\colon t\ge0\}$ with respect to monotonic additive 
convolution.
\item[{\rm(2)}] The $n^{\rm th}$ root of $\mu$ is unique for all 
$n\in\mathbb N.$
\end{enumerate}
\end{proposition}
\begin{proof}
We shall first show that there exists an appropriate fundamental set for 
our functions $F_{\mu_{2^j}}$, $j\in\mathbb Z.$
Let us observe that if $V_{{1/n}}$ is a fundamental
set for $F_{\mu_{1/n}}$ on $\mathbb C^+,$ then $V_{{1/n}}$
is a fundamental set for $F_{\mu_{1/n}}^{\circ p}=
F_{\mu_{p/n}}$ for all $p\in\mathbb N.$ Moreover, if $V_{{
1/n}}$ satisfies the conditions of Theorem \ref{Cowen} for
$F_{\mu_{1/n}}$, then it satisfies the same conditions for 
$F_{\mu_{p/n}}$ for all $p\in\mathbb N.$
Indeed, it is clear that $F_{\mu_{p/n}}(V_{{1/n}})\subseteq
V_{{1/n}},$ and thus, by the Denjoy-Wolff Theorem (Theorem
\ref{DenjoyWolff}), for any compact set $K\subset\mathbb C^+$ there
exists $N\in\mathbb N$, $N>0,$ such that $F_{\mu_{p/n}}^{\circ N}(K)
\subset V_{{1/n}}$. Assume now that $z_1,z_2\in V_{{1/n}}$
are such that $F_{\mu_{p/n}}(z_1)=F_{\mu_{p/n}}(z_2).$ Then 
$F_{\mu_{1/n}}(F_{\mu_{(p-1)/n}}(z_1))=F_{\mu_{1/n}}(F_{\mu_{(p-1)/n}}(z_2));$
 since $F_{\mu_{(p-1)/n}}(z_1),F_{\mu_{(p-1)/n}}(z_2)\in V_{{1/n}},$
we have that $F_{\mu_{(p-1)/n}}(z_1)=F_{\mu_{(p-1)/n}}(z_2).$ Iterating this
procedure, we obtain that $z_1=z_2.$

On the other hand, observe that $F_{\mu_{1/n}}$ is injective on $V_{1}$
for any fundamental set $V_{1}$ of $F_\mu$ on $\mathbb C^+.$ Indeed, 
if $z_1,z_2\in V_\mu$ and $F_{\mu_{1/n}}(z_1)=F_{\mu_{1/n}}(z_2)$ then
$F_\mu(z_1)=F_\mu(z_2),$ and this implies $z_1=z_2.$

Fix now $n\in\mathbb N.$ Without loss of generality, we may 
assume that $V_{{1/2^j}}\subseteq V_{1/2^{j-1}},$ $1\le j\le n.$
Proposition \ref{sqrt} guarantees the uniqueness of the measures $\mu_{1/2^j},
$ $j\in\mathbb Z.$
By Theorems \ref{Pom} and \ref{PomBak}, we know that there exist $\sigma$
and $c\in\mathbb C\setminus\{0\}$ such that $\sigma(F_{\mu}(z))=\sigma(z)+c.$
For fixed $c$, Theorem \ref{Cowen} assures us of the uniqueness of $\sigma.$ 
If we denote $H(w,t)=w+tc,$ $t\in\mathbb R,$ $w\in\mathbb C,$ 
Theorem \ref{Cowensemigroup} provides the function $F_{1/2^j}(z)=
(\sigma|_V)^{-1}(H(\sigma(z),t))$ for all $z\in\mathbb C^+,$ $t>\tau(z).$
We shall find a domain on which these functions coincide with $F_{\mu_{2^{-j}}}
,$ $1\le j\le n.$


Assume $z\in V_{{1/2^n}}\cap\{z\in\mathbb C^+\colon\tau(z)=0\}.$
Then $F_{2^j}(z)=(\sigma|_V)^{-1}(H(\sigma(z),2^j))$ for all
$j\ge-n.$ Since $\sigma(F_{2^j}(z))=\sigma(F_{\mu_{1/2^j}}(z))$ and
$F_{2^j}(z),F_{\mu_{1/2^j}}(z)\in V_{\mu_{1/2^n}},$ we have 
$F_{2^j}(z)=F_{\mu_{1/2^j}}(z),$ for $z$ in an open subset of the 
upper half-plane, $j\ge-n.$ This is true for arbitrary $n\in\mathbb N,$
and thus we have $F_{2^m}=F_{\mu_{2^m}}$ for all $m\in\mathbb Z.$
By Theorem \ref{Cowensemigroup}, $F_{2^j}(z)\to z$ as $j$ tends to $-\infty$
uniformly on compact subsets of $\{z\in\mathbb C^+\colon\tau(z)=0\}.$
We conclude that $\lim_{j\to-\infty}F_{\mu_{2^j}}(z)= z$ uniformly on
compact subsets of $\mathbb C^+.$ 
To obtain $\mu_t$ for some $t\in[0,+\infty),$ consider a sequence 
$q_n$ of dyadic numbers converging to $t$. By Proposition 
\ref{Cauchyalt}, the analytic function
$F_t=\lim_{n\to\infty}F_{\mu_{q_n}}$ provides the required measure 
$\mu_t$ by $F_{\mu_t}=F_t$ (the independence of the limit $F_t$ from the 
chosen sequence $q_n$ follows from the fact that $F_{2^j}(z)\to z$ as 
$j$ tends to $-\infty$).  
This proves (1). Part 2 is an immediate consequence of the uniqueness
of the semigroup from part 1.
\end{proof}

Similar results hold for multiplicative monotone convolution of probability
measures supported on $[0,+\infty)$ (recall Definition \ref{Monot}). However, 
in this case, the analogue of Proposition \ref{compsemigr} is an immediate 
consequence of results of Cowen and Gorya{\u\i}nov. 
Let us first state the analogue of Theorem \ref{MurakiMonot}:
\begin{theorem}\label{BercoviciMonot}
Consider two probability meaures $\mu,\nu$ on $[0,+\infty).$ Then we have
$$\eta_{\mu\circlearrowright\nu}(z)=\eta_\mu(\eta_\nu(z)),\quad z\in\mathbb C
\setminus[0,+\infty).$$
\end{theorem}
For the proof of this theorem, we refer to \cite{BercoviciMonot}
Define $\theta\colon\mathbb C\setminus(-\infty,0)\longrightarrow H$ by
$\theta(z)=\sqrt{|z|}e^{\frac{i\arg z}{2}},$ where $\arg z$ is defined as in 
the beginning of Chapter 2, Section 1.2 on $\mathbb C^+$, while
$\arg z=-\arg\overline{z}$ if $z\in-\mathbb C^+,$ and $\arg z=0$ if
$z\in(0,+\infty).$ Then the function $$\mathcal V\colon
\mathbb C\setminus(-\infty,0)\longrightarrow \mathbb D,\quad
\mathcal V(z)=\frac{\theta(-z)-1}{\theta(-z)+1},$$ 
is a conformal automorphism such that $\mathcal V((-\infty,0))=(-1,1).$
Thus, for any probability measure $\mu$ on the positive half-line, the
function $\phi_\mu(z)=\mathcal V(\eta_\mu(\mathcal V^{-1}(z))),$ $z\in
\mathbb D,$ is a self-map of the unit disk that maps the interval 
$(-1,1)$ into itself. Moreover, it is trivial to see that $\eta_\mu$ 
belongs to a composition semigroup $\{\eta_{\mu_t}\colon t\ge 0\}$ if and only 
if $\phi_\mu$ belongs to a composition semigroup $\{\phi_{\mu_t}\colon t\ge0\},
$ and the correspondence is given by $\phi_{\mu_t}(z)=\mathcal V(\eta_{\mu_t}
(\mathcal V^{-1}(z))),$ $z\in\mathbb D.$ 
In \cite{Goryainov}, Gorya{\u\i}nov shows that any analytic self-map 
$\phi$ of the unit disk preserving the interval $(-1,1)$ which is the 
$n$-iteration of an analytic self-map $\phi_n$ of $\mathbb D$ satisfying
$\phi_n((-1,1))\subseteq(-1,1)$ can be embedded in a one-parameter semigroup
with respect to composition (this fact can be also derived from Theorem 5.1
and its corollary from \cite{Cowen81}). This, of course, provides the 
equivalent of Proposition \ref{compsemigr} for multiplicative monotone 
convolution on $[0,+\infty).$

\end{document}